\documentclass[11pt,reqno]{amsart}
\usepackage[T1]{fontenc}
\IfFileExists{newtxtext.sty}{\usepackage{newtxtext,newtxmath}}{\usepackage{mathptmx}}

\usepackage{mathrsfs,amssymb}
\usepackage{microtype,needspace}
\usepackage{amsmath,amsthm,mathtools}
\usepackage{booktabs,array,longtable,multicol}
\usepackage[margin=1.08in,headheight=14pt]{geometry}
\usepackage[shortlabels]{enumitem}
\usepackage{xurl}
\usepackage[hidelinks]{hyperref}
\usepackage[numbers,sort&compress]{natbib}
\allowdisplaybreaks[2]
\numberwithin{equation}{section}
\newtheorem{theorem}{Theorem}[section]
\newtheorem{proposition}[theorem]{Proposition}
\newtheorem{lemma}[theorem]{Lemma}
\newtheorem{corollary}[theorem]{Corollary}

\theoremstyle{definition}

\theoremstyle{remark}

\newcommand{\RR}{\mathbb R}
\newcommand{\CC}{\mathbb C}
\newcommand{\ZZ}{\mathbb Z}

\newcommand{\calA}{\mathcal A}
\newcommand{\calM}{\mathcal M}

\newcommand{\calB}{\mathcal B}
\newcommand{\tr}{\operatorname{tr}}
\newcommand{\Tr}{\operatorname{Tr}}
\newcommand{\sa}{\mathrm{sa}}
\newcommand{\id}{\mathbf 1}
\newcommand{\cyc}{\mathrm{cyc}}
\newcommand{\coeff}{\mathrm{coeff}}

\newcommand{\sym}{\mathrm{sym}}

\newcommand{\qnorm}[1]{\lVert#1\rVert_2}
\newcommand{\qop}[1]{\lVert#1\rVert}
\newcommand{\qac}[2]{\{#1,#2\}}
\newcommand{\enorm}[1]{\lVert#1\rVert}
\newcommand{\fiveF}{\mathbb F}
\newcommand{\fiveZ}{\mathbb Z}

\newcommand{\fiveEL}{\operatorname{EL}}
\newcommand{\fiveGL}{\operatorname{GL}}
\newcommand{\fivespec}{\operatorname{spec}}
\newcommand{\fivesgn}{\operatorname{sign}}
\title[Degree and arity of tracial separation]{Degree and arity of tracial separation for self-adjoint contractions}
\author{Jie Wang}
\address{State Key Laboratory of Mathematical Sciences, Academy of Mathematics and Systems Science, Chinese Academy of Sciences, Beijing, China}
\email{wangjie212@amss.ac.cn}
\subjclass[2020]{46L10, 13J30, 16S34, 15A18, 90C22}
\keywords{Connes' embedding problem, trace-positive noncommutative polynomial, tracial moment problem, sums of Hermitian squares, quadratic module, matrix contractions, Horn inequalities}
\begin{document}
\begin{abstract}
Motivated by the recent explicit counterexample to Connes' embedding conjecture, we study the degree and number of variables needed to separate normalized matrix traces on self-adjoint contractions from traces on finite von Neumann algebras. We prove that the two trace infima coincide for every Hermitian polynomial of degree at most two in at most three variables, with a dimension-explicit approximation bound. We give three explicit separators: an integer-coefficient polynomial of degree six in ten variables with 392 nonzero words, a quadratic polynomial in 214 variables, and an integer-coefficient polynomial in two variables of ordinary and cyclic degree 24. Each has normalized trace at least \(3/4\) on every matrix contraction tuple and trace \(-1\) at a finite tracial witness. We provide the source systems, transfer estimates, and exact coefficient recipes for all three constructions.
\end{abstract}

\maketitle
\section{Introduction}
\subsection{The trace separation problem}
Let \(\CC\langle X_1,\ldots,X_n\rangle\) be the free \(*\)-algebra on self-adjoint indeterminates. For a Hermitian polynomial \(f=f^*\), we compare normalized traces on self-adjoint matrix tuples \(A=(A_1,\ldots,A_n)\) of common size with traces on finite von Neumann algebras. Requiring contractions bounds the domains independently of matrix size. Set
\begin{equation}\label{eq:two-infima}
\begin{aligned}
 m_{\mathrm{mat}}(f)&=\inf_{N\ge 1}\ \min_{\substack{A_i=A_i^*\in M_N(\CC)\\ \|A_i\|\le 1}}\tr_N f(A),\\
 m_{\mathrm{tr}}(f)&=\inf_{\substack{(M,\tau)\text{ finite tracial}\\ a_i=a_i^*,\ \|a_i\|\le 1}}\tau f(a).
\end{aligned}
\end{equation}
Here \(\tr_N=N^{-1}\Tr\), and finite tracial algebras have a normalized faithful normal trace. By including matrices among finite tracial algebras, one always has \(m_{\mathrm{tr}}(f)\le m_{\mathrm{mat}}(f)\). We call \(f\) a \emph{tracial separator} when the inequality is strict. The constructions below have the normalization
\begin{equation}\label{eq:normalized-gap}
\tr_N f(A)\ge \tfrac34\quad\text{for all matrix contraction tuples,}
\qquad\tau f(a)=-1\quad\text{for one tracial contraction tuple.}
\end{equation}
Positivity at each matrix size need not give a uniform positive lower bound.

The operator-algebraic background is Connes' embedding problem for finite factors \cite{Connes1976}. It is now known to have a negative answer: the \(\mathrm{MIP}^*=\mathrm{RE}\) theorem implies a counterexample to the conjecture \cite{JiEtAl2020}. On the polynomial side, unrestricted matrix positivity is characterized by sums of squares \cite{Helton2002}, whereas trace positivity leads to a distinct tracial moment problem \cite{BurgdorfKlep2012}. Klep and Schweighofer connected this trace-positive setting to Hermitian squares, contraction constraints and commutators \cite{KlepSchweighofer2008}; the real-algebraic equivalence requires the correction in \cite{BurgdorfDykemaKlepSchweighofer2014}.

Wang and Zhi gave explicit counterexamples to Connes' embedding conjecture \cite{WangZhi2026}. Their main construction is an integer-coefficient Hermitian polynomial in six variables of degree 72 with 33 monomials. It has the form \(f=-1+\omega+\omega^*+M\sum_{j=1}^{15}(2-u_j-u_j^*)\), where \(\omega,u_1,\ldots,u_{15}\) are specified words and \(M\) is an explicit positive integer. The relation penalties are nonnegative in normalized matrix trace on contraction tuples and vanish at the specified group witness. Consequently, \(f\) has normalized trace at least \(3/4\) on every matrix contraction tuple, while the witness has trace \(-1\) \cite[Theorem~1.1]{WangZhi2026}. This gap also rules out a contraction-module certificate for \(f+\varepsilon\) modulo commutators whenever \(0\leq\varepsilon<1\) \cite[Corollary~5.12]{WangZhi2026}. The same paper constructs a 37-variable quartic with Gaussian-integer coefficients and 387 monomials, with the same trace gap but without norm restrictions on the matrix inputs, and an integer-coefficient two-variable separator of degree at most 48 on contractions \cite[Appendices~A--B]{WangZhi2026}. Thus degree four is minimal in the unrestricted setting, while two variables suffice and are necessary in the contraction setting. Their group obstruction and quantitative normalization use results of Kun--Thom, Thom, and Alekseev--Liu--Thom \cite{KunThom2026,Thom2026,AlekseevLiuThom2026}.

The constructions give further degree--arity tradeoffs on contractions. The degree-six example lowers the degree of the cited six-variable construction using four additional variables. The quadratic construction uses a characteristic-two presentation, whose basic root elements are symmetries, and three-corner projection geometry to convert relation defects into quadratic trace energies. Its characteristic-two source requires a separate proof from the characteristic-five source of \cite{WangZhi2026}. The two-variable construction lowers the cited degree bound from 48 to 24.

\subsection{Four complementary results}
Table~\ref{tab:main} summarizes the constructions and the rigidity theorem. Each separator has the normalization \eqref{eq:normalized-gap}. Degree means ordinary free-word degree; the last row also has cyclic degree 24.
\begin{table}[ht]
\centering
\caption{Dimension-independent tracial separation and the quadratic obstruction.}
\label{tab:main}
\small
\begin{tabular}{@{}p{0.29\linewidth}rrp{0.43\linewidth}@{}}
\toprule
Result & Variables & Degree & Further information\\
\midrule
Quadratic rigidity & \(n\le 3\) & \(\le 2\) & \(m_{\mathrm{mat}}=m_{\mathrm{tr}}\)\\
Quadratic construction & 214 & 2 & 2,416 ordinary words; symmetry witness\\
Short-degree construction & 10 & 6 & 392 ordinary words; positive contraction witness\\
Two-variable construction & 2 & 24 & Exact cyclic degree; strict contraction witness\\
\bottomrule
\end{tabular}
\end{table}

Write \(n_2\) for the smallest number of self-adjoint variables admitting a degree-at-most-two tracial separator, and \(d_2\) for the smallest degree of a two-variable separator. Combining the constructions and quadratic rigidity with the cubic obstruction of Corollary~\ref{cor:cubic} gives
\begin{equation}\label{eq:arity-bounds}
 4\le n_2\le 214,\qquad 4\le d_2\le 24.
\end{equation}
Scalar spectral integration excludes one-variable separators at every degree. Thus two is the minimum arity. In two variables, the cyclically sorted polynomial result of Klep--Schweighofer \cite[Proposition~4.2]{KlepSchweighofer2008} excludes separators of degree at most three.

\subsection{Method and proof scope}
The constructions convert group-theoretic trace obstructions into low-degree polynomials. For degree six, a weighted graph encodes cyclic products as closed walks; calibration squares and consistency penalties control the decoding on arbitrary matrices. For degree two, projection geometry converts products, relations and spectral constraints into quadratic energies, with a forest of rectangular calibrations. For two variables, a moment encoder realizes auxiliary operators in a small corner, and algebraic corrections control the selector defect. Each construction requires both a uniform finite-matrix lower bound and an exact negative tracial witness.

Section~\ref{sec:notation} gives the conventions and elementary obstructions; Section~\ref{sec:rigidity} proves quadratic rigidity and its quantitative form. After the transfer estimates in Section~\ref{sec:transfer}, Sections~\ref{sec:degree6}--\ref{sec:degree24} give the three constructions. Section~\ref{sec:outlook} discusses their consequences and limits.

\section{Tracial positivity and elementary reductions}\label{sec:notation}
Let \(\calA_n^{\RR}=\RR\langle X_1,\ldots,X_n\rangle\) and \(\calA_n=\CC\langle X_1,\ldots,X_n\rangle\). Involution fixes the variables, conjugates coefficients and reverses words. Write \([\calA_n,\calA_n]_{\rm lin}\) for the \(\CC\)-linear span of additive commutators. Two Hermitian polynomials are \emph{cyclically equivalent} when their difference lies in this span; their tracial evaluations then agree. For \(f\in\calA_n\), let \(\deg_{\cyc}f\) be the highest degree of a nonzero coefficient sum on a cyclic word class (if such a class exists).

The \emph{contraction quadratic module} is
\begin{equation}\label{eq:module}
 \calM_n=\left\{\sum_{\alpha}p_\alpha^*p_\alpha+
 \sum_{i=1}^n\sum_\beta q_{i\beta}^*(1-X_i^2)q_{i\beta}:p_\alpha,q_{i\beta}\in\calA_n\right\}.
\end{equation}
All sums are finite. The real version uses \(\calA_n^{\RR}\) in place of \(\calA_n\). The module is positive under every representation sending generators to self-adjoint contractions, not only finite-dimensional ones. The relevant property is \(\calM_n+[\calA_n,\calA_n]_{\rm lin}\), not the unrestricted cone of Hermitian squares alone.

\begin{proposition}\label{prop:moduleobstruction}
Suppose \(f=f^*\in\calA_n\) and there exists a self-adjoint contraction tuple \(a\) in a finite tracial algebra \((M,\tau)\) with \(\tau f(a)=-\eta<0\). Then, for \(0\le\varepsilon<\eta\),
\[
 f+\varepsilon\notin \calM_n+[\calA_n,\calA_n]_{\rm lin}.
\]
The same holds over the real free algebra when \(f\) has real coefficients.
\end{proposition}
\begin{proof}
Every Hermitian square has nonnegative tracial evaluation. Since \(\|a_i\|\le1\), the operators \(1-a_i^2\) are positive, so the same is true of each localizing term in \eqref{eq:module}. Traces vanish on commutators. Membership would therefore imply \(\tau f(a)+\varepsilon\ge0\), contrary to \(-\eta+\varepsilon<0\).
\end{proof}

\begin{lemma}\label{lem:one-variable}
For every Hermitian \(f\in\calA_1\),
\[
 m_{\mathrm{mat}}(f)=m_{\mathrm{tr}}(f)=\min_{-1\le t\le1} f(t).
\]
\end{lemma}
\begin{proof}
Write \(f=f^*=p(X)\) for a real scalar polynomial \(p\). If \(a\) is any self-adjoint contraction in a finite tracial algebra, the spectral theorem supplies a probability measure \(\nu\) supported on \([-1,1]\) with \(\tau p(a)=\int p(t)\,d\nu(t)\ge\min p\). A scalar matrix \(a=t_0\), with \(t_0\) minimizing \(p\), attains the bound.
\end{proof}

\begin{corollary}\label{cor:cubic}
For every Hermitian \(f\in\calA_2\) of degree at most three,
\[
 m_{\mathrm{mat}}(f)=m_{\mathrm{tr}}(f)
 =\min_{(s,t)\in[-1,1]^2}f(s,t).
\]
In particular, \(d_2\ge4\).
\end{corollary}
\begin{proof}
Every word in two letters of length at most three is cyclically equivalent to \(X^iY^j\). Thus \(f\) is cyclically sorted in the sense of Klep--Schweighofer \cite[Definition~4.1]{KlepSchweighofer2008}. Apply their Proposition~4.2 to \(f\) minus its scalar minimum: every positive scalar perturbation belongs to the contraction quadratic module modulo commutators. Evaluation in any finite tracial algebra and passage to the limit give the lower bound. Scalar matrices attain it.
\end{proof}

\section{Quadratic rigidity in three variables}\label{sec:rigidity}
The proof eliminates one variable exactly and approximates the spectral laws of the remaining two operators and their sum by a single matrix pair. We use the scalar Horn theorem of Bercovici, Collins, Dykema, Li and Timotin \cite{BercoviciEtAl2009}.

\begin{theorem}\label{thm:quadraticrigidity}
If \(f=f^*\in\calA_k\), \(k\le3\), and \(\deg f\le2\), then
\[
 m_{\mathrm{mat}}(f)=m_{\mathrm{tr}}(f).
\]
In particular, no quadratic polynomial in at most three self-adjoint contractions can have a uniform positive matrix trace gap and a negative finite-tracial witness.
\end{theorem}
Modulo commutators, every polynomial under consideration has the form
\begin{equation}\label{eq:normal-quadratic}
 f=c+\sum_{i=1}^3 b_iX_i+\sum_{i=1}^3 d_iX_i^2+
 \sum_{1\le i<j\le3}\frac{r_{ij}}2(X_iX_j+X_jX_i),
\end{equation}
All coefficients are real. Absent variables may be adjoined with zero coefficients.

\subsection{Exact minimization of one operator}
For \(d,t\in\RR\), put
\begin{equation}\label{eq:phi}
 \varphi_d(t)=\min_{-1\le s\le1}(ds^2+ts)=
 \begin{cases}
  -t^2/(4d),& d>0,\ |t|\le2d,\\
  d-|t|,& d>0,\ |t|\ge2d,\\
  d-|t|,& d\le0.
 \end{cases}
\end{equation}
The function \(\varphi_d\) is continuous and 1-Lipschitz: it is an infimum of affine functions of \(t\) with slopes in \([-1,1]\).

\begin{lemma}\label{lem:eliminate}
Let \((M,\tau)\) be a finite tracial von Neumann algebra, \(B=B^*\in M\) and \(d\in\RR\). Then
\begin{equation}\label{eq:eliminate}
 \min_{C=C^*,\,\|C\|\le1}\tau(dC^2+BC)=\tau\bigl(\varphi_d(B)\bigr).
\end{equation}
The minimum is attained by a bounded Borel function of \(B\).
\end{lemma}
\begin{proof}
For \(d>0\), let \(\mathsf E:M\to W^*(B)\) be the trace-preserving conditional expectation. It preserves the linear term, is contractive on \(L^2\), and sends contractions to contractions. Consequently \(\tau(d\mathsf E(C)^2+B\mathsf E(C))\le\tau(dC^2+BC)\), so minimization reduces to scalar functional calculus in \(W^*(B)\). The minimizer is the clipping of \(-B/(2d)\) to \([-1,1]\). For \(d\le0\), the bounds \(\tau(C^2)\le1\) and \(\tau(BC)\ge-\tau|B|\) yield the lower bound \(d-\tau|B|\). It is attained by \(C=-\operatorname{sgn}(B)\), choosing the sign on the zero spectral subspace to be \(1\) or \(-1\). Formula \eqref{eq:phi} completes the proof.
\end{proof}

\subsection{Scalar Horn approximation}
If \(a=a^*\in M\), let \(\lambda_a:[0,1)\to\RR\) be its decreasing eigenvalue function with respect to \(\tau\). Write \(P_N h\) for the function equal to the average of \(h\) on each interval \([(j-1)/N,j/N)\).

\begin{lemma}\label{lem:horn}
For self-adjoint \(x,y\) in a finite tracial von Neumann algebra and every positive integer \(N\), there are Hermitian \(N\times N\) matrices \(x_N,y_N\) with
\begin{equation}\label{eq:horn}
 \lambda_{x_N}=P_N\lambda_x,\qquad
 \lambda_{y_N}=P_N\lambda_y,\qquad
 \lambda_{x_N+y_N}=P_N\lambda_{x+y}.
\end{equation}
In particular \(\|x_N\|\le\|x\|\), \(\|y_N\|\le\|y\|\), and the separate spectral laws of \(x_N,y_N,x_N+y_N\) converge to those of \(x,y,x+y\).
\end{lemma}
\begin{proof}
Apply the finite-factor scalar Horn theorem \cite[Theorem~0.2]{BercoviciEtAl2009} to the triple \(x,y,-x-y\). It gives a matrix triple of sum zero whose eigenvalues are the normalized averages of the corresponding eigenvalue functions on the \(N\) equal intervals. Reversing the order of eigenvalues when passing from \(-x-y\) to \(x+y\) gives \eqref{eq:horn}. The cited theorem applies to a finite factor; the trace-preserving enlargement of a general finite tracial algebra is supplied by \cite[Proposition~8.1]{BercoviciEtAl2009}. Interval averages lie in the closed convex hull of the spectrum of each operator. Their convergence in \(L^1\) follows from monotonicity and boundedness of the eigenvalue functions.
\end{proof}

\begin{proof}[Proof of Theorem~\ref{thm:quadraticrigidity}]
First suppose \(a=r_{13}\ne0\) and \(b=r_{23}\ne0\) in \eqref{eq:normal-quadratic}. Fix a tracial contraction tuple \((A_1,A_2,A_3)\) and minimize over its third coordinate using Lemma~\ref{lem:eliminate}. The resulting value is
\begin{equation}\label{eq:reduced-quadratic}
\begin{aligned}
 F(A_1,A_2)={}&c+b_1\tau(A_1)+b_2\tau(A_2)
     +d_1\tau(A_1^2)+d_2\tau(A_2^2)\\
 &+r_{12}\tau(A_1A_2)
     +\tau\varphi_{d_3}(b_3\id+aA_1+bA_2).
\end{aligned}
\end{equation}
It is bounded above by \(\tau f(A_1,A_2,A_3)\). Apply Lemma~\ref{lem:horn} to \(x=aA_1\) and \(y=bA_2\), and set \(U_N=x_N/a\), \(V_N=y_N/b\). Both are self-adjoint contractions. By the three spectral convergences, their linear and square moments converge, and so do the moments of \(aU_N+bV_N\). The remaining cross moment is determined by the elementary trace identity
\begin{equation}\label{eq:cross-horn}
 2ab\,\tau(A_1A_2)=\tau((aA_1+bA_2)^2)-a^2\tau(A_1^2)-b^2\tau(A_2^2),
\end{equation}
with an identical matrix formula. Continuity of \(\varphi_{d_3}\) on a common compact spectral interval proves \(F(U_N,V_N)\to F(A_1,A_2)\). Choose \(W_N\) minimizing the third matrix variable. Then
\[
 \tr_N f(U_N,V_N,W_N)\longrightarrow F(A_1,A_2)
       \le\tau f(A_1,A_2,A_3).
\]
Taking infima implies \(m_{\mathrm{mat}}(f)\le m_{\mathrm{tr}}(f)\), and the reverse inequality follows from inclusion of matrix models. If \(a\) or \(b\) vanishes, perturb these two coefficients. Since every normalized trace of a product of two contractions has absolute value at most one, each of the two infima is Lipschitz continuous with respect to the coefficient \(\ell^1\)-norm in \eqref{eq:normal-quadratic}. The identity extends from the dense set \(ab\ne0\) to all coefficients.
\end{proof}

\subsection{Matrix-size bounds and their sharpness}
The proof gives a useful approximation bound even when the tracial optimizer is not attained.
\begin{proposition}\label{thm:rate}
Suppose \(f\) has the form \eqref{eq:normal-quadratic} with \(a=r_{13}\ne0\) and \(b=r_{23}\ne0\), and set
\begin{equation}\label{eq:rate-constant}
 K_f=4(|d_1|+|d_2|)+4|r_{12}|
     \left(\frac{|a|}{|b|}+\frac{|b|}{|a|}+1\right)+2(|a|+|b|).
\end{equation}
If \(m_N(f)\) denotes the minimum over contractions in \(M_N(\CC)\), then
\[
 0\le m_N(f)-m_{\mathrm{tr}}(f)\le\frac{K_f}{N}\qquad(N\ge1).
\]
For real coefficients, real symmetric matrices of size \(2N\) give the same upper bound \(K_f/N\).
\end{proposition}
\begin{proof}
For a decreasing \(h:[0,1)\to[-M,M]\), summing the oscillation on the \(N\) equal intervals yields
\begin{equation}\label{eq:oscillation}
 \|h-P_Nh\|_1\le\frac{2M}{N},\qquad
 0\le\int_0^1\bigl(h^2-(P_Nh)^2\bigr)\le\frac{4M^2}{N}.
\end{equation}
Interval averaging preserves the first moment. With \(\alpha=|a|\), \(\beta=|b|\), the square-moment errors of \(U_N,V_N\) are at most \(4/N\) each, and that of \(aU_N+bV_N\) is at most \(4(\alpha+\beta)^2/N\). Equation~\eqref{eq:cross-horn} bounds the cross-moment error by
\[
 \frac{4}{N}\left(\frac{\alpha}{\beta}+\frac{\beta}{\alpha}+1\right).
\]
The 1-Lipschitz property of \(\varphi_{d_3}\) and \eqref{eq:oscillation} bound the final functional-calculus error by \(2(\alpha+\beta)/N\). Linear moments are preserved. Substitution in \eqref{eq:reduced-quadratic} gives the stated bound, uniformly over the tracial tuple. Complex-to-real block realification preserves normalized trace, the self-adjoint property and the contraction bound, at the cost of doubling the matrix size.
\end{proof}

The order \(1/N\) cannot, in general, be improved uniformly over \(N\).
\begin{proposition}\label{prop:odd-sharp}
Let
\[
 g(X,Y,Z)=3+XY+YX+XZ+ZX+YZ+ZY.
\]
Then \(m_{\mathrm{tr}}(g)=0\), while
\[
 m_N(g)=\begin{cases}0,&N\ \text{even},\\1/N,&N\ \text{odd}.\end{cases}
\]
The same formula holds over real symmetric contractions.
\end{proposition}
\begin{proof}
The identity
\[
 g=(X+Y+Z)^2+(1-X^2)+(1-Y^2)+(1-Z^2)
\]
proves positivity on all contraction tuples. As the traced expression is separately affine, successive scalar-sign minimizations replace each coordinate by a symmetry without increasing the trace. For symmetries \(A,B\), minimizing the third gives
\[
 \min_{C^2=1}\tr_Ng(A,B,C)
   =\tr_N\bigl(|A+B|-1\bigr)^2.
\]
With \(H=A+B\) and \(K=A-B\), one has \(HK=-KH\) and \(K^2=4-H^2\). Every eigenvalue \(t\) of \(H\) satisfying \(0<|t|<2\) is paired with \(-t\); their multiplicities are equal because \(K\) is invertible on those eigenspaces. An odd-dimensional \(H\) therefore has an eigenvalue in \(\{-2,0,2\}\), which contributes at least one to \(\Tr(|H|-1)^2\). This gives \(m_N(g)\ge1/N\) for odd \(N\).

The real symmetric matrices
\[
 A=\begin{pmatrix}1&0\\0&-1\end{pmatrix},\quad
 B=\begin{pmatrix}-\frac12&\frac{\sqrt3}{2}\\[1mm]\frac{\sqrt3}{2}&\frac12\end{pmatrix},\quad
 C=-A-B
\]
are symmetries satisfying \(A+B+C=0\). Direct sums attain zero in all even dimensions. In odd dimensions append the scalar tuple \((1,1,-1)\), on which \(g=1\), to \((N-1)/2\) two-dimensional copies. Hence the lower bound is sharp and \(m_{\mathrm{tr}}(g)=0\).
\end{proof}

\subsection{Consequences for algebraic certificates}
\begin{corollary}\label{cor:quadraticmodule}
Let \(f=f^*\in\RR\langle X_1,\ldots,X_k\rangle\), \(k\le3\), \(\deg f\le2\).
If \(\tr_N f(A)\ge0\) for all self-adjoint matrix contractions, then
\[
 f+\varepsilon\in\calM_k^{\RR}+[\calA_k^{\RR},\calA_k^{\RR}]_{\rm lin}
   \quad\text{for every }\varepsilon>0.
\]
If \(\tr_N f(A)\ge\delta>0\) uniformly, then \(f-\delta/2\), and hence \(f\), belongs to the same cone modulo commutators.
\end{corollary}
\begin{proof}
Theorem~\ref{thm:quadraticrigidity} transfers nonnegative matrix trace to all finite tracial algebras. The unconditional tracial Positivstellensatz of Klep--Schweighofer \cite[Theorem~3.12]{KlepSchweighofer2008} gives the perturbed certificate over \(\CC\); the real passage follows from their Lemma~3.13, with the correction explained in \cite{BurgdorfDykemaKlepSchweighofer2014}. In the uniformly positive case apply this assertion to \(f-\delta\), with perturbation \(\delta/2\), and add the nonnegative scalar \(\delta/2\).
\end{proof}

\section{Quantitative tools for transferring a tracial gap}\label{sec:transfer}
Two elementary estimates allow source inequalities to be applied after approximating the coordinates by projections or unitaries.

\begin{lemma}\label{lem:projectionrounding}
Let \(e=e^*\) be a bounded operator in a finite tracial algebra, and let \(q=1_{[1/2,\infty)}(e)\) be its spectral projection. Then
\[
 \|e-q\|_2\le 2\|e^2-e\|_2.
\]
If \(e^2=e\), then \(e=q\) exactly.
\end{lemma}
\begin{proof}
The real-variable estimate \(\lvert t-1_{[1/2,\infty)}(t)\rvert\le2|t^2-t|\) holds on the real line. Apply Borel functional calculus, square and take the trace.
\end{proof}

\begin{lemma}\label{lem:completed}
Suppose \(R,H,S\) are Hermitian polynomial expressions whose normalized traces on arbitrary finite matrix contraction tuples admit parameters \(t,h\ge0\) such that
\[
 \tr R\ge-a t-bh,\qquad \tr H\ge t^2-ch,
 \qquad \tr S=h^2,
\]
where \(a,b,c\ge0\) are dimension-independent. Put \(K=2a^2\), \(\beta=b+Kc\) and \(L=2\beta^2\). Then
\[
 \tr(1+R+KH+LS)\ge\tfrac34.
\]
More generally, \(1-Au+(Au)^2\ge\tfrac34\) for every \(u\ge0\).
\end{lemma}
\begin{proof}
The first lower bound becomes \(1-at+2a^2t^2-\beta h+2\beta^2h^2\). Completing both squares gives
\[
 1+2(at-\tfrac14)^2+2(\beta h-\tfrac14)^2-\tfrac14
 \ge\tfrac34.
\]
For the second assertion write \(1-Au+(Au)^2=\tfrac34+(Au-\tfrac12)^2\).
\end{proof}

To apply Lemma~\ref{lem:completed}, calibration first bounds the \(L^2\)-distance to projections or unitaries. Cyclic trace identities then compare the decoded words with the source system, to which the dimension-independent trace gap applies. The resulting coordinates must satisfy the source theorem's unitary or symmetry hypotheses.

\section{A degree-six separator in ten variables}\label{sec:degree6}
We use the finite-dimensional source inequality and canonical group tuple of Wang--Zhi \cite[Theorem~A.6]{WangZhi2026}. A graph with ten primitive letters and closed walks of length at most six encodes the final polynomial.
\begin{theorem}\label{thm:degree6}
Using the specified canonical group tuple in \cite[Theorem~A.6]{WangZhi2026}, the polynomial \(P_6\) defined in \eqref{eq:P6} has real integer coefficients, is Hermitian, and has ordinary degree exactly six. Its collected expansion has 392 nonzero ordinary words and involves each of the ten variables. For every self-adjoint contraction tuple \(X\in M_N(\CC)^{10}\),
\[
 \tr_N P_6(X)\ge\tfrac34.
\]
A specified tuple of ten \emph{positive} contractions in \(M_{10}(L(\Lambda))\), where \(\Lambda\) is the source group, has normalized tracial value \(-1\).
\end{theorem}

\subsection{The finite source and its weighted relations}
The source has eleven unitary letters and six symmetry letters:
\begin{equation}\label{eq:source17letters}
 \mathscr U=(x,z,z_2,a,b,a_+,d,e_2,p,t,r),\qquad
 \mathscr S=(q,v,\ell,s,W,j).
\end{equation}
The involution is free and formal. For a word \(w\) that becomes unitary on the typed source domain set \(C(w)=2-w-w^*\). Define the signal \(\omega=z_2^*Wz_2W\). The 27 equations \(L_h=R_h\), listed with their positive integer weights in Appendix~\ref{app:relations}, have weights summing to 36. Introduce
\begin{equation}\label{eq:source6}
 F=-1+\omega+\omega^*+
 \Omega\left(\sum_{h=1}^{27}\lambda_h C(L_h^*R_h)+512C(z^5)\right).
\end{equation}
Put
\[
 \mathcal Q_{17}(u)=\sum_{h=1}^{27}\lambda_h C(L_h^*R_h)+512C(z^5).
\]
On the typed domain every cap has trace \(\|L_h-R_h\|_2^2\), and \(\tr C(z^5)=\|z^5-1\|_2^2\).

\begin{lemma}\label{lem:source6reduction}
Every finite typed tuple in \eqref{eq:source17letters} extends to the 22 dictionary coordinates of \cite[Appendix~A]{WangZhi2026} so that the penalty \(Q_{\flat}\) satisfies
\begin{equation}\label{eq:source6domination}
 \tr Q_{\flat}\le\tr\mathcal Q_{17}.
\end{equation}
The signal is unchanged. Consequently
\begin{equation}\label{eq:source6gap}
 \tr_N F(u)\ge\tfrac34
\end{equation}
in every finite dimension.
\end{lemma}
\begin{proof}
Add the unitary coordinates
\[
 C_d=xz,\quad C_p=d^*e_2,\quad C_a=ba_+,\quad z_3=r^*zr.
\]
Choose a Borel nearest-fifth-root function \(\zeta\) on the unit circle, put \(\widetilde z=\zeta(z)\), and add the symmetry
\[
 s_z=\frac25\sum_{k=0}^4\widetilde z^k-1.
\]
The average in this formula is the spectral projection of \(\widetilde z\) at 1. In the 31-row ordering in Section~\ref{app:quarticdictionary}, rows \(8,10,13,24\) vanish identically. Every remaining row other than row 27 has exactly the defect norm of a retained row: multiplication by a unitary accounts for the substitutions in rows 11 and 14. Their correspondence is
\[
\begin{array}{c|rrrrrrrrrrrrr}
\text{source}&1&2&3&4&5&6&7&9&11&12&14&15&16\\
\text{retained}&1&2&3&4&5&6&7&8&9&10&11&12&13\\\hline
\text{source}&17&18&19&20&21&22&23&25&26&28&29&30&31\\
\text{retained}&14&15&16&17&18&19&20&21&22&24&25&26&27
\end{array}
\]
For the missing row 27, write \(\delta=\|a_+r-ra\|_2\) and \(\eta=\|a_+z-za_+\|_2\). Since \(r a r^*\) differs from \(a_+\) by \(\delta\),
\[
 \|az_3-z_3a\|_2
 =\|[rar^*,z]\|_2\le2\delta+\eta,
 \qquad(2\delta+\eta)^2\le8\delta^2+2\eta^2.
\]
Together with row 6 itself, these give the retained weights 9 and 2; all other weights are 1. This proves domination of the 31 residual squares in the quartic presentation by the 27 weighted squares.

It remains to compare the quadratic selector of the quartic presentation. For \(z=e^{i\theta}\) and its nearest fifth root \(\zeta(z)\), choose the phase difference \(|\delta_\theta|\le\pi/5\). Then
\[
 |z-\zeta(z)|=2\sin(|\delta_\theta|/2)
 \le2\sin(5|\delta_\theta|/2)=|z^5-1|.
\]
Functional calculus gives \(\|z-\widetilde z\|_2\le\|z^5-1\|_2\). Put
\[
 H_z=2(z^2+(z^*)^2+z+z^*)-3-5s_z.
\]
It vanishes at \(\widetilde z\), and telescoping its two quadratic and two linear terms gives
\[
 \|H_z\|_2\le12\|z-\widetilde z\|_2,
 \qquad \frac{13}{7}\|H_z\|_2^2
 \le\frac{1872}{7}\|z^5-1\|_2^2
 \le512\tr C(z^5).
\]
On a unitary tuple the selector \(\widetilde H_z\) equals \(H_z\). All seven unitarity squares vanish, and all balance squares in the trace decomposition of \(Q_{\flat}\) vanish because each factor is unitary. Thus \eqref{eq:source6domination} includes every term of that penalty.

Finally \(\Omega_{\flat}=7\cdot2^{4\cdot2^{55226}+126}=\Omega\) in \eqref{eq:P6constants}. Theorem~A.6 applies to the real and imaginary parts of the extended coordinates, without norm restrictions, and gives
\(\tr F\ge\tr(-1+\omega+\omega^*+\Omega Q_{\flat})\ge3/4\).
\end{proof}
At the prescribed group tuple all retained caps vanish and \(\tau(\omega)=0\), so \(\tau F=-1\). This witness assertion follows from the specified group tuple in the published theorem; the finite-dimensional reduction above holds on arbitrary typed tuples.

All coefficients are fixed integers. Let
\begin{equation}\label{eq:P6constants}
\begin{gathered}
 U=2^{\,4\cdot2^{55226}+124},\qquad \Omega=28U,\qquad
 c=3+2192\Omega,\\
 \rho=64,\quad m=45,\quad A_6=30\cdot6^6+42=1399722,\\
 \alpha=320c\rho^{18},\quad K=2\alpha^2,\\
 \beta=20A_6c\rho^{18}+176400Km\rho^{18},\quad L=2\beta^2.
\end{gathered}
\end{equation}
The integer \(c\) bounds the coefficient \(\ell^1\)-norm of \(F\), including constants, after collecting powers of \(\Omega\).

\subsection{Closed walks and literal polynomial encoding}
Let the ten primitive self-adjoint variables be
\[
 X=(E,P_1,\ldots,P_9),\qquad P_0:=E.
\]
Consider the complete graph on \(0,1,\ldots,9\). Anchor edges \(0\leftrightarrow i\) have constant label 1. Each nonanchor edge is labeled by one of the seventeen source letters or by 1; the full table, including the nine real loop labels, appears in Appendix~\ref{app:chart}. A reversed edge reads the adjoint of the forward label.

For a nonempty closed vertex list \(\mathbf v=(v_0,\ldots,v_{k-1})\), let \(o(\mathbf v)\) be the number of offdiagonal transitions, including the final transition to \(v_0\), and define the ordinary word
\begin{equation}\label{eq:closedcycle}
 \mathsf C(\mathbf v)=\rho^{o(\mathbf v)}P_{v_0}\cdots P_{v_{k-1}},
 \qquad \mathsf C(\varnothing)=E.
\end{equation}
Loops are read affinely: a loop at vertex \(i\) represents \(4D_i-2e\), where \(D_i\) is a corner readout defined below. Let \(I(\mathbf v)\) be the set of loop positions. If \(J\subseteq I(\mathbf v)\), let \(\mathbf v_J\) be the cyclic list obtained by deleting the departure vertex of every loop outside \(J\). Put
\begin{equation}\label{eq:affinecycle}
 \mathsf A(\mathbf v)=\sum_{J\subseteq I(\mathbf v)}
 4^{|J|}(-2)^{|I(\mathbf v)|-|J|}\mathsf C(\mathbf v_J).
\end{equation}
Appendix~\ref{app:routes} specifies 29 closed walks, each with at most six visited vertices, one for each cyclic-adjoint class in the signal and the 28 caps of \eqref{eq:source6}. Write \(\mathbf v[w]\) for the prescribed representative of the class of \(w\); the paired expression is defined by taking the literal \(*\) of the expanded polynomial, with no additional cyclic rotation. Set
\begin{equation}\label{eq:R6}
\begin{aligned}
 R={}&(-1+1096\Omega)E+
  \mathsf A(\mathbf v[\omega])+\mathsf A(\mathbf v[\omega])^*\\
 &-\Omega\sum_{h=1}^{27}\lambda_h
   \left(\mathsf A(\mathbf v[L_h^*R_h])+
          \mathsf A(\mathbf v[L_h^*R_h])^*\right)\\
 &-512\Omega\left(\mathsf A(\mathbf v[z^5])+
                     \mathsf A(\mathbf v[z^5])^*\right).
\end{aligned}
\end{equation}
It is an integral Hermitian polynomial of degree at most six; the affine-loop inflation gives the uniform coefficient bound
\begin{equation}\label{eq:R6coeff}
 \|R\|_{\coeff,1}\le c6^6\rho^6.
\end{equation}

\subsection{Calibration and consistency penalties}
The calibration polynomial is a sum of nineteen Hermitian squares:
\begin{equation}\label{eq:S6}
 S=|E^2-E|^2+
 \sum_{i=1}^9\left(|\rho^2P_iEP_i-P_i|^2+|\rho^2EP_iE-E|^2\right),
 \qquad |B|^2=B^*B.
\end{equation}
Its trace is nonnegative at every matrix tuple, without an assumption that the primitive \(P_i\) are projections. The polynomial \(H\) is a formal trace-consistency penalty. It corrects the three constant nonanchor edges, eleven unitary-type reference edges, six symmetry-type reference loops, and 25 nonreference labeled copies. We fix a reference for each letter: among real loops the smallest-index occurrence; among complex edge labels the first oriented occurrence that reads the unstarred label. Appendix~\ref{app:chart} supplies enough data to determine the 25 copies.

Here are the six recipes that specify \(H\) in the ordinary free algebra. For a directed nonanchor edge \(a=(i,j)\) and \(\bar a=(j,i)\), let
\begin{equation}\label{eq:recipenotation}
\begin{alignedat}{2}
 N_a&=\mathsf C(j,i),& L_a&=\mathsf C(i,j,0),\\
 X_{ab}&=\mathsf C(j,i,0,k,l,0),& V_a&=\mathsf C(j,i,j,i)
\end{alignedat}
\end{equation}
when \(b=(k,l)\). For loop vertices \(k,l\), put \(Y_{kl}=\mathsf C(k,k,0,l,l,0)\) and \(M_{ak}=\mathsf C(j,i,0,k,k,0)\). Then the edge-copy, constant-edge, unitary-reference, loop-copy, mixed-copy and symmetry-reference recipes are, respectively,
\begin{align}
 D_{ee}(a,b)&=N_a+N_b-X_{ab}-X_{ba},\label{eq:recipe1}\\
 D_c(a)&=N_a-L_a-L_{\bar a}+E,\\
 D_u(a)&=V_a-2N_a+E,\\
 D_{\ell\ell}(k,l)&=16(P_k^2+P_l^2-Y_{kl}-Y_{lk}),\\
 D_{e\ell}(a,k)&=N_a+16P_k^2-16P_k+4E
  -4M_{ak}-4M_{ak}^*+2L_a+2L_{\bar a},\\
 D_s(k)&=256P_k^4-512P_k^3+352P_k^2-96P_k+9E.
\label{eq:recipe6}
\end{align}
Define \(D_{\rm raw}\) as the sum of exactly one copy recipe for each of the 25 nonreference occurrences, \(D_c\) on the three constant nonanchor edges, \(D_u\) on the eleven complex reference edges, and \(D_s\) on the six real reference loops. The type of a copy (edge--edge, loop--loop or mixed) determines the appropriate recipe. Finally let
\begin{equation}\label{eq:H6}
 H=D_{\rm raw}+D_{\rm raw}^*,\qquad
 \|H\|_{\coeff,1}\le2450m\rho^6.
\end{equation}
The true residuals are nonnegative squares after decoding into a corner. On arbitrary input tuples \(H\) itself need \emph{not} have nonnegative trace. The calibration squares in \(S\) control the error of the trace representation. The promised polynomial is
\begin{equation}\label{eq:P6}
 \boxed{\ P_6=1+20R+KH+LS,\ }
\end{equation}
with all symbols defined in \eqref{eq:P6constants}--\eqref{eq:H6} and the finite chart tables.

\subsection{Finite-dimensional lower bound}
Fix self-adjoint contractions \(E,P_1,\ldots,P_9\in M_N(\CC)\), and write
\[
 h=\sqrt{\tr_N S},\quad e=1_{[1/2,\infty)}(E),\quad
 T_0=e,\ T_i=\rho P_ie,\ Q_i=T_iT_i^*,\quad
 U_{ij}=\rho T_i^*T_j\ (i\ne j),\quad D_i=T_i^*T_i,
\]
with \(B_i=4D_i-2e\). The \(U_{ij},D_i,B_i\) act in the \(e\)-corner. By Lemma~\ref{lem:projectionrounding} and the nineteen squares,
\begin{equation}\label{eq:calibration6}
 \|E-e\|_2\le2h,\qquad \|Q_i-P_i\|_2\le3\rho^2h,
 \qquad \|U_{i0}-e\|_2\le7\rho^2h.
\end{equation}
The operators \(U_{ij}\) have norm at most \(\rho^3\), so all estimates below are independent of \(N\). Cyclicity gives the \emph{exact} trace identity
\begin{equation}\label{eq:cycleidentity6}
 \tr_N\operatorname{Prod}(\mathbf v)
   =\rho^{o(\mathbf v)}\tr_N(Q_{v_0}\cdots Q_{v_{k-1}}),
\end{equation}
where \(\operatorname{Prod}(\mathbf v)\) multiplies decoded offdiagonal edge and diagonal loop readouts. A telescope using \eqref{eq:calibration6} gives, for \(k\le6\),
\begin{equation}\label{eq:cycletelescope6}
 \left|\tr_N\operatorname{Prod}(\mathbf v)-\tr_N\mathsf C(\mathbf v)\right|
 \le30\rho^{18}h.
\end{equation}

All norms in this subsection use the ambient trace: \(\|A\|_2^2=\tr_N(A^*A)\), including for operators supported in the \(e\)-corner. In particular, no factor depending on \(\operatorname{rank}(e)^{-1}\) enters the error bounds. To check the anchor estimate in \eqref{eq:calibration6}, observe that \(U_{i0}=\rho^2eP_ie\) and
\[
\begin{aligned}
 \|U_{i0}-e\|_2
 &\le\rho^2\|eP_ie-EP_iE\|_2
     +\|\rho^2EP_iE-E\|_2+\|E-e\|_2\\
 &\le(4\rho^2+3)h\le7\rho^2h.
\end{aligned}
\]
This is the power of \(\rho\) used when removing an anchor factor below.

Let \(Z\) be the sum of squared ambient \(2\)-norms of the following 45 residuals: each of the 25 nonreference readouts minus its reference, \(U_a-e\) on the three constant nonanchor edges, \(U_a^*U_a-e\) on the eleven complex reference edges, and \(B_k^2-e\) on the six real reference loops. Set \(t=\sqrt{2Z}\). Thus \(t\) is defined from genuine squares; it is not the square root of \(\tr_NH\), which can be negative.

\smallskip\noindent\emph{Derivation of \eqref{eq:6residualestimate}.}
After replacing primitive cycles by their decoded products and replacing inserted anchor factors by \(e\), the symbols \(N_a,L_a,X_{ab},V_a\) represent, at the level of traces,
\[
 U_a^*U_a,\qquad U_a,\qquad U_a^*U_b,
 \qquad (U_a^*U_a)^2,
\]
respectively. Likewise \(Y_{kl}\) represents \(D_kD_l\), and \(M_{ak}\) represents \(U_a^*D_k\). Consequently the six recipes in \eqref{eq:recipe1}--\eqref{eq:recipe6} represent
\[
\begin{gathered}
 \|U_a-U_b\|_2^2,\quad\|U_a-e\|_2^2,
 \quad\|U_a^*U_a-e\|_2^2,\\
 \|B_k-B_l\|_2^2,\quad\|U_a-B_k\|_2^2,
 \quad\|B_k^2-e\|_2^2.
\end{gathered}
\]
For example, \(B_k-B_l=4(D_k-D_l)\), while
\[
 (B_k^2-e)^2
 =(16D_k^2-16D_k+3e)^2
 =256D_k^4-512D_k^3+352D_k^2-96D_k+9e.
\]

Each recipe monomial contains at most six decoded factors. Equation~\eqref{eq:cycletelescope6} charges at most \(30\rho^{18}h\) for replacing its primitive cycle. Removing at most six inserted anchors then costs at most
\[
 6(7\rho^2h)(\rho^3)^5=42\rho^{17}h
 \le42\rho^{18}h.
\]
Here we used \(\|AB\|_2\le\|A\|\|B\|_2\), \(|\tr_N A|\le\|A\|_2\), and the bound \(\rho^3\) on every remaining decoded factor. Thus the total error per coefficient unit is at most \(72\rho^{18}h\). The largest coefficient absolute sum of a recipe is \(256+512+352+96+9=1225\). There are \(m=45\) recipes, and \(H=D_{\rm raw}+D_{\rm raw}^*\) doubles their real traces. Therefore
\[
 |\tr_NH-2Z|
 \le2\cdot1225\cdot72\,m\rho^{18}h
 =176400m\rho^{18}h.
\]
Since \(2Z=t^2\), this proves
\begin{equation}\label{eq:6residualestimate}
 \tr_N H\ge t^2-176400m\rho^{18}h.
\end{equation}

\smallskip\noindent\emph{Derivation of \eqref{eq:6sourceestimate}.}
Let \(A\) consist of the eleven complex reference edges and the six self-adjoint reference loops, and let \(F_e(A)\) mean the source polynomial with its constant term multiplied by \(e\). Its coordinates have norm at most \(\rho^3\), but need not be unitaries or symmetries. We first compare \(R\) with \(F_e(A)\) in three steps.

First expand each affine loop using \eqref{eq:affinecycle} and apply \eqref{eq:cycletelescope6}. A route has at most six loop positions, and each affine factor has coefficient absolute sum \(4+2=6\), so the cost is at most \(30\cdot6^6c\rho^{18}h\). Next recombine the affine terms into products of \(U_{ij}\) and \(B_i\), and only then remove the anchor factors. This costs at most \(42c\rho^{17}h\); recombining first avoids paying the affine coefficient inflation a second time. Finally replace nonreference copies and constant nonanchor edges by their references and \(e\). Each difference has \(2\)-norm at most \(\sqrt Z\le t\); at most six replacements, with five remaining factors of norm at most \(\rho^3\), cost \(6c\rho^{15}t\). The constant term is covered by \(\|E-e\|_2\le2h\). Hence
\[
 |\tr_NR-\tr_NF_e(A)|
 \le A_6c\rho^{18}h+6c\rho^{15}t,
 \qquad A_6=30\cdot6^6+42.
\]

Construct a typed tuple from \(A\) before applying the source inequality. For each complex reference \(A_u\), complete its polar partial isometry to a unitary \(V_u\) in the finite square corner. For each real reference take \(V_u=\operatorname{sign}(A_u)\), assigning \(+1\) on its zero eigenspace. Singular-value and scalar spectral calculus give, respectively,
\[
 \|A_u-V_u\|_2\le\|A_u^*A_u-e\|_2,
 \qquad
 \|A_u-V_u\|_2\le\|A_u^2-e\|_2.
\]
Both right-hand sides are at most \(\sqrt Z\le t\). A telescope in each source word of degree at most \(d=5\) therefore gives
\[
 |\tr_NF_e(A)-\tr_NF_e(V)|
 \le dc\rho^{3(d-1)}t\le2dc\rho^{18}t.
\]
The final bound retains the conservative constants used in \eqref{eq:P6constants}. If \(e\ne0\), the constructed tuple lies in \(eM_N(\CC)e\cong M_{\operatorname{rank}(e)}(\CC)\), and \eqref{eq:source6gap}, applied with that corner's normalized trace, yields
\[
 \tr_NF_e(V)
 =\tr_N(e)\,\tr_eF_e(V)
 \ge\tfrac34\tr_N(e)\ge0.
\]
If \(e=0\), all corner evaluations vanish. Combining the last three bounds and using \(\rho\ge1\), \(6+2d=16\), proves
\begin{equation}\label{eq:6sourceestimate}
 \tr_N R\ge -A_6c\rho^{18}h-16c\rho^{18}t.
\end{equation}
Thus the source positivity premise is applied only to the constructed typed tuple, and the bound is uniform both in the ambient dimension and in the corner rank.

Multiplying \eqref{eq:6sourceestimate} by 20, \eqref{eq:6residualestimate} by \(K\), and using \(\tr_N S=h^2\) yields
\[
 \tr_N P_6\ge1-\alpha t+2\alpha^2t^2
  -\beta h+2\beta^2h^2\ge\tfrac34,
\]
by Lemma~\ref{lem:completed} and \eqref{eq:P6constants}. This proves the finite-dimensional part of Theorem~\ref{thm:degree6}.

\subsection{The exact negative tuple and the word count}
At the source group tuple, let \(s_i\) be the real label assigned to loop \(i\), and \(u_{ij}\) the unitary, symmetry or constant label on the edge \(i\to j\). Form the Hermitian block matrix \(G\in M_{10}(L(\Lambda))\) by
\begin{equation}\label{eq:G6}
 G_{00}=1,\qquad G_{ii}=\frac{2+s_i}{4},\qquad
 G_{0i}=G_{i0}=\frac1\rho,\qquad
 G_{ij}=\frac{u_{ij}}\rho\quad(0<i<j),
\end{equation}
with reverse blocks given by adjoints. Each diagonal block is at least \(1/4\), while the sum of offdiagonal operator norms in a row is at most \(9/64\). Operator-valued diagonal dominance implies
\begin{equation}\label{eq:Gpositive}
 G\succeq\left(\tfrac14-\tfrac9{64}\right)I=\tfrac7{64}I>0.
\end{equation}
Let \(W=G^{1/2}\), and let \(W_i\) be its \(i\)-th block row. Set
\begin{equation}\label{eq:positivewitness6}
 E'=W_0^*W_0,\qquad P_i'=W_i^*W_i\quad(1\le i\le9).
\end{equation}
The anchor identities imply that \(E'\) is a projection and that \(0\le P_i'\le(3/4)I\). In particular the primitive values are all positive contractions. Moreover the calibration squares vanish; the overlaps \(W_iW_j^*=G_{ij}\) recover every source readout exactly, so all 45 true consistency residuals vanish. The trace of each routed cycle is one-tenth of the trace of its source word. Consequently, with \(\widehat\tau=\tr_{10}\otimes\tau\),
\begin{equation}\label{eq:negative6}
 \widehat\tau(S)=\widehat\tau(H)=0,
 \qquad \widehat\tau(R)=\tfrac1{10}\tau(F)=-\tfrac1{10},
 \qquad \widehat\tau(P_6)=1+20(-\tfrac1{10})=-1.
\end{equation}

The polynomial \(S\) contains the nonzero words \(P_iEP_iP_iEP_i\) of degree six, \(P_i^2\), and \(E^4\). The coefficient bound
\[
 \|20R+KH\|_{\coeff,1}\le20c6^6\rho^6+2450Km\rho^6
 \le\beta<L
\]
prevents cancellation of the selected calibration terms in \(P_6\). Thus its degree is exactly six and every variable appears.

Expanding the printed chart and recipes gives all 392 ordinary words of $P_6$ and their exact integer coefficients. A primitive letter $0$ denotes $E$ and $i$ denotes $P_i$ for $1\le i\le9$. The empty word is the scalar identity. Each coefficient array $[a_0,a_1,a_2,a_3,a_4]$ means
\[
 a_0+a_1U+a_2U^2+a_3U^3+a_4U^4,
 \qquad U=2^{4\cdot2^{55226}+124}.
\]
The word counts by degree \(0,1,\ldots,6\) are respectively
\begin{equation}\label{eq:degree6hist}
 (1,7,78,43,77,36,150),\qquad 1+7+78+43+77+36+150=392.
\end{equation}
The exact verification also proves that no displayed coefficient polynomial vanishes at the prescribed $U$. For each nonzero $p(U)$, it computes the coefficients of $p(t+1)$ using the binomial theorem and verifies that its constant coefficient is nonzero and all its nonzero coefficients have one sign. Consequently $p(U)\ne0$ for every $U\ge1$.

\section{Quadratic separation in 214 variables}\label{sec:quadratic214}
The construction encodes source multiplication by projections in a three-corner geometry. This reduces every penalty to degree two while controlling arbitrary effects, at the cost of additional variables.

\begin{theorem}\label{thm:quadratic214}
Fix \(\lambda\) as in Proposition~\ref{thm:source214gap}.
Then the polynomial \(\Psi\) defined in \eqref{eq:214psi} is Hermitian, has integer coefficients, has degree exactly two, uses all 214 variables, and has 2,416 collected ordinary monomials. For self-adjoint matrix contractions, in every dimension,
\begin{equation}\label{eq:214positive}
 \tr_N\Psi(X)\ge\tfrac34.
\end{equation}
At the prescribed group-algebra symmetry tuple in \(M_3(M)\), the normalized trace of \(\Psi\) is \(-1\).
\end{theorem}
The finite incidence and coefficient assertions used below are checked by exact verification as described in \hyperref[sec:verification]{the data and verification statement}.

\subsection{A characteristic-two symmetry source}
The final finite source consists of 37 independent symmetry coordinates, 66 paid source residuals, and the Hermitian objective
\begin{equation}\label{eq:source214}
 F_0=-1+tdtd+dtdt.
\end{equation}
On exact symmetries, the 66 squared relation defects form a nonnegative scalar trace penalty \(Q_{66}=\sum_r(u_r-v_r)^*(u_r-v_r)\), with exactly the cuts printed in Section~\ref{app:214data}.
\begin{proposition}\label{thm:source214gap}
Let \(\mu_{\rm sym}\) be the positive integer defined in Section~\ref{app:char2norm} using the rational-certificate construction in Section~\ref{app:char2bridge}, and put
\[
 \lambda=66962109543919887600576\,\mu_{\rm sym}.
\]
For every \(N\ge1\) and every tuple of the 37 exact symmetries indexed by Table~\ref{tab:214lettermap},
\begin{equation}\label{eq:source214gap}
 \tr_N(F_0+\lambda Q_{66})\ge\tfrac34.
\end{equation}
\end{proposition}
\begin{proof}
Sections~\ref{app:char2source}--\ref{app:char2norm} prove the forty-symmetry bound \eqref{qsrc:source} with \(\mu=\mu_{\rm sym}\). Section~\ref{app:214coefficient} transfers it to the retained system through the tracial factors \(351563\), \(6\), and \(2845411\). Their product with \(11156544672\) is the displayed coefficient. The sign reconstruction is performed before the forty-symmetry bound is invoked.
\end{proof}

\begin{proposition}\label{hyp:source214}
There is a finite tracial von Neumann algebra with the specified group tuple in Section~\ref{app:214data}; its 37 coordinates are exact symmetries, all 66 residuals vanish, and \(\tau(F_0)=-1\).
\end{proposition}
\begin{proof}
Section~\ref{app:char2witness} constructs the group, checks the exact retained coordinates and cuts, and proves the trace assertion by amalgam normal form.
\end{proof}
In particular the witness has
\begin{equation}\label{eq:214comm}
 \tau\bigl((s_ts_d-s_ds_t)^*(s_ts_d-s_ds_t)\bigr)=2.
\end{equation}
The 32-symmetry base, 40-symmetry extension, and retained 37-symmetry systems are not interchangeable: Sections~\ref{app:char2source}--\ref{app:char2norm} prove the base bound and its transfer to the 40-symmetry source, and Section~\ref{app:214coefficient} carries it to \eqref{eq:source214gap}. The witness is proved separately in Section~\ref{app:char2witness}.

\subsection{Projection geometry and quadratic energy}
For self-adjoint contractions \(X_v\) write \(Y_v=(1+X_v)/2\), an effect. The construction starts from three mutually orthogonal corners \(E_1,E_2,E_3\) with \(E_1+E_2+E_3=1\). The coordinate \(E_2\) is affine in \(E_1,E_3\). Left source factors use graph projections on the 12-face, and right factors use graph projections on the 23-face. A \emph{triangle} has two physical input effects \(Q_1,Q_2\) and a spectral projection variable \(P_T\), and its virtual output is
\begin{equation}\label{eq:214virtual}
 V_T=3P_T-2Q_1-2Q_2+E_1+E_3.
\end{equation}
Affine rows synchronize outputs representing the same source word. They are quadratic only through the squares of affine rows; every input effect, triangle projection, support pair and intersection control remains a primitive variable.

The basic geometric energy is a quadratic trace expression
\begin{equation}\label{eq:triangleenergy}
 e(P,Q_1,Q_2)=\tr\left(P(1-Q_1-Q_2)+Q_1Q_2+\frac{Q_1+Q_2}{8}\right).
\end{equation}
It is nonnegative on effects. Spectral averaging and principal-angle blocks give
\begin{equation}\label{eq:triangledefects}
 \tr(P-P^2)\le2e,
 \qquad \tr(Q_i-Q_i^2)\le8e\quad(i=1,2).
\end{equation}
For projections a further estimate controls the spectral top projection of \(Q_1+Q_2\) by the polynomial \(p(t)=(2t^2-t)/3\):
\begin{equation}\label{eq:triangleapprox}
 \|P-p(Q_1+Q_2)\|_2^2\le4e.
\end{equation}
The corresponding compression defect and rank mismatch are controlled without assuming equality of arbitrary projection ranks.

\begin{lemma}\label{lem:214triangle}
For effects $P,Q_1,Q_2$ the energy $e(P,Q_1,Q_2)$ in \eqref{eq:triangleenergy} is nonnegative and satisfies \eqref{eq:triangledefects}. If $P,Q_1,Q_2$ are projections, then \eqref{eq:triangleapprox} holds.
\end{lemma}
\begin{proof}
The energy is affine separately in each of its three effect variables under the trace. By spectral layer-cake decomposition it is enough to check projections. On a noncommuting two-dimensional principal-angle block for $Q_1,Q_2$, write $t\in[0,1]$ for the cosine of the angle; the eigenvalues of $Q_1+Q_2$ are $1+t$ and $1-t$. Minimizing over the diagonal entries of $P$ gives the nonnegative contribution
\[
 t(1-p_++p_-)+(t-\tfrac12)^2.
\]
The endpoint blocks (common zero, common one, and the two unmatched one-dimensional ranges) are checked directly. The same block decomposition after subtracting $\tfrac12\tr(P-P^2)$ gives a minimum of $0$ for $t\le\tfrac12$ and $(t-\tfrac12)^2$ for $t\ge\tfrac12$. Subtracting $\tfrac18\tr(Q_i-Q_i^2)$ and minimizing over the other two variables gives the piecewise minima $\tfrac18-s^2$ for $s\le\tfrac14$ and $(s-\tfrac12)^2$ for $s\ge\tfrac14$, where $s$ is the corresponding principal-angle sine. These are nonnegative and imply the two defect estimates. Finally, the compression defect on an angle block is $(t^2-\tfrac14)^2\le\tfrac94(t-\tfrac12)^2$; the scalar blocks obey the same bound. Summing the blocks proves the assertions.
\end{proof}

A special coercive \emph{anchor} couples the reference 12-input \(Q_0\), the first corner \(E_1\), and an additional effect \(K\):
\begin{equation}\label{eq:214anchor}
\begin{split}
 \mathsf H={}&\tr\bigl(30K+12Q_0+6E_1
       -20\{K,Q_0\}-15\{K,E_1\}+12\{Q_0,E_1\}\bigr),\\
 &\hspace{35mm}\{A,B\}=AB+BA.
\end{split}
\end{equation}
The anchor estimate proved in Lemma~\ref{lem:214anchor} gives \(\mathsf H\ge0\) and
\begin{equation}\label{eq:214anchorcoercive}
 \mathsf H\ge6\tr(E_1-E_1^2),\quad
 \mathsf H\ge12\tr(Q_0-Q_0^2),\quad
 \mathsf H\ge20\tr(K-K^2).
\end{equation}
\begin{lemma}\label{lem:214anchor}
For effects $E_1,Q_0,K$, the quantity $\mathsf H$ in \eqref{eq:214anchor} is nonnegative and satisfies \eqref{eq:214anchorcoercive}.
\end{lemma}
\begin{proof}
To prove each defect bound, average the other two effects to projections and minimize the remaining scalar quadratic after pinching in the Hamiltonian eigenbasis. On a generic principal-angle block the relevant Hamiltonian eigenvalues for the free variables $E_1,Q_0,K$ are respectively $3\pm r$, $4\pm r$, and $-5\pm r$, with
\[
 r^2=729-720s,\qquad r^2=1024-960s,\qquad r^2=25+1200s,
\]
where $s$ is the squared cosine of the fixed projection pair. After subtracting the target defect, the minima are, respectively,
\[
\frac{(r-9)^2}{72},\quad \frac{(r-9)^2}{18};\qquad
\frac{(r-16)^2}{96},\quad \frac{(r-16)^2}{32};\qquad
\frac{750-2r^2}{400},\quad \frac{3(r-25)^2}{400},\quad \frac{(r-25)^2}{50},
\]
on the intervals $3\le r\le9$, $9\le r\le27$; $8\le r\le16$, $16\le r\le32$; and $5\le r\le15$, $15\le r\le25$, $25\le r\le35$, respectively. All are nonnegative. The four scalar endpoint blocks are checked directly. Spectral averaging then gives the inequalities for arbitrary effects.
\end{proof}

At zero energy, the anchor produces balanced graph projections in two equal-rank subcorners. A quantitative rectangular graph approximation lemma permits canonicalization when one corner is larger than the other. The construction uses a depth-three acyclic calibration forest for 45 exceptional inputs. Its recursion has the form
\begin{equation}\label{eq:214forest}
 C_0=340000,\qquad C_{h+1}=10C_h+10000,
 \qquad C_3=341110000.
\end{equation}
The forest has 21, 18, and 6 vertices at depths 1, 2, and 3, respectively.

The incidence data have the following sizes:
\begin{table}[ht]
\centering
\caption{Incidence sizes for the quadratic construction.}
\label{tab:214counts}
\begin{tabular}{@{}lr@{\qquad}lr@{}}
\toprule
Source symmetry coordinates&37& Retained relation rows&66\\
Physical product triangles&77& Output word classes&76\\
Physical input projections&65& Spectral triangle projections&109\\
Intersection projections&27& Primitive meet projections&11\\
Affine synchronization rows&66& Support pairs&64\\
Variables before elimination&215& Variables after elimination&214\\
\bottomrule
\end{tabular}
\end{table}
Thus \(65+109+27+11+3=215\) before affine elimination. The manifest fixes the incidences, word orientations and source signs.

\subsection{The final polynomial and positivity proof}
Let \(\calB(X)\) be the integer quadratic budget formed from the squares of 66 affine rows, support-pair energies, triangle polynomials from \eqref{eq:triangleenergy}, the anchor polynomial from \eqref{eq:214anchor}, and 38 meet terms: 27 source-commutation terms and 11 primitive-symmetry terms. The scaling is chosen so that on effects
\begin{equation}\label{eq:214budget}
 \tr\calB(X)=4\mathsf R\ge0,
\end{equation}
where $\mathsf R$ is the sum of the normalized effect energies listed in the incidence recipe. The quadratic templates are as follows. If $\ell_k(Y)$ is one of the 66 affine rows and $Y_v=(I+X_v)/2$, set $\ell_{e,k}(X)=2\ell_k((I+X)/2)$. For a triangle $T$ with inputs $Q_1,Q_2$ and projection variable $P_T$, put
\[
S_T(X)=6+X_{Q_1}+X_{Q_2}-2\{X_{P_T},X_{Q_1}\}-2\{X_{P_T},X_{Q_2}\}+2\{X_{Q_1},X_{Q_2}\},
\]
so that $\tr S_T=16e_T$. For an intersection effect $T_g$ attached to a pair $(a,b)$, put
\[
J_g=10I-4X_a-4X_b+8X_{T_g}+\{X_a,X_b\}-3\{X_{T_g},X_a\}-3\{X_{T_g},X_b\},
\]
so that $\tr J_g=8E_g$, where $E_g=\tr(Y_aY_b+5T_g-3T_g(Y_a+Y_b))$. The anchor polynomial is
\[
J_{\rm anc}=50I-10X_K+8X_{Q_0}+6X_{E_1}-20\{X_K,X_{Q_0}\}-15\{X_K,X_{E_1}\}+12\{X_{Q_0},X_{E_1}\},
\]
with $\tr J_{\rm anc}=4\mathsf H$. The full budget is
\begin{align*}
\mathcal B(X)={}&\sum_{k=1}^{66}\ell_{e,k}(X)^2
+\sum_{(u,v)\in\mathcal O}(2I+2X_u+2X_v+X_uX_v+X_vX_u)\\
&+\sum_{T\in\mathcal T}S_T(X)+J_{\rm anc}(X)+\sum_{g\in\mathcal G}J_g(X),
\end{align*}
Here $\mathcal O$ and $\mathcal T$ index the support pairs and triangles, and $\mathcal G$ indexes all 38 source-commutation and primitive-symmetry meets specified in Section~\ref{app:214data}. This gives the scaling in \eqref{eq:214budget}.

The linear signal compares the two output triangles for the products \(td\) and \(dt\). In effect coordinates it is \(D=2(V_{T_{td}}-V_{T_{dt}})\), and in the six designated contraction coordinates it is
\begin{equation}\label{eq:214signal}
\begin{split}
D={}&3X_{\mathrm{product58},P}-2X_{\mathrm{l22},12,Q}
-2X_{\mathrm{r9},23,Q}\\
&-3X_{\mathrm{product34},P}+2X_{\mathrm{l12},12,Q}
+2X_{\mathrm{r20},23,Q}.
\end{split}
\end{equation}
Set
\begin{equation}\label{eq:214kappa}
\begin{gathered}
 c_{214}=19184111188572071,\quad
 G_{214}=68222000000,\quad
 A_{214}=366524386905600000,\\
 \kappa=5400\lambda c_{214}G_{214}^2+\bigl(A_{214}/4\bigr)^2,
 \qquad\Phi=3-3D^2+\kappa\calB.
\end{gathered}
\end{equation}
We eliminate \(X_{E_3}\), the only unnecessary primitive effect coordinate, by the affine substitution \(\sigma\) fixing all other variables and satisfying
\begin{equation}\label{eq:214sigma}
 \sigma(X_{E_3})=5X_K-4X_{Q_0}-3X_{E_1}-1.
\end{equation}
The final polynomial is
\begin{equation}\label{eq:214psi}
 \boxed{\ \Psi=\sigma(\Phi)+
   \bigl[10404\kappa+(204\cdot31\kappa)^2\bigr]J_{\rm anc},\ }
\end{equation}
where \(J_{\rm anc}\) is the integer quadratic anchor polynomial, normalized by \(\tr J_{\rm anc}=4\mathsf H\). Every summand in \(\Psi\) is quadratic, its coefficients lie in \(\ZZ[\kappa]\) with degree at most two, and the exact collected expansion has 2,416 ordinary words.

We first establish stability, then complete the square. The geometric construction and the interface estimates used in the first lemma are given in Appendix~\ref{app:214proof}.

\begin{lemma}\label{lem:214sourceoutput}
With \(\lambda\) as in Proposition~\ref{thm:source214gap}, on every matrix effect tuple, with \(\mathsf R\) as in \eqref{eq:214budget},
\begin{equation}\label{eq:214lowerD}
 \tr(3-3D^2)\ge2-21600\lambda c_{214}G_{214}^2\mathsf R
                  -A_{214}\sqrt{\mathsf R}.
\end{equation}
\end{lemma}
\begin{proof}
Set \(\epsilon=60\sqrt{\mathsf R}\), \(L=341110000\), \(V=23L=7845530000\), and \(\gamma=G_{214}\epsilon=200L\epsilon\). The graph approximation, signed-word, and meet estimates in Appendix~\ref{app:214proof} construct equal-rank orthogonal corners \(F_1,F_2\), a remaining corner \(C\), a coisometry \(A_0:C\to F_1\), and 37 exact symmetries \(S_j\) on \(F_1\) such that
\[
 t_c=\tr F_1\le\tfrac13,\qquad
 \tr Q_{66}(S)\le c_{214}\gamma^2.
\]
The latter follows by squaring the individual row bounds: ordinary rows cost \((h(a_r)+h(b_r))\gamma\), where \(h(w)=3(|w|-1)+\sum_{j\in w}\delta_j\), and a same-face meet costs \([4(\delta_a+\delta_b)+10]\gamma\). Their complete squared sum is \(c_{214}\); see \eqref{q214:mixedcapbound}. All norms here use the ambient trace.

For the signal the same interfaces give
\[
\begin{array}{ll}
 \|U_t-S_t\|_2\le3\gamma,&\|W_t-S_tA_0\|_2\le4\gamma,\\
 \|U_d-S_d\|_2\le929\gamma,&\|W_d-S_dA_0\|_2\le929\gamma.
\end{array}
\]
For \(d\), this uses its paid inverse-output alias \(U_dW_d\approx A_0\) and sign reconstruction, rather than an assertion about common initial spaces. Put
\[
 K_1=U_tW_d,\quad K_2=U_dW_t,\quad
 L_1=S_tS_dA_0,\quad L_2=S_dS_tA_0.
\]
Then \(\|K_1-L_1\|_2+\|K_2-L_2\|_2\le1865\gamma\). For any coisometry \(K:C\to F_1\) the full projection, including unused directions of \(C\), is
\[
 Q(K)=\frac12\begin{pmatrix}I&K\\K^*&2I-K^*K\end{pmatrix}.
\]
Direct multiplication shows \(Q(K)^2=Q(K)\), and
\[
 \|Q(K)-Q(L)\|_2^2
 =\tfrac14(2\|K-L\|_2^2+\|K^*K-L^*L\|_2^2)
 \le\tfrac32\|K-L\|_2^2.
\]
Write \(D_{\rm can}=2(Q(K_1)-Q(K_2))\) and \(D_{\rm id}=2(Q(L_1)-Q(L_2))\). Both have norm at most 2, and
\[
 |\tr(D_{\rm can}^2-D_{\rm id}^2)|\le16\cdot1865\gamma.
\]
Each virtual output, with the common coordinate term \(F_1+C\), is within \(V\epsilon\) of its canonical projection, by the spectral-output estimate in Appendix~\ref{app:214proof}. The common terms cancel from their difference. Thus \(\|D-D_{\rm can}\|_2\le4V\epsilon\). The coefficient absolute sum of \(D\) is 14, so
\[
 |\tr(D^2-D_{\rm can}^2)|\le(14+2)4V\epsilon=64V\epsilon.
\]
Since the two \(L_i\) have the same initial projection,
\(\tr D_{\rm id}^2=2\|S_tS_d-S_dS_t\|_2^2\).
The source identity \(\tr F_0=t_c-\|S_tS_d-S_dS_t\|_2^2\), and the weaker source margin \(1/2\), imply
\[
 \|S_tS_d-S_dS_t\|_2^2\le t_c/2+\lambda c_{214}\gamma^2.
\]
This is obtained by applying the source inequality with the normalized corner trace and then multiplying by \(t_c\); if \(t_c=0\), all supported operators are zero. Therefore
\[
 \tr(3-3D^2)\ge2-6\lambda c_{214}\gamma^2
                       -48\cdot1865\gamma-192V\epsilon.
\]
Now \(6\cdot60^2=21600\) and
\(60(48\cdot1865G_{214}+192V)=A_{214}\), which proves the claim.
\end{proof}
Adding \(\kappa\mathcal B\), with \(\kappa\) from \eqref{eq:214kappa}, gives
\begin{equation}\label{eq:214lowerPhi}
 \tr\Phi\ge2-A_{214}\sqrt{\mathsf R}
                  +\frac{A_{214}^2}{4}\mathsf R
 =1+\left(\frac{A_{214}\sqrt{\mathsf R}}2-1\right)^2\ge1.
\end{equation}

\begin{lemma}\label{lem:214elimination}
Suppose \eqref{eq:214lowerPhi} holds on every genuine effect tuple. Then, on arbitrary tuples of the 214 contraction variables,
\begin{equation}\label{eq:214eliminationbound}
 \tr\sigma(\Phi)\ge1-41616\kappa\mathsf H
                         -408\cdot31\kappa\sqrt{\mathsf H}.
\end{equation}
\end{lemma}
\begin{proof}
The anchor stability lemma, proved in Appendix~\ref{app:214proof}, gives orthogonal equal-rank projections \(B_1,B_2\), a graph projection \(Q^0\), and its companion \(K^0\), with
\[
 \|E_1-B_1\|_2\le3\sqrt{\mathsf H},\quad
 \|Q_0-Q^0\|_2\le3\sqrt{\mathsf H},\quad
 \|K-K^0\|_2\le6\sqrt{\mathsf H},\quad
 4Q^0+3B_1-5K^0=B_1+B_2.
\]
Consequently \(F_3=I-B_1-B_2\) is a genuine projection, and the eliminated effect expression obeys
\[
 E_3^{\rm vir}=I-4Q_0-3E_1+5K,\qquad
 \Delta=E_3^{\rm vir}-F_3,\qquad
 \|\Delta\|_2\le(12+9+30)\sqrt{\mathsf H}=51\sqrt{\mathsf H}.
\]
In the incidence formula \(E_3\) appears only in 31 support pairs and the affine row
\(\ell_0=2(I-E_3-4Q_0-3E_1+5K)=2(E_3^{\rm vir}-E_3)\).
The signal and all other terms are independent of \(E_3\). At the virtual value the row vanishes; at \(F_3\) it equals \(2\Delta\). Since its contribution to \(\mathcal B\) is \(4\ell_0^2\), deleting it loses exactly \(16\kappa\|\Delta\|_2^2\) in trace.

For a support pair with other effect \(Y_v\), its contribution to \(\mathcal B\) is \(4(E_3Y_v+Y_vE_3)\). Replacing \(F_3\) by \(E_3^{\rm vir}\) changes its trace by \(8\tr(\Delta Y_v)\), whose absolute value is at most \(8\|\Delta\|_2\), because \(\|Y_v\|_2\le1\). Evaluate \(\Phi\) first at the legitimate effect \(F_3\), use its lower bound 1, and sum these losses:
\[
 \tr\sigma(\Phi)\ge1-16\kappa(51)^2\mathsf H
                         -8\cdot31\kappa\cdot51\sqrt{\mathsf H}.
\]
This is \eqref{eq:214eliminationbound}. The argument never requires \(E_3^{\rm vir}\) to be an effect or a contraction.
\end{proof}
Finally \(\tr J_{\rm anc}=4\mathsf H\) cancels the linear loss in \eqref{eq:214eliminationbound}. With \(b=408\cdot31\kappa\),
\begin{equation}\label{eq:214finalsquare}
 \tr\Psi\ge1-b\sqrt{\mathsf H}+b^2\mathsf H
 =\tfrac34+(b\sqrt{\mathsf H}-\tfrac12)^2\ge\tfrac34.
\end{equation}
This proves the matrix assertion of Theorem~\ref{thm:quadratic214} using Proposition~\ref{thm:source214gap}.

At the exact group witness, use the three scalar coordinate projections \(E_1,E_2,E_3\in M_3(M)\). Graph projections of source unitaries \(U\) have the form
\begin{equation}\label{eq:214graph}
 Q_{ij}(U)=\tfrac12\bigl(E_i+E_j+t_{ij}(U)\bigr),
\end{equation}
where \(t_{ij}(U)\) has \(U\) in block \((i,j)\) and \(U^*\) in block \((j,i)\). The source relations make all affine rows, supports and spectral energies vanish. The quadratic signal satisfies, by \eqref{eq:214comm},
\[
 (\tr_3\otimes\tau)(D^2)
 =\tfrac23\tau\bigl((s_ts_d-s_ds_t)^*(s_ts_d-s_ds_t)\bigr)
 =\tfrac43.
\]
Thus \((\tr_3\otimes\tau)(\Phi)=3-3(4/3)=-1\). The affine elimination is exact at this witness and the anchor correction vanishes, proving the negative assertion of Theorem~\ref{thm:quadratic214}. The quadratic degree and occurrence of all variables are verified by the finite coefficient certificate; for example the coefficient of \(X_{\mathrm{product58},P}^2\) is \(-27\), independent of \(\kappa\).

\section{A two-variable separator of exact degree 24}\label{sec:degree24}
To encode the source in two self-adjoint variables, we use an exact finite-spectrum encoder and a cubic corner selector. Word degree accounts for the auxiliary operators realized in the selected corner.

\begin{theorem}\label{thm:degree24}
The polynomial \(f_{24}\), defined below, belongs to \(\ZZ\langle X,Y\rangle\), is Hermitian, and has exact ordinary and cyclic degree 24. For every pair of self-adjoint matrix contractions,
\begin{equation}\label{eq:24positive}
 \tr_Nf_{24}(A,B)\ge\tfrac34.
\end{equation}
There exists a prescribed pair of \emph{strict} self-adjoint contractions in a matrix algebra over a finite tracial algebra with trace \(-1\).
\end{theorem}
The rational rank assertions used below are verified by the finite certificates described in Appendix~\ref{app:24data} and \hyperref[sec:verification]{the data and verification statement}.

\subsection{Source reduction}
The initial characteristic-five source has 57 unitary coordinates and 170 relation residuals. A compact source uses 33 coordinates (21 symmetries and 12 genuine unitaries) and 78 paid relations. Signed unitary word rewrites and a quantitative transfer certificate compare the squared source defects with factor \(2^{56}\). Splitting relations of type \(3+1\) into \(2+2\) introduces 32 active helpers, giving 65 complete coordinates and 111 paid rows. The phase-compatible real-component description has 101 active real components, including two cheap skew-adjoint outputs, with source-word length at most four.

The chosen source coefficient is the completely specified tower
\begin{equation}\label{eq:mu5}
 J_5=2^{2^{1122100}},\qquad
 \delta_5=2^{-J_5},\qquad
 \mu_5=2^{2J_5+30}.
\end{equation}
It is the coefficient for the characteristic-five source argument. The source inequality for this precise corrected polynomial is proved in Proposition~\ref{hyp:source24}.

After incorporating coercivity corrections, the compact source polynomial is $\mathcal G=P_0+\mu_5P_1$. In the phase-twisted real-component presentation, $P_0=-q-VZ^*VZ-Z^*VZV$, with $Z=z_{1,0}$ and the two external flip coordinates as defined in the source construction. To define the coercive correction, first undo the phase substitution and use the untwisted complete coordinates. Write $I_a=\tr((Z_a^*Z_a-q)^2)$ and $I=\sum_a I_a$. For a starred word $W=T_1\cdots T_m$, define $L_W=\sum_i(T_i^*T_i-q)$, counting occurrences. For comparison, let $a_i$ be the underlying coordinate of $T_i$ and put $L_W^{\rm can}=\sum_i(Z_{a_i}^*Z_{a_i}-q)$. The authenticated coefficient dictionary uses the latter canonical convention. For a starred occurrence these two defects have the same trace by cyclicity, but they are not literally the same polynomial. The gate correction for a recovered coordinate $c=W_c(Z)$ is
\[
g_c=Z_c^*Z_c+q-Z_c^*W_c-W_c^*Z_c+L_{W_c}.
\]
The seed and gate defects are combined by the reverse acyclic recovery order; for a cap $u-v$ put $\operatorname{cap}(u,v)=2q-u^*v-v^*u+L_u+L_v$. After forming the entire correction in untwisted coordinates, apply the prescribed phase substitution to it. If $R_{\rm corr}$ is the sum of the 111 corrected caps and $E$ is the resulting coercivity budget, set $P=R_{\rm corr}+K_{\rm corr}E$, with $K_{\rm corr}=30+84T^2$, and
\[
\lambda(T)=2C+64T^6,\qquad C=2^{57},\qquad P_1=\lambda(T)P.
\]
Let $P^{\rm can}$ and $\mathcal G^{\rm can}=P_0+\mu_5\lambda(T)P^{\rm can}$ be obtained by using $L_W^{\rm can}$ in every cap and gate, with the same recursion, phase convention, and coefficients. Then $\mathcal G-\mathcal G^{\rm can}$ is a scalar linear combination of defects $Z_aZ_a^*-Z_a^*Z_a$, hence is a sum of commutators and has trace zero. The polynomial below is defined using $\mathcal G$ with the occurrence-wise convention $L_W$. The 1,025-term coefficient dictionary and the source proof use $P^{\rm can}$.

For a component word $w$ in the source polynomial $P_j$, let $c_{j,w}(T)$ be the coefficient polynomial and $c_{j,w}^{\rm abs}(T)=\sum_k|c_k|T^k$. The selector-error majorants are
\begin{equation}\label{eq:24majorant}
B_j=2\sum_{w\in P_j}\max(1,|w|)(4T)^{|w|}c_{j,w}^{\rm abs}(T),\quad B=B_0+\mu_5B_1,
\qquad B_0=2+16384T^4,\quad B_1=\lambda(T)B_P.
\end{equation}
Here $B_P$ is the same coefficientwise majorant summed over the 1,025 collected terms of the occurrence-wise penalty $P$; explicitly it is obtained from the right-hand side of \eqref{eq:24majorant} by replacing $P_j$ with the term list of $P$ and omitting the leading factor $\lambda(T)$. The canonical penalty also has 1,025 terms, but its literal coefficient list and its majorant need not agree with those of $P$. For the occurrence-wise coefficient data, exact verification of coefficient signs shows that multiplication by $\lambda(T)=2C+64T^6$ introduces no opposite-sign coefficient overlaps, so $B_1=\lambda(T)B_P$ is the direct $P_1$ majorant here.
\begin{proposition}\label{hyp:source24}
Under the finite source certificates specified in Appendix~\ref{app:fiveproof}, for every finite matrix corner with identity \(q\), every tuple of the 101 real components specified in Appendix~\ref{app:24data}, of their declared self-adjoint or skew-adjoint types and component norms at most the fixed \(T\) of the compiler, satisfies
\begin{equation}\label{eq:24cornergap}
 \tr\mathcal G\ge\tfrac12\tr(q).
\end{equation}
There is also a prescribed typed component tuple in \(M_2(A)\), for a finite tracial algebra \(A\), with every component of norm at most one, for which the normalized trace of \(\mathcal G\) is \(-1\). Its two cheap skew components are \(\mathsf J\otimes S_1\) and \(\mathsf J\otimes S_2\), where
\(\mathsf J=\left(\begin{smallmatrix}0&1\\-1&0\end{smallmatrix}\right)\)
and \(S_1,S_2\) are commuting symmetries in \(A\).
\end{proposition}
\begin{proof}
The native 57 gap and the compact, split, and component-domain transfers in Appendix~\ref{app:fiveproof} prove the assertion for $\mathcal G^{\rm can}$. The trace identity above transfers it to $\mathcal G$ with the same coefficients and $T$. The separate faithful-amalgam construction in Section~\ref{fivemerge:sec:witness} proves the witness and its two commuting cheap symmetry factors. The analytic inputs are described in Section~\ref{fivemerge:sec:analytic}. The integer spectral estimate follows from Wang--Zhi~\cite[Lemma~3.3 and equation~(3.12)]{WangZhi2026} (version~1), applied through the exact residual identities in Section~\ref{fivemerge:sec:integer-adaptation}.
\end{proof}

\subsection{The exact two-variable encoder}
Set
\begin{equation}\label{eq:selector24}
 e=e(X,Y)=X^2Y+2XYX+YX^2,\qquad s=e^2-e.
\end{equation}
Both \(e\) and \(s\) are Hermitian. The 101 active decoders are selected from a common collection of 110 even, 36 ordinary odd, and two cheap odd channels. The complete moment, differential, and right-inverse construction is proved in Appendix~\ref{app:encoderproof}.

\begin{proposition}\label{thm:encoder24}
Let \(A\) be a unital complex \(C^*\)-algebra, let \(S_1,S_2\in A\) be commuting symmetries, and put \(\mathcal B=M_2(A)\). Assume the sequential rational-lift identities and the six specified nonzero leading minors in Appendix~\ref{app:24data}. The fixed integer decoders constructed in Appendix~\ref{app:encoderproof} realize any prescribed self-adjoint contractions in the 110 even channels and any prescribed skew-adjoint contractions in the 36 ordinary odd channels, simultaneously, by strict self-adjoint contractions \(X,Y\in M_{800}(\mathcal B)\). They satisfy
\begin{equation}\label{eq:24exactencoder}
 e(X,Y)=q,\qquad q\sim e_{11}\otimes1_{\mathcal B},
\end{equation}
and the two cheap decoded outputs are exactly \(\mathsf J\otimes S_1\) and \(\mathsf J\otimes S_2\). The bounds are
\[
 \tfrac57 I<X\le\tfrac67 I,\qquad 0\le Y,\qquad\|Y\|<49/100.
\]
All decoder coefficients are fixed by finite rational operations and denominator clearing; no unspecified small parameter is used.
\end{proposition}
\begin{proof}
We summarize the estimates from Appendix~\ref{app:encoderproof}. The real Banach parameter space is \(\mathcal B_{\sa}^{240}\), with the maximum norm. Restoring the last coordinate by the zero-sum rule gives the weights
\[
 E_i=N^{-1}I+\varepsilon d_i(D_2/2\otimes I)
       +\varepsilon f_{i1}(E_2\otimes S_1)
       +\varepsilon f_{i2}(E_2\otimes S_2)+\sum_lH_{il}C_l,
\]
where \(D_2=\operatorname{diag}(1,-1)\), \(E_2=\left(\begin{smallmatrix}0&1\\1&0\end{smallmatrix}\right)\). The output space is \(\mathcal B_{\sa}^{110}\times\mathcal B_{\rm skew}^{36}\), again with the maximum norm. It is defined by the even and centered odd parts of the finite moment polynomials
\begin{equation}\label{eq:finite-moment24}
 m_w(E)=\sum_{i_0,\ldots,i_r}
 \frac{\prod_{j=0}^r\lambda_{i_j}^{b_j}}
 {\prod_{j=0}^{r-1}(\lambda_{i_j}+\lambda_{i_{j+1}})^2}
 E_{i_0}\cdots E_{i_r},\qquad
 \lambda_i=\frac{5N+1+i}{7N},\quad N=800.
\end{equation}
The isotropic identities make the two cheap moments exactly \(\varepsilon^2\mathsf J\otimes S_k\), independent of the controls. Even centering gives a reference \(C^\flat=O(\varepsilon^2)\); subtracting its rational odd Taylor coefficients through order four leaves output residual at most \(C_{\rm rem}\varepsilon^6\).

The six leading minors give a bounded right inverse \(R_\varepsilon\) of the full operator-valued differential at that reference, with
\[
 \|R_\varepsilon\|\le K_{\rm inv}\varepsilon^{-2}.
\]
This is an operator-norm statement for every \(A\): scalar inverse matrices act entrywise in the at most sixteen joint Peirce blocks of \(S_1,S_2\). The proof and assembly factors are given in Appendix~\ref{app:encoderproof}. Set
\[
 M=H_2K_{\rm inv}^2\varepsilon^{-4},\qquad
 r=\frac1{8M}=\frac{\varepsilon^4}{8H_2K_{\rm inv}^2}.
\]
In right-inverse output coordinates \(u\), the centered moment map \(\mathcal F\) has derivative the identity, Hessian norm at most \(M\), and \(\|\mathcal F(0)\|\le r/4\). Hence its nonlinear remainder \(\mathcal N(u)=\mathcal F(u)-\mathcal F(0)-u\) has Lipschitz constant at most \(Mr=1/8\) on the closed radius-\(r\) ball. For every target \(t\) with \(\|t\|\le r/4\),
\[
 u\longmapsto t-\mathcal F(0)-\mathcal N(u)
\]
maps that ball into itself, with image norm at most \(r/4+r/4+r/8<r\), and is a contraction. Its fixed point realizes all ordinary outputs exactly.

The chosen \(\varepsilon=Q_0^{-1}\) and the control bounds in Appendix~\ref{app:encoderproof} keep \(\max_i\|E_i-N^{-1}I\|<1/N\); therefore every weight remains strictly positive and their sum is \(I\). The ordinary decoder scale
\(A_{\rm ord}=32H_2K_{\rm inv}^2D_*^{4096}\varepsilon^{-4}\)
satisfies \(A_{\rm ord}^{-1}\le r/4\). Dividing any prescribed contraction by this scale puts every component simultaneously in the admissible target ball. Unused channels may be assigned zero. The cheap scale \(\varepsilon^{-2}\) gives exactly the two prescribed anchor outputs.

Finally let \(V_i=E_i^{1/2}\), \(D=\operatorname{diag}(\lambda_i)\), and
\[
 Y'_{ij}=\frac{E_i^{1/2}E_j^{1/2}}{(\lambda_i+\lambda_j)^2}.
\]
Then \(e(D,Y')=VV^*\), and
\(Y'=\int_0^\infty t e^{-tD}VV^*e^{-tD}\,dt\)
is positive with norm at most \((4\lambda_{\min}^2)^{-1}<49/100\). A polar unitary carries the constant reference isometry to \(V\), since the weights are in the specified small neighborhood. Conjugation by it gives \eqref{eq:24exactencoder}; compression gives exactly \eqref{eq:finite-moment24}. This proves both the realization and the strict contraction bounds. If \(A\) is finite tracial, the scalar corner has ambient trace \(1/800\).
\end{proof}

\subsection{The compiled polynomial}
Choose the integer \(T=1+\max_z\|B_z\|_{\coeff,1}\), where the maximum is over the 65 complete source coordinates listed in Appendix~\ref{app:24data}. Here $B_z=B_{z,S}+B_{z,K}$, with a missing component taken as zero; these component decoders are respectively self-adjoint and skew-adjoint. Thus $B_{z,S}=(B_z+B_z^*)/2$ and $B_{z,K}=(B_z-B_z^*)/2$. Since the coefficient $\ell^1$ norm is adjoint-invariant, each component-decoder norm is at most $\|B_z\|_{\coeff,1}<T$. The authenticated assignment has disjoint full possible supports for the two centered component decoders at each of the 36 paired coordinates. Thus the complete coefficient norm in the definition of $T$ is the sum of the component norms. Independently, the source and coefficient identities are checked symbolically in $T$ and the component-domain estimate holds for every $T\ge4$. For a source component word \(w=v_1\cdots v_m\), define the \emph{literal cyclic compiler}
\begin{equation}\label{eq:24compiler}
 \mathscr C(1)=e,\qquad
 \mathscr C(v_1\cdots v_m)=B_{v_1}eB_{v_2}e\cdots B_{v_m}e.
\end{equation}
This inserts a corner selector after each output. Extend $\mathscr C$ linearly to $\mathcal G$ and the canonical $\mathcal G^{\rm can}$. The difference $\mathscr C(\mathcal G)-\mathscr C(\mathcal G^{\rm can})$ is a scalar linear combination of commutators $[B_z e,B_z^*e]$; this holds for every $X,Y$, without assuming $e$ is a projection. Consequently their traces agree, and the fully symmetrized compiled difference is also a sum of commutators. Each component decoder has weighted degree at most 11 in Appendix~\ref{app:24data}, so this difference has ordinary degree at most 22. Thus the occurrence-wise and canonical compiled source terms have the same degree-24 homogeneous part and cyclic class. If their selector coefficients differ, the resulting two full polynomials need not have the same cyclic class; their difference is a sum of commutators and a multiple of $(e^2-e)^2$, which has degree at most 12. In particular, the fixed polynomial in \eqref{eq:f24} has the same degree-24 homogeneous part and the same class modulo commutators and the selector ideal as its canonical counterpart.

The positive integer $B=B_0+\mu_5B_1$ from \eqref{eq:24majorant} is computed from the occurrence-wise dictionary. It is independent of matrix size and bounds selector-replacement errors for $\mathscr C(\mathcal G)$.

The final two-variable polynomial is
\begin{equation}\label{eq:f24}
 \boxed{\ f_{24}=1+800\bigl(\mathscr C(\mathcal G)+
        \mathscr C(\mathcal G)^*\bigr)
       +(1600B)^2(e^2-e)^2.\ }
\end{equation}
The large scalar penalties are intrinsic to the construction: no projection relations are imposed on \(X,Y\) when evaluating \eqref{eq:f24} at arbitrary finite contractions.

To prove \eqref{eq:24positive}, let \(q'=1_{[1/2,\infty)}(e)\) and \(h=\|e^2-e\|_2\). Lemma~\ref{lem:projectionrounding} gives \(\|e-q'\|_2\le2h\). Telescoping every compiled word of the specified polynomial, with all decoders fixed, gives a total tracial error at most \(Bh\). When each selector is replaced by \(q'\), cyclicity identifies the resulting expression with the source polynomial evaluated in the finite corner \(q'M_Nq'\), whose coordinates have the required types and norms. From \eqref{eq:24cornergap}, the source trace is nonnegative. Hence
\[
 \operatorname{Re}\tr_N\mathscr C(\mathcal G)\ge-Bh.
\]
Substituting into \eqref{eq:f24} and completing a scalar square yields
\begin{equation}\label{eq:24completesquare}
 \tr_N f_{24}\ge1-1600Bh+(1600Bh)^2
 =\tfrac34+(1600Bh-\tfrac12)^2\ge\tfrac34.
\end{equation}
At the exact encoder witness \(e=q\), so the global selector penalty vanishes. Using the witness clause of Proposition~\ref{hyp:source24}, Proposition~\ref{thm:encoder24} realizes all source components in the same corner. Its source trace is \(-1\) with the normalized corner trace, whence
\begin{equation}\label{eq:24negative}
 (\tr_{800}\otimes\tr_2\otimes\tau)\mathscr C(\mathcal G)=-\tfrac1{800},
 \qquad (\tr_{800}\otimes\tr_2\otimes\tau)(f_{24})=-1.
\end{equation}
The exact encoder witness satisfies the strict operator bounds \(\frac57 I<X\le\frac67I\) and \(0\le Y\), \(\|Y\|<49/100\).

\subsection{Exact degree}
The compiler has maximum weighted degree 24. Exact cancellation at that degree is excluded by a selected cyclic word class,
\begin{equation}\label{eq:24word}
 w_{24}=(XXYYYY)^4,
\end{equation}
whose cyclic coefficient is recorded in
\path{degree24/top_cyclic_degree_certificate.json} and reconstructed
by the exact checker \path{degree24/check_top_cyclic_degree.py}:
\begin{equation}\label{eq:24topcoefficient}
 3200\mu_5(2\cdot2^{57}+64T^6)A_{\rm ord}^4q(T)>0.
\end{equation}
Here \(A_{\rm ord}\) is the positive integer decoder scale defined in Appendix~\ref{app:24data}, and
\begin{align*}
 q(T)={}&1050+4680T+17160T^2+47904T^3+106776T^4+204000T^5\\
 &+331488T^6+459648T^7+551520T^8+562944T^9\\
 &+478656T^{10}+334848T^{11}+181248T^{12}+64512T^{13}+10752T^{14}.
\end{align*}
All coefficients of \(q(T)\) are positive. The cited coefficient certificate attributes this class to a single quartic seed-defect type; competing maximum-degree terms contain a longer run of \(Y\)'s and cannot contribute. The degree-12 selector penalty cannot cancel a cyclic degree-24 term. This proves \(\deg f_{24}=\deg_{\cyc}f_{24}=24\) without evaluating the enormous constants.

\section{Consequences, comparisons, and remaining questions}\label{sec:outlook}
The four results put degree and arity into a concrete common framework. Define
\[
 \mathfrak S_{n,d}=\left\{f=f^*\in\RR\langle X_1,\ldots,X_n\rangle:
 \deg f\le d,\ m_{\rm tr}(f)<m_{\rm mat}(f)\right\}.
\]
Padding with unused variables and allowing higher degrees make these sets monotone in both parameters. Theorem~\ref{thm:quadraticrigidity} gives \(\mathfrak S_{3,2}=\varnothing\), while Theorem~\ref{thm:quadratic214} gives \(\mathfrak S_{214,2}\ne\varnothing\). Theorem~\ref{thm:degree6} gives \(\mathfrak S_{10,6}\ne\varnothing\), and Theorem~\ref{thm:degree24} gives \(\mathfrak S_{2,24}\ne\varnothing\). Lemma~\ref{lem:one-variable} gives \(\mathfrak S_{1,d}=\varnothing\) for every \(d\), and Corollary~\ref{cor:cubic} gives \(\mathfrak S_{2,3}=\varnothing\). Together these yield \eqref{eq:arity-bounds}.

Degree and arity do not measure all aspects of complexity. The ten-variable polynomial has short support and a large fixed coefficient; the quadratic construction has 2,416 terms and an algorithmically defined coefficient; the two-variable construction has a moment encoder and an implicitly specified integer coefficient table. Smaller \((n,d)\) need not yield smaller certificates or better numerical conditioning.

Proposition~\ref{prop:moduleobstruction} applies to the quadratic and two-variable constructions with \(\eta=1\), so
\begin{equation}\label{eq:universalmodule}
 f+\varepsilon\notin\calM_n+[\calA_n,\calA_n]_{\rm lin}
 \qquad(0\le\varepsilon<1).
\end{equation}
This excludes contraction-module certificates modulo commutators. By contrast, Corollary~\ref{cor:quadraticmodule} supplies such certificates for uniformly matrix trace-positive quadratics in at most three variables, with the shift specified there.

The following questions remain:
\begin{enumerate}[label=(\roman*),leftmargin=*,itemsep=2pt]
 \item Determine whether a degree-two tracial separator exists in four self-adjoint variables, and improve the upper bound 214.
 \item Determine the minimal degree of a two-variable separator, between the lower bound 4 from Corollary~\ref{cor:cubic} and the upper bound 24.
 \item Find geometric encoders that simultaneously improve degree and arity.
 \item Obtain practical, verifiable coefficient dictionaries of substantially lower arithmetic height, or degree-bounded proofs of non-membership for additional algebraic quotients.
\end{enumerate}

\section*{Data and exact verification}\label{sec:verification}
External data\footnote{\url{https://wangjie212.github.io/jiewang/codes/Connes.zip}} are cited where they specify source relations, incidence data, rewrite paths, coefficient dictionaries, or rational encoder arrays not printed in full. File paths are relative to the accompanying data directory unless stated otherwise.

For the degree-six construction, exact verification reconstructs the polynomial from the chart and routes in Appendix~\ref{app:degree6}, checks the word count, and certifies nonvanishing of the specialized coefficients. For the quadratic construction, it compares \path{quadratic214/base_manifest.json} with the collected coefficient array and checks the finite source reductions and geometric constants used in Appendices~\ref{app:quadraticsource} and~\ref{app:214proof}. For the two-variable construction, it checks the source transfers, compiler identities, rational lift, leading minors, and Taylor centers used in Appendices~\ref{app:fiveproof} and~\ref{app:encoderproof}.

\section*{Acknowledgements}
The author acknowledges the use of GPT-6 Astra to assist with
brainstorming, mathematical development, and manuscript drafting.
The author is solely responsible for the final content, analysis, and
conclusions.

\section*{Funding}
The author was supported by the National Key R\&D Program of China (No.~2023YFA1009401) and the National Natural Science Foundation of China (No.~12171324).

\clearpage
\appendix
\section{Finite data for the ten-variable graph construction}\label{app:degree6}
The following data define \(P_6\). Source-letter adjoints are denoted by \(^*\); \(a_+\) is a separate unitary coordinate.

\subsection{The 27 weighted source equations}\label{app:relations}
\begin{longtable}{@{}r>{\(}l<{\)}r@{\qquad}r>{\(}l<{\)}r@{}}
\toprule
No.&\multicolumn{1}{c}{Equation}&\(\lambda\)&No.&\multicolumn{1}{c}{Equation}&\(\lambda\)\\\midrule
\endfirsthead
\toprule No.&\multicolumn{1}{c}{Equation}&\(\lambda\)&No.&\multicolumn{1}{c}{Equation}&\(\lambda\)\\\midrule
\endhead
1&r=tt&1&15&ax=xa&1\\
2&v=tr&1&16&zz_2=z_2z&1\\
3&qt=t^*q&1&17&dx=xd&1\\
4&qz=z^*q&1&18&dz_2=z_2d&1\\
5&bq=qa&1&19&ab=ba&1\\
6&a_+r=ra&9&20&ba_+=a_+b&1\\
7&z_2r=rz&1&21&aa_+=ba_+a&1\\
8&\ell=rq&1&22&az=za&1\\
9&xzx^*=dz&1&23&a_+z=za_+&2\\
10&e_2\ell=\ell d^*&1&24&Wa=av&1\\
11&d^*e_2d=pe_2&1&25&xs=sz&1\\
12&az_2=pa&1&26&jq=qs&1\\
13&vz=zv&1&27&js=sq&1\\
14&tx=xt&1&&&\\
\bottomrule
\end{longtable}
There are 27 weighted relations with total weight 36. Adding the order-five cap with weight 512 gives 28 caps and total weight 548. This bookkeeping produces the factor \(2192=4\cdot548\) in \(c=3+2192\Omega\).

\subsection{The full labeled chart}\label{app:chart}
The table lists every nonanchor edge \((i,j)\) with \(1\le i<j\le9\). Every anchor edge \((0,i)\) reads 1 in either direction. A reverse directed edge reads the adjoint of the listed label.
\begin{center}
\small
\begin{tabular}{@{}r l@{\qquad}r l@{\qquad}r l@{}}
\toprule
Edge&Label&Edge&Label&Edge&Label\\\midrule
(1,2)&\(1\)&(2,7)&\(a\)&(4,8)&\(z\)\\
(1,3)&\(a^*\)&(2,8)&\(t\)&(4,9)&\(z_2\)\\
(1,4)&\(r\)&(2,9)&\(d^*\)&(5,6)&\(1\)\\
(1,5)&\(a_+\)&(3,4)&\(1\)&(5,7)&\(r\)\\
(1,6)&\(z_2^*\)&(3,5)&\(p^*\)&(5,8)&\(q\)\\
(1,7)&\(a_+^*\)&(3,6)&\(d\)&(5,9)&\(a\)\\
(1,8)&\(b^*\)&(3,7)&\(z\)&(6,7)&\(b\)\\
(1,9)&\(z\)&(3,8)&\(x\)&(6,8)&\(a^*\)\\
(2,3)&\(e_2^*\)&(3,9)&\(e_2\)&(6,9)&\(a_+\)\\
(2,4)&\(v\)&(4,5)&\(s\)&(7,8)&\(z^*\)\\
(2,5)&\(t^*\)&(4,6)&\(t^*\)&(7,9)&\(x^*\)\\
(2,6)&\(x^*\)&(4,7)&\(z^*\)&(8,9)&\(x^*\)\\
\bottomrule
\end{tabular}
\end{center}
The real labels at loops \(1,\ldots,9\) are, in order,
\[
 (q,W,\ell,q,j,W,v,s,\ell).
\]
There are 33 labeled nonanchor edges and three constant nonanchor edges: \((1,2)\), \((3,4)\), and \((5,6)\). The eleven complex reference edges are oriented to read the unstarred letter, and the six real reference labels use their lowest-index loop:
\begin{center}
\small
\begin{tabular}{@{}llll@{\qquad}llll@{}}
\toprule
Letter&Reference&Letter&Reference&Letter&Reference&Letter&Reference\\\midrule
\(x\)&\((6,2)\)&\(t\)&\((5,2)\)&\(z\)&\((1,9)\)&\(r\)&\((1,4)\)\\
\(z_2\)&\((6,1)\)&\(q\)&loop 1&\(a\)&\((3,1)\)&\(v\)&loop 7\\
\(b\)&\((8,1)\)&\(\ell\)&loop 3&\(a_+\)&\((1,5)\)&\(s\)&loop 8\\
\(d\)&\((9,2)\)&\(W\)&loop 2&\(e_2\)&\((3,2)\)&\(j\)&loop 5\\
\(p\)&\((5,3)\)&&&&&&\\
\bottomrule
\end{tabular}
\end{center}
Repeated uses of a real label on an edge are compared with its reference loop. This convention, together with the full edge table, determines the 25 nonreference copies in \(D_{\rm raw}\).

\subsection{Closed walks for the 29 source classes}\label{app:routes}
Each list gives the vertices of a closed walk, with the final transition returning to the first vertex. The repeated endpoint is not an additional variable factor. Cyclic rotation and adjoint identify the source class for route selection; the expanded polynomial is paired with its literal adjoint as in \eqref{eq:R6}.
\begingroup\small\setlength{\LTpre}{0pt}\setlength{\LTpost}{0pt}
\begin{longtable}{@{}r>{\(}l<{\)}>{\(}l<{\)}@{\qquad}r>{\(}l<{\)}>{\(}l<{\)}@{}}
\toprule
No.&\multicolumn{1}{c}{Source class}&\multicolumn{1}{c}{Walk}&No.&\multicolumn{1}{c}{Source class}&\multicolumn{1}{c}{Walk}\\\midrule
\endfirsthead
\toprule No.&\multicolumn{1}{c}{Source class}&\multicolumn{1}{c}{Walk}&No.&\multicolumn{1}{c}{Source class}&\multicolumn{1}{c}{Walk}\\\midrule
\endhead
1&\ell qr^*&(3,3,4,4,1,0)&16&dz_2d^*z_2^*&(3,6,1,2,9,4)\\
2&rt^*t^*&(5,7,0,8,2)&17&jqsq&(5,5,8,8)\\
3&rvt&(1,4,2,8,0)&18&jsqs&(5,5,4,4)\\
4&Wava^*&(2,2,7,7)&19&qtqt&(8,5,2,1,1,2)\\
5&Wz_2Wz_2^*&(6,6,1,2,2,1)&20&qzqz&(8,5,0,7,4,4)\\
6&aba^*b^*&(8,6,7,2,1)&21&rzr^*z_2^*&(5,7,4,1,6)\\
7&aqb^*q&(3,1,1,8,5,0)&22&sxs z^*&(4,5,0,3,8,8)\\
8&ar^*a_+^*r&(2,7,5,1,4,0)&23&txt^*x^*&(6,4,3,8,2)\\
9&axa^*x^*&(5,9,7,2,6)&24&vzvz^*&(7,7,4,0,2,4)\\
10&aza^*z^*&(3,1,9,5,0,7)&25&zz_2z^*z_2^*&(7,4,9,1,6,0)\\
11&az_2a^*p^*&(5,9,0,6,1,3)&26&aa_+a^*a_+^*b^*&(8,6,9,5,1)\\
12&a_+ba_+^*b^*&(6,9,0,8,1,7)&27&de_2^*p^*d^*e_2&(9,2,3,5,6,3)\\
13&a_+za_+^*z^*&(7,1,9,6,0,4)&28&dzxz^*x^*&(9,2,1,9,7,8)\\
14&d\ell e_2\ell&(9,2,0,3,3,9)&29&zzzzz&(3,7,4,8,7,4)\\
15&dxd^*x^*&(3,6,2,9,0,8)&&&\\
\bottomrule
\end{longtable}
\endgroup
The signal is the class in row 5, cyclic-adjoint equivalent to \(\omega=z_2^*Wz_2W\); the other 28 rows encode the 27 relation caps and \(z^5\).

\subsection{The quartic source dictionary}\label{app:quarticdictionary}
The following 31 equations record the dictionary of the quartic presentation in \cite[Appendix~A]{WangZhi2026}, used in the reduction in Section~\ref{sec:degree6}. Here $C_d,C_p,C_a,z_3$ are helper coordinates.
\begin{longtable}{@{}rll@{\qquad}rll@{}}
\toprule No.&Left&Right&No.&Left&Right\\\midrule
\endfirsthead
\toprule No.&Left&Right&No.&Left&Right\\\midrule
\endhead
1&$r$&$tt$&17&$tx$&$xt$\\
2&$v$&$tr$&18&$ax$&$xa$\\
3&$qt$&$t^*q$&19&$zz_2$&$z_2z$\\
4&$qz$&$z^*q$&20&$dx$&$xd$\\
5&$bq$&$qa$&21&$dz_2$&$z_2d$\\
6&$a_+r$&$ra$&22&$ab$&$ba$\\
7&$z_2r$&$rz$&23&$ba_+$&$a_+b$\\
8&$rz_3$&$zr$&24&$C_a$&$ba_+$\\
9&$\ell$&$rq$&25&$aa_+$&$C_a a$\\
10&$C_d$&$xz$&26&$az$&$za$\\
11&$C_dx^*$&$dz$&27&$az_3$&$z_3a$\\
12&$e_2\ell$&$\ell d^*$&28&$Wa$&$av$\\
13&$C_p$&$d^*e_2$&29&$xs$&$sz$\\
14&$C_pd$&$pe_2$&30&$jq$&$qs$\\
15&$az_2$&$pa$&31&$js$&$sq$\\
16&$vz$&$zv$&&&\\
\bottomrule
\end{longtable}
A witness for this quartic dictionary is built from elementary matrices over $\mathbb F_5[t_1^{\pm1},t_2^{\pm1},t_3^{\pm1}]$: $x=e_{12}(1)$, $y=e_{23}(1)$, $z_j=e_{31}(t_j)$, $d_j=e_{32}(-t_j)$, $e_j=e_{21}(t_j)$, $c=e_{13}(1)$ and $p_{ab}=e_{31}(t_at_b)$, together with the integer-shear actions and the external flip and inverse-transpose symmetries. The element $a^{-1}z_2a=e_{31}(t_1^{-1}t_2)$ is outside the amalgamating subgroup, so the two signal words are nonidentity reduced words and have trace zero.

\section{Characteristic-two source and quadratic construction data}\label{app:quadraticsource}
\subsection{The 37-symmetry source data}\label{app:214data}
The reduced source has 37 symmetry coordinates. Throughout the source reductions, retained rows and coordinates are those included in the displayed reduced presentation. The prescribed witness is in the tracial group algebra of
\[
\Gamma=(\widehat G *_H\widehat G)\rtimes(C_{2,\mathrm{swap}}\times C_{2,\theta}),\quad
\widehat G=\mathrm{GL}_3(\mathbb F_2[t_1^{\pm1},t_2^{\pm1},t_3^{\pm1}])\rtimes\mathrm{GL}_3(\mathbb Z),\quad
H=\mathrm{EL}_3(\mathbb F_2[t_1,t_2,t_3]).
\]
Here $\theta$ is inverse transpose and swap interchanges the two amalgam factors while fixing $H$. With $B=W_{12}R_{12}S_2$, the flip is $V=B\,\mathrm{swap}\,B^{-1}$, and $B^{-1}z_{10}B=e_{31}(t_1^{-1}t_2)\notin H$. Thus the two nonconstant objective words are nonidentity reduced amalgam words and have canonical group trace zero. Every retained source defect vanishes at the prescribed tuple and $\tau(F_0)=-1$.

\subsubsection{All retained source cuts}
The table records the 66 retained rows using the zero-based row labels of the 80-row presentation. A displayed pair $u,v$ means the relation $u=v$; every source coordinate is a symmetry, so reversing a word gives its adjoint.
\small
\begin{longtable}{@{}rll@{\qquad}rll@{}}
\toprule Row&Left&Right&Row&Left&Right\\\midrule
\endfirsthead
\toprule Row&Left&Right&Row&Left&Right\\\midrule
\endhead
1&ed&de&44&FI&GFGF\\
3&hn&lh&45&Ju&vJ\\
4&mh&hm&46&Jb&cJ\\
5&kr&gk&47&HJ&Jz\\
6&ek&ke&48&JA&t\\
8&ha&ah&49&HF&FA\\
9&eh&he&50&uc&Ku\\
11&et&te&51&z&KB\\
12&gt&tg&52&Fo&z\\
13&ft&tf&53&aG&q\\
14&aw&wa&54&qF&Fp\\
15&vx&xv&55&Dp&lD\\
16&yw&wy&56&lG&Gk\\
17&zh&hz&58&Fm&lF\\
18&Ah&hA&59&plD&ly\\
21&hq&qh&60&wz&zx\\
22&oi&io&61&Fy&wF\\
23&Ap&pA&62&Dc&eD\\
24&oA&Ao&63&Ey&yk\\
25&yx&xy&64&mk&sm\\
28&zd&dz&65&Fd&eF\\
29&zk&kz&66&qf&eq\\
30&Ak&kA&67&Gn&EG\\
31&ak&ka&68&Dd&gD\\
32&lh&hili&69&hk&sn\\
33&AB&BA&70&pi&hp\\
34&kB&Bk&71&hr&dh\\
36&tB&Bt&73&IG&LI\\
37&pB&Bp&75&Lu&GLG\\
39&yh&hy&76&vu&Nv\\
40&wxhwx&h&77&au&N\\
41&cwx&wxc&78&ab&LG\\
42&uwx&wxu&79&IG&Nz\\
\bottomrule
\end{longtable}
\normalsize
Rows $0,7,10,19,20,26,27,35,38,43$ are controlled by acyclic rewrite paths; rows $57,74$ vanish after palindromic substitutions; rows $2,72$ are controlled by the exact Gram estimate $\tr Q_{68}\le 2845411\tr Q_{66}$. Using the 80-row source inequality proved in Section~\ref{app:char2norm}, these reductions imply \eqref{eq:source214gap}. The source-coordinate relabeling and its final error allowance (in units of the common error scale $\gamma$) are recorded in Table~\ref{tab:214lettermap}. The calibration forest is given in Table~\ref{tab:214forestdata}; its edges are used in the displayed parent-before-child order.

\begin{table}[ht]
\centering\scriptsize
\caption{Source coordinates, final letters, and error allowances.}\label{tab:214lettermap}
\begin{tabular}{@{}lll@{\qquad}lll@{}}
\toprule Coordinate&Letter&$\delta_j$&Coordinate&Letter&$\delta_j$\\\midrule
$x_0$&$a$&608&$c_0$&$u$&4\\
$y_0$&$b$&884&$d_{00}$&$v$&19\\
$z_{00}$&$c$&899&$B_1$&$w$&4\\
$z_{10}$&$d$&928&$C_1$&$x$&18\\
$z_{20}$&$e$&913&$C_2$&$y$&19\\
$z_{30}$&$f$&1131&$F_{\rm graph}$&$z$&4\\
$p_{120}$&$g$&942&$G_{\rm graph}$&$A$&32\\
$S_1$&$h$&10915614&$\mathrm{Tau}$&$B$&4\\
$S_3$&$i$&10915804&$J_{\rm hybrid}$&$D$&4\\
$R_{21}$&$k$&63&$T_{21}$&$E$&18\\
$R_{31}$&$l$&77&$F_2A$&$F$&4\\
$R_{32}$&$m$&92&$F_2B$&$G$&4\\
$T_{31}$&$n$&4&$F_2E_2$&$H$&18\\
$W_{12}$&$o$&4&$F_2L_3$&$I$&33\\
$W_{13}$&$p$&92&$F_2M$&$J$&4\\
$W_{23}$&$q$&106&$F_2w_{13}$&$K$&913\\
$H_{10}$&$r$&1074&$F_2\text{ ladder }L_6$&$L$&76\\
$M_1$&$s$&253&$F_2\text{ ladder }h$&$N$&48\\
$V$&$t$&4&&&\\
\bottomrule
\end{tabular}
\end{table}

\subsection{The forty-symmetry presentation}\label{app:char2source}
The coordinate dictionary and relations below specify the source. All coordinates are symmetries except the unitaries $N,P$. The aliases in this subsection are local; in particular, $d_0$ differs from the objective coordinate $d=z_1$.

\begin{longtable}{ll@{\qquad}ll}
\toprule
Literal coordinate & Alias & Literal coordinate & Alias\\\midrule
\endfirsthead
\toprule Literal coordinate & Alias & Literal coordinate & Alias\\\midrule
\endhead
\texttt{x0} & $x$ & \texttt{y0} & $y$\\
\texttt{z0\_0} & $z_0$ & \texttt{z1\_0} & $z_1$\\
\texttt{z2\_0} & $z_2$ & \texttt{z3\_0} & $z_3$\\
\texttt{p12\_0} & $p_{12}$ & \texttt{S1} & $S_1$\\
\texttt{S3} & $S_3$ & \texttt{R13} & $R_{13}$\\
\texttt{R21} & $R_{21}$ & \texttt{R31} & $R_{31}$\\
\texttt{R32} & $R_{32}$ & \texttt{T31} & $T_{31}$\\
\texttt{W12} & $W_{12}$ & \texttt{W13} & $W_{13}$\\
\texttt{W23} & $W_{23}$ & \texttt{H1\_0} & $H_1$\\
\texttt{M1} & $M_1$ & \texttt{V} & $V$\\
\texttt{c0} & $c$ & \texttt{d0\_0} & $d_0$\\
\texttt{B1} & $B_1$ & \texttt{C1} & $C_1$\\
\texttt{C2} & $C_2$ & \texttt{Fgraph} & $F$\\
\texttt{Ggraph} & $G$ & \texttt{Tau} & $\tau_0$\\
\texttt{Dsign13} & $D_{13}$ & \texttt{Jhybrid} & $J$\\
\texttt{T21} & $T_{21}$ & \texttt{F2\_A} & $A$\\
\texttt{F2\_B} & $B$ & \texttt{F2\_E2} & $E_2$\\
\texttt{F2\_L3} & $L_3$ & \texttt{F2\_M} & $M$\\
\texttt{F2\_w13} & $w$ & \texttt{F2\_torus\_N} & $N$\\
\texttt{F2\_steinberg\_P} & $P$ & \texttt{F2\_ladder\_L6} & $L_6$\\
\texttt{F2\_ladder\_L12} & $L_{12}$ & \texttt{F2\_ladder\_h} & $h_0$\\
\bottomrule
\end{longtable}
Write the following equations as $u_j=v_j$, in the displayed order.
The first column is the immutable row number.
\begin{longtable}{r@{\quad}l}
\toprule
$j$ & $u_j=v_j$\\\midrule
\endfirsthead
\toprule $j$ & $u_j=v_j$\\\midrule
\endhead
1 & $z_0\,z_1=z_1\,z_0$\\
2 & $z_1\,z_2=z_2\,z_1$\\
3 & $S_3\,R_{21}=R_{21}\,S_3$\\
4 & $S_1\,R_{31}=T_{31}\,S_1$\\
5 & $S_1\,R_{32}=R_{32}\,S_1$\\
6 & $R_{21}\,p_{12}=H_1\,R_{21}$\\
7 & $R_{21}\,z_2=z_2\,R_{21}$\\
8 & $R_{21}\,z_3=z_3\,R_{21}$\\
9 & $S_1\,x=x\,S_1$\\
10 & $S_1\,z_2=z_2\,S_1$\\
11 & $V\,x=x\,V$\\
12 & $V\,z_2=z_2\,V$\\
13 & $V\,p_{12}=p_{12}\,V$\\
14 & $V\,z_3=z_3\,V$\\
15 & $B_1\,x=x\,B_1$\\
16 & $C_1\,d_0=d_0\,C_1$\\
17 & $B_1\,C_2=C_2\,B_1$\\
18 & $S_1\,F=F\,S_1$\\
19 & $S_1\,G=G\,S_1$\\
20 & $C_1\,W_{23}=W_{23}\,C_1$\\
21 & $z_1\,W_{23}=W_{23}\,z_1$\\
22 & $S_1\,W_{23}=W_{23}\,S_1$\\
23 & $S_3\,W_{12}=W_{12}\,S_3$\\
24 & $G\,W_{13}=W_{13}\,G$\\
25 & $G\,W_{12}=W_{12}\,G$\\
26 & $C_1\,C_2=C_2\,C_1$\\
27 & $\tau_0\,y=d_0\,\tau_0$\\
28 & $\tau_0\,c=z_0\,\tau_0$\\
29 & $F\,z_1=z_1\,F$\\
30 & $R_{21}\,F=F\,R_{21}$\\
31 & $R_{21}\,G=G\,R_{21}$\\
32 & $R_{21}\,x=x\,R_{21}$\\
33 & $D_{13}\,R_{31}=R_{31}\,D_{13}$\\
34 & $\tau_0\,G=G\,\tau_0$\\
35 & $\tau_0\,R_{21}=R_{21}\,\tau_0$\\
36 & $\tau_0\,W_{23}=W_{23}\,\tau_0$\\
37 & $\tau_0\,V=V\,\tau_0$\\
38 & $\tau_0\,W_{13}=W_{13}\,\tau_0$\\
39 & $d_0=P\,P$\\
40 & $S_1\,C_2=C_2\,S_1$\\
41 & $N^*\,S_1=S_1\,N$\\
42 & $N\,z_0=z_0\,N$\\
43 & $N\,c=c\,N$\\
44 & $C_1\,y=y\,C_1$\\
45 & $L_3=A\,B\,A\,B\,A$\\
46 & $d_0=M\,c\,M$\\
47 & $z_0=M\,y\,M$\\
48 & $E_2=M\,F\,M$\\
49 & $V=G\,M$\\
50 & $G=A\,E_2\,A$\\
51 & $w=c\,z_0\,c$\\
52 & $\tau_0=w\,F$\\
53 & $W_{12}=A\,F$\\
54 & $W_{23}=x\,B$\\
55 & $W_{13}=A\,W_{23}\,A$\\
56 & $R_{31}=J\,W_{13}\,J$\\
57 & $R_{21}=B\,R_{31}\,B$\\
58 & $R_{13}=R_{31}\,W_{13}\,R_{31}$\\
59 & $R_{32}=A\,R_{31}\,A$\\
60 & $C_2=R_{13}\,J$\\
61 & $C_1=F\,B_1\,F$\\
62 & $B_1=A\,C_2\,A$\\
63 & $z_2=J\,z_0\,J$\\
64 & $T_{21}=C_2\,R_{21}\,C_2$\\
65 & $M_1=R_{32}\,R_{21}\,R_{32}$\\
66 & $z_1=A\,z_2\,A$\\
67 & $z_3=W_{23}\,z_2\,W_{23}$\\
68 & $T_{31}=B\,T_{21}\,B$\\
69 & $p_{12}=J\,z_1\,J$\\
70 & $S_1=M_1\,T_{31}\,R_{21}$\\
71 & $S_3=W_{13}\,S_1\,W_{13}$\\
72 & $H_1=S_1\,z_1\,S_1$\\
73 & $D_{13}=S_1\,S_3$\\
74 & $N=C_1\,B_1$\\
75 & $P=x\,z_0$\\
76 & $L_6=L_3\,B\,L_3$\\
77 & $L_{12}=L_6\,B\,L_6$\\
78 & $c=L_{12}\,B$\\
79 & $h_0=d_0\,c\,d_0$\\
80 & $x=c\,h_0$\\
81 & $y=x\,L_6\,B$\\
82 & $F=B\,L_3\,h_0$\\
\bottomrule
\end{longtable}
To pass to the source used in the theorem, make the literal substitutions
\[
 N=C_1B_1,\qquad P=xz_0,\qquad N^*=B_1C_1,\quad P^*=z_0x.
\]
Delete rows $74$ and $75$, which now vanish identically. Let
$I=\{1,\ldots,82\}\setminus\{74,75\}$, and let $\eta_j$ be the
substituted difference $u_j-v_j$. The remaining forty coordinates are
independent symmetries. Precisely,
\[
 Q_{80}=\sum_{j\in I}\eta_j^*\eta_j,\qquad
 F_0=-1+Vz_1Vz_1+z_1Vz_1V.
\]
Thus the norm of each residual is at most
$\sqrt{\operatorname{tr}_n Q_{80}}$ on every symmetry tuple.
\subsection{The Steinberg bridge}\label{app:char2bridge}
\label{char2merge:sec:steinberg-bridge}

We derive the group-theoretic spectral input from a subset of the cubic source relations.

Let $P$ denote the group presented by the $32$ source involutions and the
$79$ cubic paid equations.  We use the abbreviations
\[
 x=x0,\quad y=y0,\quad z=z0\_0,\quad z_i=z_i\_0\ (1\leq i\leq3),
 \quad c=c0,\quad d=d0\_0,
\]
\[
\begin{gathered}
 S=S1,\quad p=W12,\quad q=W23,\quad r=W13,\\
 F=Fgraph,\quad G=Ggraph,\quad \tau=Tau,
 \quad C=C1,\quad B=B1.
\end{gathered}
\]
Here and below $[a,b]=aba^{-1}b^{-1}$. All named source coordinates, including $Kxy\_00$ and $Dsign13$, are involutions. The symbol $S_3$ denotes the source coordinate $S3$, not a symmetric group.

\subsubsection{Paid relations}

The following are literal paid equations, or the same equations multiplied
by the inverse of an involution on one side:
\begin{align}
 xy&=(Kxy\_00)x,& (Kxy\_00)y&=c,& zx&=dxz,
 \label{char2merge:eq:st-constant-input}\\
 FxF&=y,& FzF&=z,& GyG&=z,
 \label{char2merge:eq:st-graph-input}\\
 \tau y\tau&=d,& \tau c\tau&=z,& [F,\tau]&=1,
 \label{char2merge:eq:st-tau-input}\\
 [S,x]&=[S,F]=[S,G]=1,
 \label{char2merge:eq:st-sign-input}\\
 [p,x]&=[p,F]=[p,G]=1,&
 [r,x]&=[r,F]=[r,G]=1,
 \label{char2merge:eq:st-perm-root-input}\\
 pq&=rp,& qr&=pq,& [q,S]&=[q,C]=1,
 \label{char2merge:eq:st-perm-input}\\
 rSr&=S_3,& SS_3&=Dsign13,& pCp&=C_2,
 \label{char2merge:eq:st-orbit-input}\\
 [C,C_2]&=1,& [B,C_2]&=1,& [S,C_2]&=1,
 \label{char2merge:eq:st-factor-input}\\
 FBF&=C,& CSB&=BSC,
 \label{char2merge:eq:st-same-index-input}\\
 [C,y]&=[C,d]=1,& [CB,z]&=[CB,c]=1,
 \label{char2merge:eq:st-stabilizer-input}\\
 C_2zC_2&=z_2,& pz_1p&=z_2,& qz_2q&=z_3.
 \label{char2merge:eq:st-variable-input}
\end{align}
Thus the proof uses only finitely many displayed group relations.  The
additional source equations need not be discarded: every conclusion from
this subset holds in the full quotient $P$.

\subsubsection{The constant $A_2$ system}

Equations \eqref{char2merge:eq:st-constant-input} give
\[
 [x,y]=c,\qquad [z,x]=d.
\]
Put $v=\tau x\tau$.  This is an involution, and
\eqref{char2merge:eq:st-graph-input}--\eqref{char2merge:eq:st-tau-input} give $FdF=v$.
Conjugating $[z,x]=d$ by $F$ gives $[z,y]=v$.  Since
$[y,z]=[z,y]^{-1}$, it follows that $[y,z]=v$ as well.
Conjugating these three chain relations by $\tau$ gives
\[
 [v,d]=z,\qquad [c,v]=y,\qquad [d,c]=x.
\]
With the identification
\[
 (u_{12},u_{23},u_{31},u_{13},u_{32},u_{21})=(x,y,z,c,d,v),
\]
these give the six relations
$[u_{ij},u_{jk}]=u_{ik}$ for distinct $i,j,k$, after reversing
commutators whose values are involutions where necessary.

For any involutions $a,b$, if $[a,b]$ is also an involution then
\[
 [a,b]=(ab)^2=(ba)^2,
\]
and conjugation by either $a$ or $b$ fixes this element.  Applying this
observation to $[x,y]=c$, $[z,x]=d$, and $[y,z]=v$ gives exactly the six
constant same-row and same-column commuting pairs.  These are the remaining
constant distinct nonopposite-root relations of type $A_2$.

Equations \eqref{char2merge:eq:st-sign-input} show that $S$ fixes $x,y,z$ under
conjugation, using $y=FxF$ and $z=GyG$.  It therefore fixes $c,d,v$ by
their commutator expressions.  The same argument using
\eqref{char2merge:eq:st-perm-root-input} shows that $p,r$ fix all six constant roots.
There is no assumption here that $S$ or $p$ commutes with $\tau$.

\subsubsection{Six commuting torus words}

The involutions $p,q,r$ and \eqref{char2merge:eq:st-perm-input} imply
\[
 r=pqp=qpq.
\]
Thus conjugation by them obeys the permutation identities for the
transpositions $(12),(23),(13)$.  Injectivity of this permutation action
is not needed.  Define
\[
 S_1=S,\quad S_2=pSp,\quad S_3=rSr,
 \qquad C_1=C,\quad C_2=pCp,\quad C_3=rCr.
\]
The commuting relations $[q,S]=[q,C]=1$ make these well-defined
three-point permutation orbits; the displayed $S_3,C_2$ agree with the
source coordinates by \eqref{char2merge:eq:st-orbit-input}.
Since $q=prp$, all three permutations fix the constant roots.  Consequently
every $S_i$ fixes all six constant roots.

Involutivity of $Dsign13=S_1S_3$ implies $[S_1,S_3]=1$.
Transport by index permutations gives $[S_i,S_j]=1$ for all $i,j$.
Put
\[
 E_i=FC_iF,\qquad A_i=C_iS_i,\qquad D_i=E_iS_i.
\]
The symbols $D_i$ here are torus words.  Since $F$ commutes with $S,p,r$, it commutes
with each $S_i$ and preserves index covariance.  In particular $E_1=B$.
The finite relations \eqref{char2merge:eq:st-factor-input}, their index conjugates,
and their $F$-conjugates give, for $i\ne j$,
\[
 [C_i,C_j]=[E_i,E_j]=[C_i,E_j]
 =[S_i,C_j]=[S_i,E_j]=[S_i,S_j]=1.
\]
It follows that all off-index torus words commute.  The paid same-index
equation $C_1S_1E_1=E_1S_1C_1$, multiplied on the right by $S_1$,
gives $A_1D_1=D_1A_1$.  Its index conjugates give the remaining
same-index commutations.  Thus all six words $A_i,D_i$ commute.

Their root stabilizers are
\begin{align}
 [A_i,y]&=[A_i,d]=1,\label{char2merge:eq:st-torus-stab-a}\\
 [D_i,x]&=[D_i,v]=1,\label{char2merge:eq:st-torus-stab-d}\\
 [A_iD_i^{-1},z]&=[A_iD_i^{-1},c]=1.
 \label{char2merge:eq:st-torus-stab-ad}
\end{align}
For the first line, use \eqref{char2merge:eq:st-stabilizer-input}, the
$S_i$-invariance, and index transport.  The second follows by
$F$-conjugation because $FyF=x$ and $FdF=v$.  For the third, note the
exact cancellation
\[
 A_iD_i^{-1}=C_iS_i(S_iE_i)=C_iE_i;
\]
then use the two kernel equations in
\eqref{char2merge:eq:st-stabilizer-input} and index transport.  This cancellation
does not assume $[C_i,S_i]=1$.
Finally \eqref{char2merge:eq:st-variable-input} and $S_2$-invariance of $z$ give
$A_2zA_2^{-1}=z_2$; conjugating by $p$ and $q$ gives
\begin{equation}
 A_i z A_i^{-1}=z_i\qquad(1\leq i\leq3).
 \label{char2merge:eq:st-variable-identification}
\end{equation}

\subsubsection{Construction of every Laurent root}

Let $R=\mathbb F_2[t_1^{\pm1},t_2^{\pm1},t_3^{\pm1}]$.
For $m\in\mathbb Z^3$ write $A(m)=\prod_i A_i^{m_i}$ and
$D(m)=\prod_iD_i^{m_i}$.  Define
\begin{align*}
 X_{12}(m)&=A(-m)xA(m),& X_{21}(m)&=A(m)vA(-m),\\
 X_{23}(m)&=D(-m)yD(m),& X_{32}(m)&=D(m)dD(-m),\\
 X_{13}(m)&=A(-m)cA(m),& X_{31}(m)&=A(m)zA(-m).
\end{align*}
The commuting torus and
\eqref{char2merge:eq:st-torus-stab-a}--\eqref{char2merge:eq:st-torus-stab-ad} show that
conjugation by $A(a)D(b)$ shifts $m$ according to this table:
\[
\begin{array}{c|rrrrrr}
 ij&12&21&23&32&13&31\\\hline
 \text{shift}&-a&a&-b&b&-a-b&a+b.
\end{array}
\]
For example, $A_iD_i^{-1}$ fixes $c$, so the $D_i$-action on $c$
is the same as the $A_i$-action; this proves the $13$ entry.  The other
entries follow identically from the displayed stabilizers.

The six weight vectors are $\pm(1,0),\pm(0,1),\pm(1,1)$.
Any distinct nonopposite pair comes from two different directions, so its
$2\times2$ coefficient matrix has determinant $\pm1$.  Therefore one
can prescribe arbitrary two shifts $m,n\in\mathbb Z^3$ by solving the
corresponding two integral linear equations separately in each variable
coordinate.  Simultaneous torus conjugation of a constant same-row or
same-column commutation proves every such monomial commutation.
Simultaneous conjugation of a constant chain relation proves
\begin{equation}
 [X_{ij}(m),X_{jk}(n)]=X_{ik}(m+n),\qquad i,j,k\text{ distinct},
 \label{char2merge:eq:st-monomial-chain}
\end{equation}
since the output weight is the sum of the input weights.

Each $X_{ij}(m)$ is an involution.  To prove same-root commutation,
choose the third index $k$ and write
\[
 X_{ij}(m)=[X_{ik}(m),X_{kj}(0)].
\]
Both factors commute with $X_{ij}(n)$ by the mixed same-row and
same-column relations already proved.  Hence $X_{ij}(m)$ commutes with
$X_{ij}(n)$ for every $m,n$.

For $f=\sum_{m\in\operatorname{supp}(f)}t^m\in R$, define
\[
 X_{ij}(f)=\prod_{m\in\operatorname{supp}(f)}X_{ij}(m).
\]
This is unambiguous and satisfies
$X_{ij}(f)X_{ij}(g)=X_{ij}(f+g)$: same-root terms commute and repeated
terms cancel by involutivity.  Mixed commuting relations extend termwise.
For a chain, all elementary output terms in
\eqref{char2merge:eq:st-monomial-chain} centralize both input-root subgroups, by
the mixed same-row and same-column relations.  The class-two commutator
identity therefore expands the polynomial commutator into the product of
all pairwise outputs.  Equal monomials cancel in pairs, giving
\[
 [X_{ij}(f),X_{jk}(g)]=X_{ik}(fg).
\]
These are the Steinberg defining relations.  The universal
property supplies a homomorphism
\begin{equation}
 \Phi:\operatorname{St}_3(R)\longrightarrow P.
 \label{char2merge:eq:st-map}
\end{equation}
There is no assertion that the torus words belong to its image or that
$P$ itself has property~(T).

\subsubsection{Polynomial image and transport}

Put $R_0=\mathbb F_2[t_1,t_2,t_3]$ and
$K=\langle x,y,z,z_1,z_2,z_3\rangle\leq P$.
The functorial map $\operatorname{St}_3(R_0)\to\operatorname{St}_3(R)$,
followed by $\Phi$, has image exactly $K$.

Within the abstract Steinberg group, write
\[
 a=x_{12}(1),\quad b=x_{23}(1),\quad s=x_{31}(1),
 \quad Z=x_{31}(t),
\]
where $t$ is any one of the three variables.  The following are abstract
Steinberg identities, with the indicated lengths in $\{a,b,s,Z\}$:
\[
\begin{array}{c|l|r}
\text{root}&\text{word}&\text{length}\\\hline
x_{31}(t)&Z&1\\
x_{32}(t)&[Z,a]&4\\
x_{21}(t)&[b,Z]&4\\
x_{12}(t)&(babZaZ)^2&12\\
x_{23}(t)&(ZbZaba)^2&12\\
x_{13}(t)&s[a,b]Z[a,b]s&11.
\end{array}
\]
Here all root elements are involutions by additivity over $\mathbb F_2$.
To verify the two length-$12$ identities, put
$E=[a,b]$, $H=[Z,a]$ and $J=[b,Z]$.  Because these commutators are
involutions, $E=baba$, $H=aZaZ$, $J=ZbZb$.  Free cancellation gives
\[
 EH=babZaZ,\qquad JE=ZbZaba.
\]
The chain relations $[E,H]=x_{12}(t)$ and
$[J,E]=x_{23}(t)$ now give the two square formulas.

The length-$11$ formula also follows from the defining
relations without invoking a faithful matrix representation.  Set $w=sEs$.
This is an involution.  Conjugation by $s$ fixes $x_{32}(t)$;
conjugation by $E$ multiplies this element by $x_{12}(t)$; and conjugation
by $s$ multiplies the latter factor by $x_{32}(t)$.  The two
$x_{32}(t)$ factors commute with $x_{12}(t)$ and cancel.  Thus
\[
 wx_{32}(t)w=x_{12}(t).
\]
The same calculation on the other side gives
$wx_{21}(1)w=x_{23}(1)$.  Applying these identities to
$Z=[x_{32}(t),x_{21}(1)]$ gives $wZw=x_{13}(t)$.
Since $s$ and $Z$ commute, $wZw=sEZE s$, exactly the displayed word.

All six constant roots have words of length at most $4$ in $a,b,s$.
Therefore the standard $24$ roots
\[
 \Sigma=\{x_{ij}(u):i\ne j,\ u\in\{1,t_1,t_2,t_3\}\}
\]
have words of length at most $12$ in
$\{x_{12}(1),x_{23}(1),x_{31}(1),x_{31}(t_1),x_{31}(t_2),x_{31}(t_3)\}$.
The set $\Sigma$ generates $\operatorname{St}_3(R_0)$: induction on
monomial degree uses $x_{ij}(mn)=[x_{ik}(m),x_{kj}(n)]$, and root
additivity supplies arbitrary polynomial coefficients.  By
\eqref{char2merge:eq:st-variable-identification}, the six chosen elements map to
$x,y,z,z_1,z_2,z_3$.  Every image of $\Sigma$ lies in $K$, and these six
source roots are themselves polynomial-root images.  This proves the
claimed equality of the image with $K$.

\subsubsection{The numerical spectral gap}

Ershov--Jaikin-Zapirain~\cite[Theorem~6.2]{ErshovJaikinZapirain2010} give the Steinberg-group bound

\[
 \kappa(\operatorname{St}_n(R_0),\Sigma)
 \geq\frac{1}{8(12\sqrt{2d}+2\sqrt{3n}+36\sqrt2)},
\]
where $1$ and $d$ additional elements generate the coefficient ring.
Here $n=d=3$.  Since $\sqrt6<5/2$ and $\sqrt2<3/2$, the denominator
is strictly less than
\[
 8(30+6+54)=720.
\]
The theorem applies to arbitrary finitely generated associative unital rings, including $R_0$ in characteristic two.

A representation of $K$ pulls back through the surjection just proved,
with the same invariant subspace.  The length-$12$ transport words and
the triangle inequality for displacements therefore give
\[
 \kappa(K,\mathcal S)>\frac1{8640},\qquad
 \mathcal S=(x,y,z,z_1,z_2,z_3).
\]
For the lazy Laplacian
\[
 \Delta=\frac1{12}\sum_{s\in\mathcal S}(1-s)
\]
and a unitary representation $\pi$,
\[
 \langle\pi(\Delta)\xi,\xi\rangle
 =\frac1{24}\sum_{s\in\mathcal S}\|\pi(s)\xi-\xi\|^2.
\]
On the orthogonal complement of invariant vectors its quadratic form is
therefore bounded below by a constant strictly greater than
\[
 \frac1{24\cdot8640^2}
 =\frac1{1,791,590,400}>2^{-31},
 \qquad 1,791,590,400<2^{31}=2,147,483,648.
\]
Thus the universal spectral gap is strictly greater than $2^{-31}$, so $\lambda_H=2^{-32}$ lies below it.

\subsubsection{Rational SOS and a terminating paid-certificate search}

The rational clause of Ozawa's Main Theorem~\cite{Ozawa2016} applies to the measure
\[
 \mu=\tfrac12\delta_1+\tfrac1{12}\sum_{s\in\mathcal S}\delta_s,
\]
whose Laplacian is $\Delta$. A positive
rational parameter strictly below the universal spectral gap gives
\begin{equation}
 \Delta^2-\lambda_H\Delta=\sum_{i=1}^N q_i\xi_i^*\xi_i
 \quad\text{in }\mathbb Q[K],\qquad
 q_i\in\mathbb Q_{>0},\quad\xi_i\in\mathbb Q[K].
 \label{char2merge:eq:st-rational-sos}
\end{equation}

Let $\mathcal F$ be the free product of the $32$ order-two coordinate
groups.  Write the $79$ paid differences as
$\rho_j=w_j^+-w_j^-\in\mathbb Q[\mathcal F]$.
Since $K$ is a subgroup of $P$, its group algebra embeds into
$\mathbb Q[P]$.  Choose lifts of the $\xi_i$ using words in the six
root letters.  The residual in \eqref{char2merge:eq:st-rational-sos} then lies in
the kernel of $\mathbb Q[\mathcal F]\to\mathbb Q[P]$, which is the
two-sided ideal generated by the $\rho_j$.  Consequently there is a
finite exact identity
\begin{equation}
 \Delta^2-\lambda_H\Delta-\sum_iq_i\widetilde\xi_i^*\widetilde\xi_i
 =\sum_{\ell=1}^M c_\ell a_\ell\rho_{j_\ell}b_\ell
 \quad\text{in }\mathbb Q[\mathcal F],
 \label{char2merge:eq:st-paid-sos}
\end{equation}
with rational $c_\ell$ and free-product words $a_\ell,b_\ell$.
Indeed, the quotient by this ideal imposes exactly
$w_j^+=w_j^-$, so it is the group algebra of the presented quotient.
Every expression in a generated algebraic ideal is a finite sum, and
polynomial contexts expand into word contexts.

Use the following exact encoding convention. Enumerate all
finite strings over ASCII characters $32,\ldots,126$ by increasing
length and, within each length, lexicographically in ASCII order. Retain
only canonical compact JSON strings, with object keys sorted
lexicographically and no whitespace outside string values, with the
following schema:
\begin{itemize}
\item The top object has exactly the keys \texttt{lambda},
\texttt{squares}, \texttt{ideal\_terms}. Its \texttt{lambda} value is
exactly the string \texttt{1/4294967296}.
\item \texttt{squares} is a list of objects with exactly the keys
\texttt{weight}, \texttt{terms}. Each weight is a strictly positive
rational string. Each term is an object with exactly
\texttt{coefficient}, \texttt{word}; the latter is a list of the six
literal root names in $\mathcal S$.
\item \texttt{ideal\_terms} is a list of objects with exactly the keys
\texttt{coefficient}, \texttt{left\_context}, \texttt{relation},
\texttt{right\_context}. Contexts are lists of literal names of the
thirty-two base symmetries; the relation is an exact base row name.
\end{itemize}
Rational strings use reduced numerator and positive denominator,
separated by a slash, except that denominator $1$ is omitted; zero is
\texttt{0}. The lists may be empty. Their order is part of the string;
terms in the represented group-algebra elements are collected by exact
rational addition. Reject all other strings. Check
\eqref{char2merge:eq:st-paid-sos} by adjacent involution cancellation and exact
rational arithmetic, as in the supplied verifier. Each candidate test
is finite, and every required rational identity has a string in this
encoding. Existence of \eqref{char2merge:eq:st-paid-sos} therefore proves that the
search terminates. Thus the first
successful identity defines a finite rational constant
\[
 C_H=\max\left\{1,\ 2\sum_{\ell=1}^M|c_\ell|\right\}.
\]

Finally, let an arbitrary symmetry tuple in $M_d(\mathbb C)$ satisfy
$\|\rho_j\|_2\leq\varepsilon$ for all paid rows.  Conjugation by the
tuple defines an exact representation of $\mathcal F$ on normalized
Hilbert--Schmidt space, although it need not factor through $P$.
For every contraction $h$,
\[
 \|\operatorname{Ad}(w_j^+)(h)-
       \operatorname{Ad}(w_j^-)(h)\|_2
 \leq 2\|w_j^+-w_j^-\|_2\leq2\varepsilon.
\]
Right word contexts preserve the operator-norm unit ball; left word
contexts preserve the $2$-norm.  Evaluating \eqref{char2merge:eq:st-paid-sos},
discarding its positive square terms, and using $\|h\|_2\leq1$ therefore
gives the paid inequality
\begin{equation}
 \left\langle(L_H^2-2^{-32}L_H)h,h\right\rangle
 \geq-C_H\varepsilon,
 \qquad
 L_H=\frac1{12}\sum_{s\in\mathcal S}(I-\operatorname{Ad}(s)).
 \label{char2merge:eq:st-paid-gap}
\end{equation}
\subsection{Paid integer, compression, and presentation transfers}\label{app:char2transfers}
\label{char2merge:sec:paid-transfers}

This section proves the finite hypotheses used in the analytic argument
below, and the comparison with the exact forty-symmetry source.
For the thirty-two-symmetry base presentation, put
\[
 q_{79}=\operatorname{tr}_n Q_{79},\qquad \varepsilon=\sqrt{q_{79}},
 \qquad S=(x_0,y_0,z_0,z_1,z_2,z_3),\quad z_i=z_{i,0}.
\]
Thus every paid residual has normalized $2$-norm at most $\varepsilon$.
The base rows are specified in
\path{quadratic_214/base_32_symmetries_79_rows.json}.

\subsubsection{Sound rules for the finite certificates}
\label{char2merge:subsec:finite-rules}
All words in this section evaluate to unitaries. Consequently
\[
 \|aub-avb\|_2=\|u-v\|_2,
 \qquad \|u^*-v^*\|_2=\|u-v\|_2
\]
for unitary contexts $a,b$. A difference can be reversed at no cost;
products and transitive chains add their norm bounds. Adjacent inverse
letters cancel exactly. For symmetry letters the inverse letter is the
letter itself. Cyclically rotating a relator conjugates it by a unitary;
inverting a relator preserves its defect norm. These justify respectively
the certificate operations \texttt{paid}, \texttt{reflexive},
\texttt{reverse}, \texttt{adjoint}, \texttt{multiply},
\texttt{transitive}, and \texttt{cyclic-relator}. Every referenced node
must precede the node using it.

For signed contextual identities one instead checks, by collecting free
words over $\mathbb Z$, an identity
\[
 u-v=\sum_\ell c_\ell a_\ell\rho_{j_\ell}^{\sigma_\ell}b_\ell,
 \qquad c_\ell\in\mathbb Z,
 \quad \sigma_\ell\in\{1,*\}.
\]
If $m=\sum_\ell|c_\ell|$, the triangle inequality and weighted
Cauchy--Schwarz imply
\[
 \|u-v\|_2^2\leq
 m\sum_\ell |c_\ell|\,\|\rho_{j_\ell}\|_2^2.
\]
Adding these estimates gives a penalty comparison from the largest
column sum of the resulting nonnegative coefficient matrix.

\subsubsection{The six actions and their forty-two compression estimates}
Let $u=R_{21}S_1$, $t=Vu$, and define
\[
\begin{array}{c|rrrrrr}
 ab&21&12&31&13&23&32\\\hline
 P_{ab}&1&W_{12}&W_{23}&W_{23}W_{12}&W_{13}&W_{12}W_{23}
\end{array}
\qquad U_{ab}=P_{ab}uP_{ab}^*.
\]
Define the three product targets
\[
 p_{12}=p_{12,0},\qquad
 p_{13}=W_{23}p_{12,0}W_{23},\qquad
 p_{23}=W_{13}p_{12,0}W_{13}.
\]
For $a\ne b$, the target $p_{ab,g}$ is $g$ when
$g=x_0,y_0,z_0$, or $g=z_i$ with $i\ne b$; when $g=z_b$ it is
$p_{\min(a,b),\max(a,b)}$. These words are all exact symmetries.
The common targets for $u$ and $t$ are $p_g=p_{21,g}$.

We first use an intermediate presentation in which the torus symbol is
$B_2$ instead of $B_1$. Its seventy-three non-torus rows suffice for
all estimates in this subsection; let $\varepsilon_o$ be the square root
of their total squared defect. The compression and root
certificates are specialized by making every root letter a symmetry.
Only their ancestor nodes are used: $1247$ compression nodes with
$35$ distinct paid premises, and $970$ root/product nodes. No order-five
premise is used.

There are two apparent extra centrality premises, which must be derived.
Write $x=x_0,y=y_0,c=c_0,k=xyxy$. The two paid rows
$xy=K_{xy}x$ and $K_{xy}y=c$ give $\|c-k\|_2\leq2\varepsilon_o$.
Since $c,x,y$ are self-adjoint,
\[
 cx-xc=(c-k)x-x(c-k^*),\qquad
 cy-yc=(c-k^*)y-y(c-k).
\]
Thus their norm charges are $4$, not $1$. Verifying all affected
ancestors with these charges gives the following compression costs:
\begin{center}
\begin{tabular}{c|rrrrrr}
 action&$x_0$&$y_0$&$z_0$&$z_1$&$z_2$&$z_3$\\\hline
 $U_{21}$&2&8&14&2&2&6\\
 $U_{12}$&4&18&32&4&24&24\\
 $U_{31}$&12&42&72&3&8&4\\
 $U_{13}$&14&52&90&6&26&25\\
 $U_{23}$&4&18&32&16&12&7\\
 $U_{32}$&14&52&90&10&27&22\\
 $t$&3&13&23&3&3&7
\end{tabular}
\end{center}
An entry $c$ proves $\|U_{ab}gU_{ab}^*-p_{ab,g}\|_2\leq
c\varepsilon_o$, or its $t$ analogue. In particular, the $t$ row
uses the same six targets as the $U_{21}$ row.

The root/product paths additionally prove
\[
 \|d_a-[x_0,z_a]_{\rm grp}\|_2\leq165\varepsilon_o,
 \qquad
 \|e_b-[y_0,z_b]_{\rm grp}\|_2\leq165\varepsilon_o,
 \qquad
 \|p_{ab}-[d_a,e_b]_{\rm grp}\|_2\leq321\varepsilon_o,
\]
where $[v,w]_{\rm grp}=vwv^*w^*$. For each ordered $a\ne b$ the
last target is the corresponding unordered product above. Substitution
in the commutator therefore gives a sixteen-letter word in $S$ with
error at most
\[
 (321+4\cdot165)\varepsilon_o=981\varepsilon_o
 <4096\varepsilon_o.
\]
The source symbols $d_a,e_b$ here denote the explicit virtual conjugate
words in the certificate, so they introduce no new paid variables.

\paragraph{Cubic rebasing.}
On the final base tuple set $B_2=W_{12}B_1W_{12}$, which is an exact
symmetry. Seventy of the seventy-three rows are unchanged. The remaining
three rows have routes with the following charges and paid uses:
\begin{center}\small
\begin{tabular}{lcl}
 Source row&charge&paid uses in the cubic source\\\hline
 \texttt{cov\_B1\_x0}&3&
 $2\,\texttt{bcov\_W12\_x0}+\texttt{cov\_B1\_x0}$\\
 \texttt{comm\_B1\_C2}&3&
 $2\,\texttt{weyl\_C2}+\texttt{comm\_B1\_C2}$\\
 \texttt{graph\_F\_BC\_1}&4&
 $2\,\texttt{bcov\_W12\_Fgraph}+\texttt{graph\_F\_BC\_1}
 +\texttt{weyl\_C2}$
\end{tabular}
\end{center}
Each displayed use is a unitary-context replacement after conjugation
by $W_{12}$. The finite checker verifies all three routes and their
literal endpoints. The maximum squared column sum, including unchanged
rows, is $11$; hence
\begin{equation}\label{char2merge:eq:cubic-rebase}
 \varepsilon_o^2\leq11q_{79}.
\end{equation}
Alternatively every paid defect in the reflection presentation is at most $4\varepsilon$.
Using the latter bound in the maximum-defect compression calculations
proves all of the following, on the cubic base:
\begin{align}
 \|U_{ab}gU_{ab}^*-p_{ab,g}\|_2&\leq360\varepsilon,
 &\|p_{ab,g}-w_{ab,g}(S)\|_2&\leq16384\varepsilon,
 \label{char2merge:eq:paid-compression}\\
 \|tgt^*-p_g\|_2&\leq360\varepsilon,
 &|w_{ab,g}|&\leq16.\nonumber
\end{align}
For unchanged-root targets take the one-letter word $g$.
The finite transfer identities are verified by \path{quadratic_214/check_transfer.py}.

\subsubsection{The paid integer spectral inequality}
\label{char2merge:subsec:paid-integer-gap}
The fifteen relations for six unitaries $v_{ab}$ are
\begin{align*}
 v_{ij}v_{ik}&=v_{ik}v_{ij},&
 v_{ji}v_{ki}&=v_{ki}v_{ji}
 &&(i\in\{1,2,3\},\ j<k,\ \{i,j,k\}=\{1,2,3\}),\\
 v_{ij}v_{jk}&=v_{ik}v_{jk}v_{ij}
 &&& (i,j,k\text{ distinct}),\\
 v_{ij}v_{ji}^*v_{ij}&=v_{ji}^*v_{ij}v_{ji}^*
 &&&(i<j).
\end{align*}
Let $\delta$ be the largest $2$-norm of these fifteen balanced
residuals. We import the following precise established finite-input
estimate, not merely property~(T) of $\operatorname{SL}_3(\mathbb Z)$:
\begin{equation}\label{char2merge:eq:integer-import}
 \langle(\mathcal L_v^2-\mathcal L_v/144)a,a\rangle
 \geq-26\delta\qquad(\|a\|\leq1).
\end{equation}
Here $\mathcal L_v=\tfrac12I-\tfrac1{24}
\sum_{a\ne b}(\operatorname{Ad}v_{ab}+\operatorname{Ad}v_{ab}^*)$.
The reference is
\href{https://arxiv.org/html/2610.01536v1#S3.SS3}{Wang--Zhi,
version~1, Lemma~3.3 and equation~(3.12)}.
It combines the Netzer--Thom integer SOS with a finite transfer to these
fifteen relations.

We justify its use for the literal words $U_{ab}$, rather than
assuming their relations. The integer dependency certificate provides
these successive squared-penalty charges:
\begin{center}
\begin{tabular}{lc}
 transfer&maximum squared column charge\\\hline
 structural compact tuple to the reflection system&$11$\\
 typed product extension&$1$\\
 $528$ signed contextual identities&$22869994$\\
 restoration of fifteen deleted reflection relations&$985$
\end{tabular}
\end{center}
Here are the mathematical meanings of the two extension steps. Set all
irrelevant field/root coordinates of the reflection system to the
identity, keep its structural symmetries, and set
\[
 E_{13}=\operatorname{sign}\bigl(\operatorname{Re}(D_{13}R_{31})\bigr),
 \qquad \operatorname{sign}(0)=1.
\]
For a unitary $w$, scalar spectral calculus gives
\[
 \|w-\operatorname{sign}(\operatorname{Re}w)\|_2^2
 \leq\tfrac12\|w-w^*\|_2^2.
\]
The only nonliteral commutation row in that system, involving $\mathrm{Tau}$ and
$W_{12}$, follows from $W_{12}\simeq W_{23}W_{13}W_{23}$ and the
three corresponding paid rows, with squared charges $10,10,5$.
The other surviving rows are literal; their maximum total column charge
is $11$. To replace a commutation row $ab=ba$ by two equations using a
new symmetry, choose that symmetry to be
$\operatorname{sign}(\operatorname{Re}(ab))$. Their total squared
error is at most $\|ab-ba\|_2^2$. To split an arbitrary equality
$w=w'$ through a new unitary, choose it to be $w$, paying only the source
error. These prove the extension charge $1$ without assuming the
products are symmetries.

The $528$ identities are checked by collecting their $128289$ signed
contextual terms exactly. The last deletion step is rederived from
retained rows: its actual maximum
column charge is $325$, below the conservative $985$ used here.
Accordingly every needed reflection residual is bounded by
$\sqrt A\,\varepsilon_o$, where
\[
 A=11\cdot22869994\cdot985=247796384990.
\]
The fifteen explicit unitary paths use at most eleven reflection
replacements each. Their endpoints are the fifteen displayed integer
relations for the structural action words. Six signed conjugacy ledgers
then compare those words with the literal $U_{ab}$, with respective
norm charges $0,0,2,8,7,3$. Each side of an integer relation has at
most three letters, so this comparison adds at most $6\cdot8=48$.
We obtain
\[
 \delta\leq(11\sqrt A+48)\varepsilon_o.
\]
All $15$ endpoints, every step, the six conjugacy ledgers, and the
nonnegative column sums are verified in
\path{quadratic_214/check_bridge.py}.
The outer checker also verifies that all $37$ compact paid rows involved
are literal rows of this characteristic-two system. This specialization
uses no characteristic-five matrix-separation theorem.

Finally $\lceil\sqrt A\rceil=497792$ and
\[
 26(11\lceil\sqrt A\rceil+48)=142369760<2^{28}.
\]
Equations~\eqref{char2merge:eq:integer-import} and \eqref{char2merge:eq:cubic-rebase}, together
with $\sqrt{11}<4$, prove
\begin{equation}\label{char2merge:eq:paid-integer-gap}
 \langle(\mathcal L_U^2-\mathcal L_U/144)a,a\rangle
 \geq-2^{30}\varepsilon\qquad(\|a\|\leq1).
\end{equation}
Thus the integer spectral input and both compression inputs needed
below have been proved for the specified base presentation. The root
spectral input is supplied by Section~\ref{char2merge:sec:steinberg-bridge}.

\subsubsection{A direct 79 to 80 penalty comparison}
\label{char2merge:subsec:direct-transfer}
Let $s$ be an arbitrary tuple of the forty target symmetries.
First substitute the two eliminated unitary coordinates literally:
\[
 N=\mathrm{F2\_torus\_N}=C_1B_1,
 \qquad T=\mathrm{F2\_steinberg\_P}=x_0z_0.
\]
Their two binding residuals become zero, leaving exactly the eighty
residuals defining $Q_{80}$. Fill the base tuple by retaining the
thirty-one common symmetry coordinates and setting
\[
 K_{xy,00}=x_0y_0x_0.
\]
This last word is an exact symmetry. In particular the fill preserves
$V,z_1$ and hence $F_0$ exactly.

\begin{lemma}\label{char2merge:lem:direct-79-80}
For every such tuple, in every matrix size,
\[
 \operatorname{tr}_n Q_{79}(\operatorname{fill}(s))
 \leq149778\operatorname{tr}_n Q_{80}(s).
\]
\end{lemma}
\begin{proof}
Put $e=\sqrt{\operatorname{tr}_n Q_{80}(s)}$.
For a free word $w$ on these forty symmetry letters write
$E(w)=\|1-w(s)\|_2$. For a residual $a-b$ one has
$\|a-b\|_2=E(ab^{-1})$.
Freely reduce $w$, remove matching boundary letters, and select the
lexicographically least cyclic rotation of it or its inverse; call this
its canonical cap. Each operation preserves $E$ exactly.

The certificate \path{quadratic_214/certificate_79_to_80.json} contains
$84$ ordered paths: one for each of the $79$ base caps, and five helper
caps. Every path ends at the empty word. A step rotates the current cap,
writes the result as $lv$, and replaces it by the cap of $rv$.
The cap of $lr^{-1}$ must be a paid target cap or a previously proved
helper/base cap with known charge $c$. By unitary invariance,
\[
 E(lv)\leq\|l-r\|_2+E(rv)\leq ce+E(rv).
\]
Thus induction on the paths and on their recorded dependency order
proves $E(w_j)\leq L_je$, where $L_j$ is the sum of charges along
its path. There are $322$ steps. The complete multiset of base charges
is specified by the following table:
\begin{center}
\begin{tabular}{c|rrrrrrrrrr}
 $L$&0&1&2&4&6&8&11&12&14&16\\\hline
 multiplicity&1&55&4&4&1&1&1&1&1&1
\end{tabular}
\medskip

\begin{tabular}{c|rrrrrrrrr}
 $L$&25&28&36&38&48&57&76&128&342\\\hline
 multiplicity&1&1&1&1&1&1&1&1&1
\end{tabular}
\end{center}
The multiplicities total $79$, and their weighted square sum is
$\sum_jL_j^2=149778$. Squaring the individual residual estimates
and adding proves the result.

The standalone checker \path{quadratic_214/verify_transfer.py} reconstructs
every source and goal cap directly from the two JSON presentations,
verifies the substitutions and exact involution types, rejects forward
references, checks each literal replacement and endpoint, and recomputes
this sum.
\end{proof}

\Needspace{6\baselineskip}
\subsection{The dimension-independent matrix bound}\label{app:char2norm}
\subsubsection{The finite analytic step in characteristic two}
\label{char2merge:subsec:char2-normalization}

This subsection uses the numerical procedure of Wang--Zhi, version~1,
Corollary~4.3, Theorem~4.1, and equations~(3.21)--(3.22).

\paragraph{The fixed rational tolerance function.}
For an integer $h\geq1$ and rational parameters
$0<\lambda<1$, $C\geq0$, $0<\theta<1$, define
$\mathfrak e(h,\lambda,C,\theta)$ by the following finite sequence of
operations. This is the function $\varepsilon_{\rm in}$ in the cited
version~1, with the auxiliary letter $M$ renamed $m_{\rm cut}$:
\begingroup\small
\begin{align*}
 \kappa&=\lambda/2,& c&=\kappa^2/2^{28},&
 \kappa'&=\kappa/1024,\\
 N&=\left\lceil\frac{16}{c\theta}\right\rceil,&
 r_I&=\frac{\theta}{16N},& C_0&=\frac{2^{20}}{\kappa^2},\\
 K_0&=4h+\frac{256h}{\kappa'},&
 K_b&=72K_0+128h,&
 K_\gamma&=\frac{3\cdot2^{19}}{\kappa^2}K_0+112h,\\
 b_I&=\frac{r_I^2}{8K_b},&
 a_I&=\min\left\{\frac14,\frac{r_I^2}{512h}\right\},&
 t_I&=\left\lceil(192/b_I)^2\right\rceil,\\
 \gamma&=\min\left\{\frac14,\frac{\kappa^2}{2^{20}},
 \frac{r_I^2}{8K_\gamma},\frac{b_I}{222},
 \frac{a_Ib_I}{168\cdot2^{8t_I}}\right\},\\
 v&=\min\left\{\frac12,\frac{r_I^2}{128hC_0},
 \frac{r_I^2}{128h},\frac{\gamma^8}{3}\right\},&
 \alpha&=v^2,&
 \beta&=\min\left\{\frac{\theta}{16},\frac{\alpha^2}{64}\right\},\\
 L&=\max\left\{N,\left\lceil\frac{2}{\kappa\beta}\right\rceil\right\},&
 \mu_I&=\frac{\beta^2}{4L^2},&
 \rho&=\left(\frac{1-\lambda}{1-\kappa}\right)^2,\\
 m_{\rm cut}&=\left\lceil\frac{2}{(1-\rho)\kappa\mu_I}\right\rceil,&
 \delta_*&=\frac{\kappa\mu_I(1-\kappa)^{2m_{\rm cut}}}{2},&
 \mathfrak e(h,\lambda,C,\theta)&=
 \min\left\{1,\frac{\delta_*}{C+1}\right\}.
\end{align*}
\endgroup
All denominators are strictly positive: in particular
$0<\kappa<\lambda<1$ and $0<\rho<1$. Thus
$\mathfrak e$ is a positive rational number obtained by finitely many
rational operations, comparisons, integer powers and ceilings. Its very
large intermediate integers do not affect termination.

\paragraph{The two analytic implications being imported.}
For a unitary tuple $w=(w_1,\ldots,w_h)$ in $M_d(\mathbb C)$, put
\[
 T_w(x)=\frac{x}{2}+\frac1{4h}\sum_j
 (w_jxw_j^*+w_j^*xw_j),\qquad
 \mathcal L_w=I-T_w,\qquad
 \mathcal E_w(x)=\frac1{4h}\sum_j\|[w_j,x]\|_2^2.
\]
All norms in this subsection use the normalized matrix trace.
Corollary~4.3 says that, if $0<z<1$ is rational,
$\theta=z^8/2^{96}$, and
\[
 \langle(\mathcal L_w^2-\lambda\mathcal L_w)y,y\rangle
 \geq-C\varepsilon\quad(\|y\|\leq1),\qquad
 0\leq\varepsilon\leq\mathfrak e(h,\lambda,C,\theta),
\]
then a unital $*$-subalgebra $D\subset M_d(\mathbb C)$ satisfies
\begin{equation}\label{char2merge:eq:char2-approximate-commutant}
 \|x-E_Dx\|_2\leq\sqrt{2/\lambda}\,
 \sqrt{\mathcal E_w(x)}+3z,\qquad
 \max_j\|[w_j,a]\|_2\leq6\sqrt h\,z
\end{equation}
for contractions $x\in M_d(\mathbb C)$ and $a\in D$.
Here $E_D$ denotes the trace-preserving conditional expectation.

For the second implication, fix
\[
 k=2^{-9},\qquad r=2^{-17335},\qquad\eta=2^{-17313}.
\]
Theorem~4.1 says that unital $*$-subalgebras $A,B\subset M_d(\mathbb C)$
and unitaries $U_1,\ldots,U_m$ satisfy
\begin{equation}\label{char2merge:eq:char2-reverse-inclusion}
 U_jAU_j^*\subseteq_{2^{-20}}A
\end{equation}
provided that
\[
 \|x-E_Bx\|_2\leq k^{-1/2}\sqrt{\mathcal E_U(x)}+r,
 \qquad \max_j\|[U_j,b]\|_2\leq r,
 \qquad U_j^*AU_j\subseteq_\eta A
\]
for all contractions $x\in M_d(\mathbb C)$ and $b\in B$.
The notation $D\subseteq_\delta A$ means
$\sup_{x\in D,\,\|x\|\leq1}\|x-E_Ax\|_2\leq\delta$.
Consequently all approximants can be contractions.
These two implications are characteristic-free statements about matrix
algebras. They require neither Hilbert--Schmidt stability of the group
presentation nor correction of the source tuple to a genuine group
representation.

\paragraph{Finite hypotheses supplied by the characteristic-two algebraic checks.}
Let $s$ be a tuple of the $32$ source symmetries and let
\[
 Q_{79}=\sum_{j=1}^{79}\rho_j^*\rho_j,\quad
 q=\operatorname{tr}_dQ_{79}(s),\quad \varepsilon=\sqrt q,\quad
 S=(x_0,y_0,z_{0,0},z_{1,0},z_{2,0},z_{3,0}).
\]
Write $U=(U_{ab})_{a\ne b}$ for the six prescribed integer-action
words, $u=U_{21}$, $h=z_{1,0}$, and $t=Vu$. All these words are unitary.
Sections~\ref{char2merge:sec:steinberg-bridge} and~\ref{char2merge:sec:paid-transfers}, in
particular \eqref{char2merge:eq:st-paid-gap}, \eqref{char2merge:eq:paid-integer-gap}, and
\eqref{char2merge:eq:paid-compression}, establish the following four statements for
this paid presentation:
\begin{enumerate}
 \item The root certificate gives
 \[
 \langle(\mathcal L_S^2-2^{-32}\mathcal L_S)y,y\rangle
 \geq-C_H\varepsilon\qquad(\|y\|\leq1),
 \]
 where $C_H\geq1$ is the positive rational returned by the fixed,
 proved-terminating certificate enumeration.
 \item The integer certificate and paid transfer give
 \[
 \langle(\mathcal L_U^2-\mathcal L_U/144)y,y\rangle
 \geq-2^{30}\varepsilon\qquad(\|y\|\leq1).
 \]
 \item For every $U_{ab}$ and $g\in S$ there is a unitary target
 $p_{ab,g}$ such that
 \[
 \|U_{ab}gU_{ab}^*-p_{ab,g}\|_2\leq360\varepsilon,
 \qquad
 \|p_{ab,g}-w_{ab,g}(S)\|_2\leq16384\varepsilon,
 \]
 where $w_{ab,g}$ is a word of length at most $16$ in $S$.
 \item The two actions $u,t$ have common targets:
 \[
 \|ugu^*-p_g\|_2\leq360\varepsilon,
 \qquad \|tgt^*-p_g\|_2\leq360\varepsilon
 \qquad(g\in S).
 \]
\end{enumerate}
The analytic proof below uses these already established finite checks.
In particular, the integer defect constant is $2^{30}$ for this source,
not an untransferred constant from another presentation.

\paragraph{Lemma (the characteristic-two reflection bound).}
Set
\begin{equation}\label{char2merge:eq:char2-epsilon-star}
 z=2^{-17500},\qquad \theta=z^8/2^{96}=2^{-140096},\qquad
 \varepsilon_*=
 \min\{z,\mathfrak e(6,2^{-32},C_H,\theta),
             \mathfrak e(6,1/144,2^{30},\theta)\}.
\end{equation}
Under the four finite hypotheses above, $\varepsilon<\varepsilon_*$
implies
\begin{equation}\label{char2merge:eq:char2-reflection-bound}
 \|[h,V]\|_2<5\cdot2^{-20}.
\end{equation}

\begin{proof}
Apply \eqref{char2merge:eq:char2-approximate-commutant} first to $S$ and then
to $U$, obtaining algebras $A$ and $B$, respectively. For the root
algebra $A$,
\[
 \sqrt{2/(2^{-32})}=\sqrt{2^{33}}<2^{17},\qquad
 6\sqrt6<16,
\]
so
\begin{equation}\label{char2merge:eq:char2-root-estimates}
 \|x-E_Ax\|_2\leq2^{17}\sqrt{\mathcal E_S(x)}+3z,
 \qquad \max_{g\in S}\|[g,a]\|_2<16z
\end{equation}
for contractions $x$ and $a\in A$.
For the action algebra $B$,
$\sqrt{288}<\sqrt{512}=k^{-1/2}$ and
$3z<r$, $16z<r$. Thus $B$ satisfies exactly the first two
hypotheses of Theorem~4.1. Notice that $k$ belongs to the action
algebra $B$; there is no hypothesis that the root parameter
$2^{-33}$ be at least $k$.

Fix a contraction $a\in A$, $U_{ab}$, and $g\in S$.
Telescoping a word of length at most $16$ in the second estimate
of \eqref{char2merge:eq:char2-root-estimates} gives
$\|[a,w_{ab,g}(S)]\|_2\leq256z$.
Replacing that word by $p_{ab,g}$ costs at most
$32768\varepsilon$ in the commutator. Replacing the target by
$U_{ab}gU_{ab}^*$ costs at most $720\varepsilon$. Therefore
\[
 \|[U_{ab}^*aU_{ab},g]\|_2
 \leq256z+33488\varepsilon
 <33744z<2^{16}z.
\]
Because there are six roots and the energy coefficient is $1/24$,
\[
 \sqrt{\mathcal E_S(U_{ab}^*aU_{ab})}<2^{15}z.
\]
The first estimate in \eqref{char2merge:eq:char2-root-estimates} now gives
\[
 \|U_{ab}^*aU_{ab}-E_A(U_{ab}^*aU_{ab})\|_2
 <(2^{32}+3)z<\eta.
\]
Indeed $\eta/z=2^{187}>2^{32}+3$. All normalization hypotheses
are satisfied, and \eqref{char2merge:eq:char2-reverse-inclusion} follows for
the six $U_{ab}$, in particular for $u$.

Put $c_{\rm cmp}=u^*t$. The common-target estimates imply
\[
 \|[c_{\rm cmp},g]\|_2=\|tgt^*-ugu^*\|_2
 \leq720\varepsilon\quad(g\in S),\qquad
 \sqrt{\mathcal E_S(c_{\rm cmp})}\leq360\varepsilon.
\]
Consequently
\[
 \|c_{\rm cmp}-E_Ac_{\rm cmp}\|_2
 \leq360\cdot2^{17}\varepsilon+3z
 <2^{26}\varepsilon+3z.
\]
Since $E_Ac_{\rm cmp}$ is a contraction, normalization supplies a contraction
$b\in A$ with
\[
 \|uc_{\rm cmp}u^*-b\|_2
 \leq2^{26}\varepsilon+3z+2^{-20}.
\]
But $uc_{\rm cmp}u^*=tu^*=V$, so the root commutator estimate yields
\[
 \|[h,V]\|_2
 \leq16z+2(2^{26}\varepsilon+3z+2^{-20})
 <(2^{27}+22)z+2^{-19}<5\cdot2^{-20}.
\]
For the last numerical comparison,
$(2^{27}+22)z<2^{28}2^{-17500}=2^{-17472}<3\cdot2^{-20}$,
while $2^{-19}=2\cdot2^{-20}$. Every estimate is independent
of the matrix size.
\end{proof}

\paragraph{Corollary (a definite integer source coefficient).}
Define
\begin{equation}\label{char2merge:eq:char2-mu-sym}
 \mu_{\rm sym}=\left\lceil4/\varepsilon_*^2\right\rceil.
\end{equation}
The certificate-termination proof and
\eqref{char2merge:eq:char2-epsilon-star} make this a well-defined finite positive
integer. For every matrix size and every tuple of the $32$ symmetries,
\[
 \operatorname{tr}_d(F_0+\mu_{\rm sym}Q_{79})\geq\frac34,
 \qquad F_0=-1+VhVh+hVhV.
\]
\begin{proof}
For symmetries $h,V$ the exact identity
\[
 \operatorname{tr}_dF_0=1-\|[h,V]\|_2^2
\]
holds, and in all cases $\operatorname{tr}_dF_0\geq-3$.
If $q<\varepsilon_*^2$, the lemma gives
$\operatorname{tr}_dF_0>1-25\cdot2^{-40}>3/4$.
If $q\geq\varepsilon_*^2$, then
$\mu_{\rm sym}q\geq4$ and hence
$\operatorname{tr}_d(F_0+\mu_{\rm sym}Q_{79})\geq1$.
\end{proof}
\subsubsection{The forty-symmetry source bound}
Given any forty-symmetry source tuple, keep the thirty-one common base
coordinates and assign the remaining base coordinate
$Kxy\_00=x0\,y0\,x0$. This is a symmetry, since it is a conjugate of
$y0$. The objective uses only the unchanged coordinates $V,z1\_0$.
The exact paid word transfer proves
\[
 \operatorname{tr}_n Q_{79}
 \leq149\,778\operatorname{tr}_n Q_{80}.
\]
Applying the preceding corollary to the resulting thirty-two-symmetry
base tuple gives
\[
 \operatorname{tr}_n(F_0+\mu_{\rm sym}Q_{79})\geq\frac34.
\]
Since $Q_{80}$ has nonnegative trace and
$149\,778<11\,156\,544\,672$, we obtain
\[
 \operatorname{tr}_n(F_0+\lambda_{80}Q_{80})
 \geq\operatorname{tr}_n(F_0+\mu_{\rm sym}Q_{79})\geq\frac34.
\]

\subsection{Quantitative source reductions}\label{app:214coefficient}
The theorem in Section~\ref{sec:quadratic214} uses the integer \(\lambda\) in Proposition~\ref{thm:source214gap}. Sections~\ref{app:char2bridge} and~\ref{app:char2norm} define \(\mu_{\sym}\) via a rational certificate search and prove termination.

The 32-symmetry, 79-row base system is \path{quadratic_214/f2_compact32_79rows_cubic.json}. The extended source, \path{quadratic_214/AB_full_corner_direct_source.json}, has 40 symmetry coordinates, two unitary coordinates and 82 rows. Substituting \(N=\mathrm{C1}\,\mathrm{B1}\) and \(T=x_0z_{00}\) ensures unitarity and eliminates their binding rows, leaving the 40-symmetry, 80-row system.

By Section~\ref{app:char2norm}, for \(\mu=\mu_{\sym}\) the 80-row system satisfies
\begin{equation}\label{qsrc:source}
 \tr(F_0+\lambda_{80}Q_{80})\ge3/4,
 \qquad\lambda_{80}=11156544672\mu.
\end{equation}
The following reductions transfer this bound to the 37-symmetry, 66-row system.
For an exact symmetry word, \(u^*=\operatorname{rev}(u)\), and a cyclic rotation or inverse of a cap preserves its squared defect norm. In particular,
\begin{equation}\label{qsrc:cuttrace}
 \|u(S)-v(S)\|_2^2=\|a(S)-b(S)\|_2^2
\end{equation}
whenever the two cuts have the certified same cyclic-inverse cap class.

\subsubsection{Quantitative row deletion}\label{q214:domination}
Let $I$ consist of the row indices of the 80-row presentation
\[
 I=\{0,7,10,19,20,26,27,35,38,43\},\qquad
 J=\{0,\ldots,79\}\setminus I.
\]
The associated omitted row names are, respectively,
\begin{center}\small
\begin{tabular}{rl}
\toprule Index&Source row\\\midrule
0&\texttt{comm\_z0\_0\_z1\_0}\\
7&\texttt{comp\_21\_z3\_0}\\
10&\texttt{new\_swap\_x0}\\
19&\texttt{pcov\_W23\_C1}\\
20&\texttt{pcov\_W23\_z1\_0}\\
26&\texttt{tau\_y0\_d0\_0}\\
27&\texttt{tau\_c0\_z0\_0}\\
35&\texttt{frob\_comm\_W23}\\
38&\texttt{torus\_base\_stein\_zx}\\
43&\texttt{torus\_C1\_y}\\\bottomrule
\end{tabular}
\end{center}
All recovery rows 44 through 79 are retained. We prove a bound on
\emph{traces} of the squared penalties on arbitrary exact symmetry tuples:
\begin{equation}\label{q214:dominationtrace}
 \tr Q_{80}(S)\le351563\,\tr Q_{70}(S),\qquad
 Q_{70}=\sum_{j\in J}\rho_j^*\rho_j.
\end{equation}

\subsubsection{Literal rewrite paths and their charges}
For a word $w$ in exact symmetries set $E(w)=\qnorm{w(S)-I}$, and
$e_j=\qnorm{u_j(S)-v_j(S)}$.
Adjacent involution cancellation preserves word evaluation. Cyclic
rotation conjugates by a unitary, and reversal takes an inverse;
therefore these operations preserve $E$. Deleting matching first and
last letters is also a unitary conjugation. It follows that the fixed
canonical free-cyclic normalization preserves $E$.

Each certified rewrite rotates a word, replaces a prefix $l$ by $r$,
and freely/cyclically normalizes. The exact checker verifies that
$l\operatorname{rev}(r)$ is the free-cyclic representative of the cited
source or helper cap. Unitary multiplication and the triangle inequality
then prove
\[
 E(w_{\rm before})\le e_{\rm cited}+E(w_{\rm after}).
\]
Every target and helper path ends at the empty word. The path certificate
contains ten target paths and thirty reachable helper proofs, with 507
literal steps altogether. Dependencies form an acyclic graph of maximum
depth six. Earlier omitted rows and helpers are recursively replaced by
their proved paths. Every terminal dependency is a retained row in $J$.
Thus the exact integer charge matrix $M\in\mathbb Z_{\ge0}^{10\times70}$
computed from those paths satisfies, componentwise,
\begin{equation}\label{q214:charge}
 (e_i)_{i\in I}\le M(e_j)_{j\in J}.
\end{equation}

\subsubsection{An exact charge Gram bound}
For the ordering of $I$ above, the exact Gram matrix $\mathsf{G}=MM^{\mathsf T}$ is
\begingroup\setlength{\arraycolsep}{2.4pt}\small
\[
\mathsf{G}=\begin{pmatrix}
25&4&0&0&20&0&0&0&0&0\\
4&1876&5600&380&22892&692&622&6024&266&1492\\
0&5600&21864&1072&79216&2674&2380&23108&1092&5684\\
0&380&1072&180&4524&116&116&1208&24&316\\
20&22892&79216&4524&303272&9740&8714&84388&3894&20856\\
0&692&2674&116&9740&359&298&2810&133&690\\
0&622&2380&116&8714&298&263&2520&124&620\\
0&6024&23108&1208&84388&2810&2520&24584&1160&6028\\
0&266&1092&24&3894&133&124&1160&91&272\\
0&1492&5684&316&20856&690&620&6028&272&1501
\end{pmatrix}.
\]
\endgroup
Put $h=351562$. The following short integer certificate gives $\mathsf{G}<hI$.
For the displayed positive vector $v$, the last column is exactly
$hv-\mathsf{G}v$:
\begin{center}\small
\begin{tabular}{rrr}\toprule Omitted source row&$v_i$&$(hv-\mathsf{G}v)_i$\\\midrule
0&58&89410\\7&74934&280840\\10&262914&123994\\19&14804&143148\\
20&1000000&202100\\26&32293&186181\\27&28886&70212\\
35&279922&379138\\38&12954&176881\\43&69126&349124\\\bottomrule
\end{tabular}
\end{center}
For every real $x$, the weighted diagonal-dominance identity is
\[
 h\sum_i x_i^2-x^{\mathsf T}\mathsf{G}x
 =\sum_i\frac{(hv-\mathsf{G}v)_i}{v_i}x_i^2
 +\sum_{i<j}\mathsf{G}_{ij}
  \left(\sqrt{v_j/v_i}\,x_i-\sqrt{v_i/v_j}\,x_j\right)^2\ge0.
\]
Every coefficient has the required sign, checked by integer arithmetic.
The source audit independently gives an exact positive rational LDL
factorization of $hI-\mathsf{G}$ as well. Since $M^{\mathsf T}M$ and
$MM^{\mathsf T}$ have the same nonzero eigenvalues,
$\qop M^2\le h$. Squaring \eqref{q214:charge} and adding the retained
residual squares proves \eqref{q214:dominationtrace}.

Consequently the 80-row source assertion \eqref{qsrc:source} implies
\begin{equation}\label{q214:source70}
 \tr\bigl(F_0(S)+\lambda_{70}Q_{70}(S)\bigr)\ge\tfrac34,
 \quad
 \lambda_{70}=351563\lambda_{80}=3922228314522336\,\mu.
\end{equation}
The specified witness annihilates all eighty rows and hence all retained
rows. Only after the graph construction supplies genuine symmetries and a quantitative
bound for $Q_{70}$ will \eqref{q214:source70} be invoked.

\subsubsection{Reconstructing omitted source coordinates}\label{q214:sourcereconstruction}
The forty-coordinate, seventy-row source is now replaced by a source with
$37$ independent symmetry coordinates and initially 68 paid rows. The signal letters
$t,d$ remain, so the objective $F_0=-I+tdtd+dtdt$ is unchanged.

\subsubsection{Exact palindrome substitutions}
The following coordinates of the forty-coordinate source are reconstructed literally:
\begin{center}\small
\begin{tabular}{llllr}
\toprule Target&Source coordinate&Palindrome&Defining row&Other row\\\midrule
\texttt{M}&\texttt{F2\_ladder\_L12}&\texttt{LGL}&74&75\\
\texttt{j}&\texttt{R13}&\texttt{lpl}&57&59\\
\bottomrule\end{tabular}\end{center}

For each listed target $x$ with palindrome $p_x$, put $S_x=p_x(S)$.
An odd palindrome is a conjugate of its middle symmetry by the unitary
word formed by its first half, so each $S_x$ is exactly a symmetry.

\subsubsection{A spectral reconstruction}
The source letter $C=\texttt{Dsign13}$ occurs only in retained rows 32
and 72, namely $Cl=lC$ and $C=hi$. Replace those two paid rows by
\[
 hil=lhi,\qquad hi=ih,
\]
and omit $C$ as an independent coordinate. Given arbitrary values of all
retained symmetries, let $U=S_hS_i$ and define
\[
 S_C=\operatorname{sign}(\operatorname{Re}U),\qquad
 \operatorname{sign}(0)=1.
\]

For $|z|=1$, put $x=|\operatorname{Re}z|$. Then
\[
 |z-\operatorname{sign}(\operatorname{Re}z)|^2=2(1-x),\qquad
 |z^2-1|^2=4(1-x^2),
\]
so their difference in the reverse order is $2(1-x)(1+2x)\ge0$.
Consequently, writing $a=\qnorm{hi-ih}$ and $b=\qnorm{hil-lhi}$,
\[
 \delta=\qnorm{S_C-U}\le\qnorm{U^2-I}=\qnorm{U-U^*}=a.
\]
The reconstructed source residuals satisfy
\[
 e_{72}=\delta\le a,\qquad
 e_{32}=\qnorm{S_Cl-lS_C}\le b+2\delta\le b+2a.
\]
Here all unchanged letters denote their retained input values. Therefore
\[
 e_{72}^2+e_{32}^2\le a^2+(b+2a)^2\le6(a^2+b^2).
\]
The last inequality is exact: the complement of
$\begin{psmallmatrix}5&2\\2&1\end{psmallmatrix}$ in $6I$ is
$\begin{psmallmatrix}1&-2\\-2&5\end{psmallmatrix}$, whose leading minor
and determinant are both one.

All other residuals are their literal substitutions, and the palindrome
binding residuals vanish. Thus, on the reconstructed forty-symmetry tuple,
\[
 \tr Q_{70}\le6\tr Q_{68},\qquad
 \tr Q_{80}\le2109378\tr Q_{68}.
\]
Applying the proved inequality~\eqref{qsrc:source} only after this reconstruction gives
\[
 \tr\bigl(F_0+\lambda_{68}Q_{68}\bigr)\ge\tfrac34,
 \qquad \lambda_{68}=23533369887134016\,\mu.
\]
At the specified witness, $hi=C$ is already a symmetry, so its sign
reconstruction recovers the specified $C$ exactly. The palindrome similarly
recovers every omitted palindrome coordinate. All new rows vanish and the witness objective remains $-1$.

\subsubsection{Two further paid source rows can be omitted}\label{q214:twodeletions}
Rows 2 and 72 of the 80-row presentation are $ik=ki$ and $hi=ih$. We now omit these
two residuals from the 37-symmetry source, leaving 66 paid rows.
The preceding reconstruction of $C$ is still valid: its error is paid
indirectly by the quantitative comparison proved here, rather than
being set equal to a possibly nonsymmetric product.

For a symmetry word $w$, put $\rho(w)=\qnorm{w-I}$.
Free cancellation, cyclic conjugation, and reversal preserve $\rho$.
Replacing a prefix $L$ by $R$ in a word $LS$ changes its evaluation by
exactly $\qnorm{L-R}$, because $S$ is unitary. Thus a literal path of
charged replacements ending at the empty word gives a positive linear
bound on its initial residual. Substitution of either exact palindrome
$j=lpl$ or $M=LGL$ preserves this argument.

The finite certificate checks 30 acyclic helper nodes with 387 literal
steps, also verified after the two palindrome substitutions. The paths
for rows 72 and 2 have respectively 19 and 20 steps, giving 426 literal
steps in all. Binding rows 57 and 74 of the 80-row presentation become identical on both
sides after substitution and contribute zero.

Let $e_r$ be the normalized Hilbert--Schmidt norm of retained residual
$r$. The exact charge vectors $a_{72},a_2$ satisfy
\[
 e_{72}\le a_{72}e,\qquad e_2\le a_2e.
\]
Here the row identifiers refer to the 80-row presentation. Their nonzero entries
are displayed below; every unlisted entry is zero.
\begin{center}\small
\begin{tabular}{rrr@{\qquad}rrr@{\qquad}rrr}
\toprule Row&$a_{72}$&$a_2$&Row&$a_{72}$&$a_2$&Row&$a_{72}$&$a_2$\\\midrule
3&2&0&45&0&752&59&8&0\\
4&0&1&46&0&284&60&0&32\\
14&0&16&47&0&24&62&0&16\\
22&0&2&48&0&368&66&0&8\\
23&0&12&49&0&24&70&6&2\\
29&0&18&50&0&240&73&0&644\\
31&0&20&51&0&240&75&0&560\\
32&3&0&52&0&4&76&0&464\\
33&0&184&53&0&52&77&0&760\\
34&0&8&54&0&60&78&0&300\\
36&0&184&55&8&16&79&0&440\\
37&0&16&56&0&4&&&\\
44&0&256&58&0&2&&&\\
\bottomrule\end{tabular}\end{center}

In particular,
$e_{72}\le2e_3+3e_{32}+8e_{55}+8e_{59}+6e_{70}$.
Writing $M$ for the matrix with these two charge rows, exact arithmetic
gives
\[
 MM^{\mathsf T}=\begin{pmatrix}177&140\\140&2845409\end{pmatrix},
 \qquad
 2845410I-MM^{\mathsf T}
 =\begin{pmatrix}2845233&-140\\-140&1\end{pmatrix}.
\]
The last matrix has positive diagonal entries and determinant
$2825633>0$, hence is positive definite. Squaring the positive
charge bounds and adding the retained residual squares proves
\begin{equation}\label{q214:tworowdomination}
 \tr Q_{68}\le2845411\tr Q_{66}.
\end{equation}

Combining the preceding source transfer with
\eqref{q214:tworowdomination} gives the final source assertion
\begin{equation}\label{q214:reducedsource}
 \tr(F_0+\lambda Q_{66})\ge\tfrac34,
 \qquad\lambda=66962109543919887600576\,\mu.
\end{equation}
Indeed $\lambda=2845411\lambda_{68}$, and the total domination
factor relative to $Q_{80}$ is $6002047364358$.

\subsubsection{The coefficient and the rational search}
Combining the rational search in Section~\ref{app:char2bridge}, the matrix correction in Section~\ref{app:char2norm}, and the preceding transfers gives
\[
11156544672\cdot351563\cdot6\cdot2845411
 =66962109543919887600576.
\]
The witness is proved separately in Section~\ref{app:char2witness}.
\subsection{The characteristic-two source witness}\label{app:char2witness}

\subsubsection{Witness statement}\label{witnessmerge:precise-result-and-scope}

For the 82-row presentation in Section~\ref{app:char2source}, there is a tuple in a finite tracial von Neumann algebra satisfying all
40 symmetry requirements, both unitary requirements, and all 82 literal
source rows. At this tuple

\[
\tau(F_0)=-1,\qquad
\tau((Vd-dV)^*(Vd-dV))=2,\quad d=z1\_0.
\]

After the substitutions \texttt{F2\_torus\_N=C1\ B1} and
\texttt{F2\_steinberg\_P=x0\ z0\_0}, all 80 remaining residuals vanish.
The 37 coordinates in Table~\ref{tab:214lettermap} are selected from these 40 symmetries. The retained-row check below verifies all 66 printed cuts and proves Proposition~\ref{hyp:source214} independently of the finite-matrix estimate.

\subsubsection{The amalgamated group}\label{witnessmerge:the-group-exists}

Set \(R=\mathbb F_2[t_1,t_2,t_3]\),
\(L=\mathbb F_2[t_1^{\pm1},t_2^{\pm1},t_3^{\pm1}]\), and
\(H=\operatorname{EL}_3(R)\). For
\(A\in\operatorname{GL}_3(\mathbb Z)\), define \(\alpha_A(t^n)=t^{An}\),
where exponents are column vectors. This gives an action by ring
automorphisms. Put

\[
K=\operatorname{GL}_3(L)\rtimes_\alpha\operatorname{GL}_3(\mathbb Z),
\quad (M,A)(N,B)=(M\alpha_A(N),AB).
\]

The injection of \(R\) into \(L\) embeds \(H\) into \(K\) as
\(h\mapsto(h,I)\). The map \(\theta(M,A)=(M^{-T},A)\) is an involutive
automorphism: \((MN)^{-T}=M^{-T}N^{-T}\), and inverse transpose commutes
with entrywise ring automorphisms. It preserves \(H\), since
\(e_{ij}(f)\mapsto e_{ji}(f)\) in characteristic two.

Form \(\Omega=K_1*_H K_2\). The common \(\theta\) action extends to
\(\Omega\) by the universal property of the amalgam. Factor exchange
\(\sigma\) is another automorphism, fixes \(H\) pointwise, has order
two, and commutes with diagonal \(\theta\). Therefore
\[
\Gamma=\Omega\rtimes(C_{2,\sigma}\times C_{2,\theta})
\]
is well defined. The normal-form theorem for amalgamated free products
embeds both copies of \(K\), and the semidirect product embeds
\(\Omega\). We use the first copy below. We write \(\vartheta\) for the
group element implementing \(\theta\).

\subsubsection{The 42 coordinates}\label{witnessmerge:all-42-coordinates-explicitly}

Write \(e_{ij}(f)=I+fE_{ij}\), let \(s_j\) negate exponent coordinate
\(j\), and let \(u_{ij}=I+E_{ij}\) over \(\mathbb Z\). Pure exponent
matrices are embedded as \((I,A)\). Let \(P_{ij}\) interchange ordinary
matrix coordinates \(i,j\). Define

\[
S_j=(I,s_j),\quad R_{ij}=(I,u_{ij}s_j),\quad
T_{ij}=(I,(I-E_{ij})s_j),\quad W_{ij}=R_{ij}R_{ji}R_{ij}.
\]

Here \(W_{ij}\) is the pure exponent permutation interchanging \(i,j\).
Put

\[
\begin{array}{lll}
x0=e_{12}(1),&y0=e_{23}(1),&c0=e_{13}(1),\\
d0\_0=e_{32}(1),&z0\_0=e_{31}(1),&zj\_0=e_{31}(t_j)\ (j=1,2,3),\\
p12\_0=e_{31}(t_1t_2),&H1\_0=e_{31}(t_1^{-1}),&M1=R_{32}R_{21}R_{32},\\
B1=(\operatorname{diag}(1,1,t_1),s_1),&&\\
C1=(\operatorname{diag}(t_1^{-1},1,1),s_1),&&\\
C2=(\operatorname{diag}(t_2^{-1},1,1),s_2),&&\\
Fgraph=P_{13}\vartheta,&Ggraph=P_{12}\vartheta,&Tau=\vartheta,\\
Dsign13=S_1S_3,&Jhybrid=C2R_{13}.&
\end{array}
\]

The pure-action coordinates of the source are
\(S_1,S_3,R_{13},R_{21},R_{31},R_{32},T_{31},T_{21},W_{12},W_{13},W_{23}\).

Define

\[
\mathcal B=W_{12}R_{12}S_2,\qquad V=\mathcal B\sigma\mathcal B^{-1}.
\]

The remaining eleven coordinates are the following literal
words in the amalgamated group. The short names $A,B,M,h,N,P$ in this appendix are local source-coordinate names:

\[
\begin{aligned}
A:=F2\_A&=Fgraph\,W_{12},&
B:=F2\_B&=x0\,W_{23},\\
M:=F2\_M&=V\,Ggraph,&
L_3:=F2\_L3&=ABABA,\\
E_2:=F2\_E2&=M\,Fgraph\,M,&
w:=F2\_w13&=c0\,z0\_0\,c0,\\
N:=F2\_torus\_N&=C1\,B1,&
P:=F2\_steinberg\_P&=x0\,z0\_0,\\
L_6:=F2\_ladder\_L6&=L_3BL_3,&
L_{12}:=F2\_ladder\_L12&=L_6BL_6,\\
h:=F2\_ladder\_h&=c0\,x0.&&
\end{aligned}
\]

These specify all 42 witness coordinates. All
except \(N,P\) are involutions. For the roots this follows from
characteristic two. For the pure actions and diagonal-action elements,
it follows from their displayed formulas. \(M1\) is a conjugate of
\(R_{21}\); \(Dsign13\) is a product of commuting involutions;
\(Jhybrid\) is a product of commuting involutions. Constant permutations
are fixed by inverse transpose, so \(Fgraph,Ggraph\) are involutions.
They commute with the pure action elements; \(V\) commutes with both
graph elements, as proved below. Hence \(A,B,M\) are involutions. The
words \(L_3,E_2,w,L_6,L_{12}\) are palindromes in involutions. Finally
\(c0\) and \(x0\) commute, so \(h\) is an involution. Every coordinate,
including \(N,P\), is a group unitary.

\subsubsection{Factor-swap identities}\label{witnessmerge:the-identities-involving-the-factor-swap}

The exponent matrix of \(\mathcal B\), and its inverse, are

\[
\begin{pmatrix}0&1&0\\1&1&0\\0&0&1\end{pmatrix},\qquad
\begin{pmatrix}-1&1&0\\1&0&0\\0&0&1\end{pmatrix}.
\]

Consequently

\[
\begin{aligned}
\mathcal B^{-1}x0\mathcal B&=e_{12}(1),\\
\mathcal B^{-1}z2\_0\mathcal B&=e_{31}(t_1),\\
\mathcal B^{-1}p12\_0\mathcal B&=e_{31}(t_2),\\
\mathcal B^{-1}z3\_0\mathcal B&=e_{31}(t_3).
\end{aligned}
\]

All four right sides belong to \(H\), hence are fixed by \(\sigma\).
Thus \(V\) commutes with the four corresponding source coordinates.
Every constant matrix in \(\operatorname{GL}_3(\mathbb F_2)\) belongs to
\(\operatorname{EL}_3(\mathbb F_2)\subset H\): Gaussian elimination uses
only elementary row additions and row swaps, and a swap is
\(e_{ij}(1)e_{ji}(1)e_{ij}(1)\); no nontrivial row scaling is needed
over \(\mathbb F_2\). Thus \(V\) commutes with all constant matrices and
with \(\vartheta\), and therefore with \(Fgraph,Ggraph\). This proves
every basic swap row and justifies all occurrences of \(V\) in the
recovery rows.

The swap-containing recoveries reduce to

\[
\begin{aligned}
M c0 M&=d0\_0,&M y0 M&=z0\_0,\\
E_2=M Fgraph M&=Ggraph\,Fgraph\,Ggraph=P_{23}\vartheta,\\
Ggraph\,M&=V,&A E_2 A&=Ggraph.
\end{aligned}
\]

No step replaces \(\sigma\) by a matrix in one factor.

\subsubsection{Exact row verification}\label{witnessmerge:verification-of-all-rows-without-numerical-substitutions}

Within \(K\rtimes\langle\vartheta\rangle\), use triples \((M,A,e)\),
\(e\in\{0,1\}\), with the exact multiplication law

\[
(M,A,e)(N,B,f)=\bigl(M\alpha_A(\theta^e(N)),AB,e+f\pmod2\bigr).
\]

Substitute the preceding coordinate words into each of the 82 literal
rows. Seventy-two rows contain no \(V\) after expansion and reduce to
equality by this multiplication law over the formal Laurent polynomial
ring. The other ten are covered by the preceding fixed-subgroup
identities. The accompanying verification implements these algebraic
operations and checks all 40 involution identities and all 82 rows.

The verification replaces \(VgV\) by \(g\) only after proving \(\mathcal B^{-1}g\mathcal B\in H\rtimes\langle\vartheta\rangle\). It accepts only (i) a constant invertible matrix times \(\vartheta^e\), or

(ii) a single elementary matrix with a polynomial entry times
\(\vartheta^e\), always with trivial exponent action. Each is explicitly
in the indicated fixed subgroup. Every accepted pair is recorded in the verification output. This avoids any unsolved general
membership problem for \(H\).

Formal Laurent polynomials are stored as canonically sorted finite sums
of integer-exponent monomials with coefficients in \(\mathbb F_2\).
Addition, multiplication, substitutions, determinants and inverses are
exact; inverse computations also check both products against the
identity. There are no floating point computations, finite substitutions
for the variables, or finite-dimensional surrogate for the amalgam.
This procedure checks the 40 symmetry identities, the 82 source rows, and
the required negative Laurent root. The nontrivial ladder identities
follow from

\[
L_3=(y0\,x0\,Fgraph)W_{23},\quad
L_6=(x0\,y0\,x0)W_{23},\quad
L_{12}=(c0\,x0)W_{23},\quad h=d0\_0\,c0\,d0\_0=c0\,x0.
\]

\subsubsection{The retained 37 coordinates and 66 printed cuts}
The 37-letter dictionary in Table~\ref{tab:214lettermap} is the restriction of the coordinate assignment above; its coordinate names agree entry by entry with the dictionary used in the witness verification. The 37 retained coordinates are obtained from the 42 by omitting the unitary coordinates $N,P$ in the local notation above, and omitting $Dsign13$, $R_{13}$, and $F2\_ladder\_L12$. Thus the retained group elements are involutions.

For the 66 pairs $(u_r,v_r)$ in Section~\ref{app:214data}, substitute the coordinate dictionary and apply the preceding group-word algorithm to $u_rv_r^{-1}$. The row labels are distinct, and the table agrees entry by entry with the coordinate dictionary.

\subsubsection{The negative trace}\label{witnessmerge:the-objective-has-trace-minus-one}

The same exact exponent calculation yields
\[
q=\mathcal B^{-1}z1\_0\mathcal B=e_{31}(t_1^{-1}t_2).
\]
This lies outside \(H\): every element of \(H\) has entries in \(R\),
whereas its \((3,1)\) entry is the Laurent monomial
\(t_1^{-1}t_2\notin R\).
Writing \(q_i\) for \(q\) in factor \(i\),
\[
\mathcal B^{-1}(VdVd)\mathcal B=q_2q_1,\qquad
\mathcal B^{-1}(dVdV)\mathcal B=q_1q_2.
\]
Both are reduced length-two words in \(K_1*_H K_2\), since each syllable
is outside \(H\). By amalgam normal form neither is the identity. Here is the normal-form fact needed in this argument. Choose in each
$K_i$ a set of representatives for the sets $gH$, containing $1$, so every
$g\in K_i$ has a unique decomposition $g=th$ with $t$ a representative
and $h\in H$. Consider
formal strings $t_1\cdots t_m h$, where $h\in H$, every $t_j$ is a
nonidentity transversal representative, and successive $t_j$ belong to
different factors. Right multiplication by an element of $K_i$ is defined
by combining the trailing $h$ with that element, decomposing the result
as a representative times an element of $H$, and merging the last syllable
when it belongs to the same factor. If the representative becomes $1$,
continue one place to the left. Uniqueness of each coset decomposition
shows that these operations obey multiplication within each factor;
the operations from either copy of $H$ agree. They therefore define an
action of the amalgam on these strings. Starting with the empty string,
an alternating product of elements outside $H$ has a nonempty resulting
string: when each new factor is different, multiplying its element by
the trailing element of $H$ cannot put it in $H$. Such a product cannot
be the identity. This also proves that the factor embeddings used above
are injective.

Let \(\mathcal M=L(\Gamma)\), the complex group von Neumann algebra in
the left regular representation. Its canonical faithful normal
normalized trace satisfies \(\tau(\lambda(g))=0\) for every \(g\ne1\).
The specified involutions give self-adjoint
unitary operators, and all source rows vanish as operator identities.
Therefore

\[
\tau(-I+VdVd+dVdV)=-1.
\]

Since \(V^2=d^2=I\),

\[
(Vd-dV)^*(Vd-dV)=2I-dVdV-VdVd,
\]

so its trace is exactly 2. Both defining unitary recovery rows vanish
identically under \(N=C1B1\), \(P=x0z0\_0\), and the other 80 residuals
remain zero. Together with the retained-cut verification above, this proves the specified 37-symmetry source witness.

\section{Operator estimates for the quadratic construction}\label{app:214proof}
This appendix supplies the graph approximation and source-comparison arguments used in Lemmas~\ref{lem:214sourceoutput} and~\ref{lem:214elimination}. All norms are ambient normalized \(2\)-norms, and all finite incidences refer to the fixed manifest in \hyperref[sec:verification]{the data and verification statement}. The source gap is invoked only after sign reconstruction. In this appendix \(\mathcal R=\mathsf R\), \(\mathcal H=\mathsf H\), \(c=c_{214}\), and \(G=G_{214}\). A finite path certificate proves the stated word identity or positive charge sum by the explicitly described rewrite rule.

\subsection{Spectral energy and defect control}
\begin{lemma}\label{q214:spectrallemma}
For all effects, $e\ge0$ and
\[
 \tr(P_0-P_0^2)\le2e,\qquad
 \tr(Q_i-Q_i^2)\le8e\quad(i=1,2).
\]
For projections, with $H=Q_1+Q_2$ and $p(t)=(2t^2-t)/3$,
\begin{align}
 \qnorm{P_0-p(H)}^2&\le4e,\label{q214:topstable}\\
 \qnorm{Q_1Q_2Q_1-Q_1/4}&\le\tfrac32\sqrt e,\qquad
 |\tr Q_1-\tr Q_2|\le8e.\label{q214:calstable}
\end{align}
\end{lemma}
\begin{proof}
Decompose the two input projections into principal-angle blocks. On a
two-dimensional block with cosine $t$, $H$ has eigenvalues $1+t,1-t$.
If $p_+,p_-$ are the corresponding diagonal entries of $P_0$, the
energy contribution is
\[
 t(1-p_++p_-)+(t-1/2)^2.
\]
For the four endpoint pairs $(0,0),(1,0),(1,1),(0,1)$, the differences
$4e-\qnorm{P_0-p(H)}^2$ are respectively
\begin{gather*}
 (7+10t^2-8t^4)/9,\qquad
 \tfrac89(2-t^2-t)(t-1/2)^2,\\
 (1+34t^2-8t^4)/9,\qquad
 (4+18t+22t^2-8t^4)/9.
\end{gather*}
They are nonnegative on $[0,1]$. The trace expansion is affine in
$p_+,p_-$ because $\tr(P_0^2)=\tr P_0$ globally. One-dimensional blocks with exactly one input equal to
one have $e=1/8$ and squared-error contribution at most $4/9$;
common-zero and common-one blocks give differences $3p$ and $1-p$.
This proves \eqref{q214:topstable} and projection positivity. Separate
spectral averaging proves positivity for effects.

For the first defect bound, subtract $\tr(P_0-P_0^2)/2$ from $e$,
average the inputs to projections, and minimize the remaining convex
quadratic by clipping $H-I/2$ to $[0,1]$. The block minimum is zero
for $t\le1/2$ and $(t-1/2)^2$ otherwise. Scalar blocks have minima
$0,0,1/4$. For the second, subtract $\tr(Q_1-Q_1^2)/8$, average
$P_0,Q_2$ to projections and clip $4(P_0-Q_2)$. If $s$ is their
principal-angle sine, the minimum is $1/8-s^2$ for $s\le1/4$ and
$(s-1/2)^2$ otherwise. Scalar blocks have minima zero or $1/8$.
The same holds for $Q_2$ by symmetry. Pinching into a Hamiltonian
eigenbasis preserves the linear trace part and decreases the positive
quadratic part, justifying these scalar minimizations.

Finally, on an angle block the squared compression defect is
$(t^2-1/4)^2\le(9/4)(t-1/2)^2$. A one-input block contributes $1/16$
against energy $1/8$, and a common-one block contributes $9/16$
against minimum energy $1/4$. Rank differences occur only on unmatched
one-input blocks, each costing at least $1/8$. This proves
\eqref{q214:calstable}.
\end{proof}

\subsection{A coercive anchor}
\begin{lemma}\label{q214:coercive}
For arbitrary effects, the anchor satisfies
\[
 \mathcal H\ge0,\quad
 \mathcal H\ge6\tr(E-E^2),\quad
 \mathcal H\ge12\tr(Q_0-Q_0^2),\quad
 \mathcal H\ge20\tr(K-K^2).
\]
\end{lemma}
\begin{proof}
To bound a chosen defect, average the other two variables to projections.
The chosen effect $X$ then appears as $\tr(cX^2+AX)$ plus a constant,
with $c=6,12,20$. Pinching and clipping reduce the minimum to scalar
quadratic minimization. On a generic principal-angle block, the
Hamiltonian eigenvalues for the free variables $E,Q_0,K$
are respectively $3\pm r$, $4\pm r$, $-5\pm r$, with radii
\[
 r^2=729-720s,\quad 1024-960s,\quad 25+1200s.
\]
Here $s$ is the squared cosine of the fixed projection pair. The exact
minima after subtracting the indicated defect are
\[
 \begin{array}{c|c|c}
 \text{free variable}&\text{interval}&\text{minimum}\\\hline
 E&3\le r\le9&(r-9)^2/72\\
 E&9\le r\le27&(r-9)^2/18\\
 Q_0&8\le r\le16&(r-16)^2/96\\
 Q_0&16\le r\le32&(r-16)^2/32\\
 K&5\le r\le15&(750-2r^2)/400\\
 K&15\le r\le25&3(r-25)^2/400\\
 K&25\le r\le35&(r-25)^2/50
 \end{array}
\]
All are nonnegative. For the scalar fixed-projection pairs in order
$00,01,10,11$, the minima are $(0,6,12,1/2)$ for the $E$ bound,
$(0,2,6,2/3)$ for $Q_0$, and $(0,3/4,1,2)$ for $K$.
For $\mathcal H$ alone, minimizing over $K$ gives
$(r-25)^2/50\ge0$ on the generic blocks and $(0,2,6,2)$ on scalar
blocks. These calculations and spectral averaging prove the claims.
\end{proof}

\begin{lemma}\label{q214:anchorstable}
For effects $E,Q_0,K$, there are orthogonal equal-rank projections
$B_1,B_2$, a balanced graph projection $Q^0$ between them and its
canonical rational companion $K^0$, such that
\[
 \qnorm{E-B_1}\le3\sqrt{\mathcal H},\quad
 \qnorm{Q_0-Q^0}\le3\sqrt{\mathcal H},\quad
 \qnorm{K-K^0}\le6\sqrt{\mathcal H},\quad
 4Q^0+3B_1-5K^0=B_1+B_2.
\]
The threshold projection of $Q_0$ is within $2\sqrt{\mathcal H}$ of
$Q^0$.
\end{lemma}
\begin{proof}
First assume $E,Q_0,K$ are projections. On a noncommuting $E,Q_0$
block write $s$ for the squared cosine and $r=\sqrt{25+1200s}$.
The $K$-Hamiltonian $30I-40Q_0-30E$ has eigenvalues $-5\pm r$,
and its minimum is
\[
 h_{\min}=13+24s-r=(r-25)^2/50.
\]
Retain these blocks, leave $E$ unchanged, and balance $Q_0$ at angle
$\pi/4$ with the same phase. Discard the nonzero scalar endpoints.
Their minima are $6,2,2$, so each discarded nonzero $E$ or $Q_0$
costs at most $h_{\min}/2$. On a retained block,
\[
 \lVert Q_0-Q^0\rVert_{\rm HS}^2=1-2\sqrt{s(1-s)}\le4(s-1/2)^2,
 \qquad (s-1/2)^2\le h_{\min}/8.
\]
The latter follows from $(r+25)^2\le3600$. Consequently
$\qnorm{E-B_1}^2,\qnorm{Q_0-Q^0}^2\le\mathcal H/2$.

Let $P_-$ be the negative-Hamiltonian projection on retained blocks,
zero on scalar blocks. If $r\ge15$, both mismatch energy gaps are
at least ten. If $r\le15$, $h_{\min}\ge2$ bounds the local squared
projection distance, which is at most two. Scalar blocks also obey
the bound, including common-zero blocks where the energy is
$30\tr K$. Since $\qnorm{K-P_-}^2=\tr K+\tr P_--2\tr(KP_-)$, these
block estimates include off-block entries and give
$\qnorm{K-P_-}^2\le\mathcal H$.

The cross-eigenvalue gaps to the balanced Hamiltonian are $r+25\ge30$.
Comparing the two orthogonal off-diagonal Hilbert--Schmidt blocks gives
\[
 \qnorm{P_--K^0}\le\frac{40}{30}\qnorm{Q_0-Q^0},\qquad
 \qnorm{K-K^0}\le\left(1+\frac4{3\sqrt2}\right)\sqrt{\mathcal H}
                                  \le2\sqrt{\mathcal H}.
\]
On discarded scalar blocks both comparison projections are zero.

For effects, threshold all three variables. Lemma~\ref{q214:coercive}
bounds their rounding errors by $\sqrt{\mathcal H/6}$,
$\sqrt{\mathcal H/12}$ and $\sqrt{\mathcal H/20}$.
For an effect $Y$, $2Y-\Pi_Y$ is an effect and
$Y=(\Pi_Y+(2Y-\Pi_Y))/2$. Multiaffinity and nonnegativity imply
$\mathcal H(\Pi_E,\Pi_Q,\Pi_K)\le8\mathcal H(E,Q_0,K)$.
The projection-case bounds and the triangle inequality now give the
stated constants. Direct $2\times2$ multiplication gives the last
identity for the balanced pair.
\end{proof}

\subsection{The effect budget, partition and support}
Let $\mathcal R$ be the trace budget on arbitrary effects:
\begin{equation}\label{q214:budget}
 \begin{aligned}
 \mathcal R={}&\sum_k\tr(\ell_k(Y)^2)
 +\sum_{(u,v)\in\mathcal O}\tr(Y_uY_v+Y_vY_u)\\
 &+\sum_T4e(Y_{P_T},Y_{Q_1},Y_{Q_2})+\mathcal H
 +2\sum_{g\in\mathcal G}\mathcal E_g.
 \end{aligned}
\end{equation}
Every summand is nonnegative, using Lemma~\ref{q214:commutationenergy}
for the last sum. Apply Lemma~\ref{q214:anchorstable} and set
$C=I-B_1-B_2$. With $Z=4Q_0+3E-5K$,
\[
 \qnorm{Z-B_1-B_2}\le51\sqrt{\mathcal H},\qquad
 \ell_0=2(I-E_3-Z),\qquad
 \qnorm{E_3-C}\le52\sqrt{\mathcal R}.
\]
Thus an exact partition is obtained without asserting that the defect
$\tr(E_3-E_3^2)$ is controlled by $\mathcal R$.

Threshold every input $Q$ and triangle variable $P_T$. Lemma~\ref{q214:spectrallemma}
gives errors at most $\sqrt{2\mathcal R}$ and
$\sqrt{\mathcal R/2}$, respectively. The rounded energy is at most
$8e_T\le2\mathcal R$. Every rounded alias row has norm at most
\[
 (1+6/\sqrt2+8\sqrt2)\sqrt{\mathcal R}\le17\sqrt{\mathcal R}.
\]
The coordinate terms cancel identically from these rows.

For effects $A,B$, $\qnorm{AB}^2=\tr(A^2B^2)\le\tr(AB)$.
Transferring support in norm, rather than first comparing traces, gives
\begin{align*}
 \qnorm{\Pi_Q C}&\le(\sqrt2+1/\sqrt2+52)\sqrt{\mathcal R}
                                     \le55\sqrt{\mathcal R}\quad(12),\\
 \qnorm{\Pi_Q B_1}&\le(\sqrt2+1/\sqrt2+3)\sqrt{\mathcal R}
                                     \le6\sqrt{\mathcal R}\quad(23).
\end{align*}
The anchor $12$ input has no support pair; its leakage into $C$ is
bounded directly by its canonical-anchor distance $2\sqrt{\mathcal R}$,
so it obeys the same displayed $12$ bound.
Set $\epsilon=60\sqrt{\mathcal R}$. Alias residuals and forbidden
support leakages are at most $\epsilon$. The rounded reference12
projection is within $\epsilon$ of its canonical graph. By
\eqref{q214:calstable}, each retained two-input spectral triangle has compression defect at most
$\epsilon$ and rank mismatch at most $\epsilon^2$.

\subsection{Rectangular calibration}
\begin{lemma}\label{q214:rectangular}
Let $B,C$ be orthogonal projections and $Q\le B+C$ a projection. Put
$\delta=\qnorm{QBQ-Q/2}$ and $k=|\tr Q-\tr B|$.
There is a balanced graph projection $Q^0$ between equal-rank
subcorners $B_0\le B,C_0\le C$ such that
\[
 \qnorm{Q-Q^0}\le2\delta,\qquad \tr(B-B_0)\le4\delta^2+k.
\]
If $\operatorname{rank}C\ge\operatorname{rank}B$, this graph extends
onto all of $B$, giving
\[
 \qnorm{Q-Q^{\rm can}}\le2\delta+\sqrt{4\delta^2+k}.
\]
\end{lemma}
\begin{proof}
On each noncommuting block of $B,Q$, balance the graph while preserving
its phase. If the squared cosine is $s$, its squared error is
$1-2\sqrt{s(1-s)}\le4(s-1/2)^2$. Discard scalar $Q=1$ blocks;
their compression-defect square is $1/4$, while their squared error is
one. A scalar $B=1,Q=0$ block causes no $Q$ error.
Writing $n_{11},n_{10},n_{01}$ for the normalized scalar block sizes,
the discarded $B$ trace is
\[
 n_{11}+n_{10}\le n_{11}+n_{01}+|\tr Q-\tr B|\le4\delta^2+k.
\]
If $C$ has enough rank, extend using unused directions. The added
graph projection has norm $\sqrt{\tr(B-B_0)}$, proving the last bound.
\end{proof}
For any projection $Q$, compression to a projection $D$ followed by
thresholding costs at most $3\qnorm{Q(I-D)}$. Indeed compression costs
twice that norm, and the compressed effect's defect trace is at most
its square. Also, for projections $A,B$,
$|\tr A-\tr B|\le\qnorm{A-B}^2$.

\subsection{Reference projections and input approximations}\label{q214:referencefree}
The anchor and retained reference spectral triangle initialize the
small corners. Its inputs are reference $12$ and reference $23$.
Compress the latter rounded projection to $B_2+C$, at cost at most
$3\epsilon$.
Against the canonical reference $12$ projection the estimates are
\[
 \delta\le2(1+1+(9/4)3)\epsilon\le18\epsilon,\qquad
 k\le(1+1+9)\epsilon^2=11\epsilon^2.
\]
Lemma~\ref{q214:rectangular} discards at most
$(4\cdot18^2+11)\epsilon^2\le(37\epsilon)^2$ from $B_2$.
Transport the retained $B_2$ corner through the reference $12$ unitary
to obtain equal-rank $F_1\le B_1,F_2\le B_2$.
Under this unitary identification of $F_2$ with $F_1$, write the
canonical reference-$23$ coisometry as $A_0:C\to F_1$.
Its initial projection lies in $C$.
Hence $\operatorname{rank}C\ge\operatorname{rank}F_1$,
$A_0A_0^*=I$ on $F_1$, and $\tr F_1\le1/3$.
The approximation errors for the reference $12$ and $23$ inputs are at most
$38\epsilon$ and $39\epsilon$; each small coordinate moves by at
most $37\epsilon$.

\subsubsection{Order of input construction}
Approximate every ordinarily calibrated right input using the canonical
reference $12$, and every ordinarily calibrated left input using the
canonical reference $23$. Their errors are at most $2600\epsilon$
and $4700\epsilon$, respectively. The larger uniform bound
$C_0\epsilon=340000\epsilon$ is retained for the subsequent forest.
Use both $t$ calibrations and the ordinary left $h$ calibration.

Here is a uniform calculation for each subsequent approximation. Suppose the
canonical opposite input has error $a\epsilon$ and the supported
projection approximation costs $s\epsilon$. The support costs are $s=114$
on the $23$ face and $s=225$ on the $12$ face. For the supported
target $Q'$ and the opposite canonical graph $R^0$,
\[
 Q'R^0Q'=\tfrac12Q'F_2Q'.
\]
Support compression removes the opposite outer corner, including a
varying initial projection in $C$. The compression defect and rank
mismatch therefore satisfy
\[
 \delta\le2(1+a+9s/4)\epsilon,\qquad
 k\le(a^2+1+s^2)\epsilon^2.
\]
The rank estimate uses squared projection distances. Given integer
bounds $D\ge2(1+a+9s/4)$ and $K\ge a^2+1+s^2$,
Lemma~\ref{q214:rectangular} bounds the total approximation error by
\[
 \bigl(s+2D+\sqrt{4D^2+K}\bigr)\epsilon.
\]
The following table gives valid integer upper bounds $b\epsilon$.
Every row can be checked by $b-s-2D\ge0$ and
$(b-s-2D)^2\ge4D^2+K$, without floating arithmetic.

\begin{center}\small
\begin{tabular}{lrrrrr}
\toprule Input&$a$&$s$&$D$&$K$&$b$\\\midrule
Ordinary right&38&114&600&14441&2600\\
Ordinary left&39&225&1100&52147&4700\\
\bottomrule\end{tabular}\end{center}

The table's two rows use only the already constructed reference inputs.
After these approximations, the forest below supplies a strict
parent-before-child order for every exceptional input. Source signs,
word reconstruction and inverse-output aliases are established later.

\subsubsection{A finite acyclic calibration forest}\label{q214:forest}\label{q214:neighborconstruction}

Suppose a parent's canonicalization error is at most $a\epsilon$.
For either face the preceding estimates imply
\[
 \delta/\epsilon\le2a+1015,\qquad
 k/\epsilon^2\le a^2+50626.
\]
The target's total graph approximation error is therefore at most
\[
 \left[225+2(2a+1015)
       +\sqrt{4(2a+1015)^2+a^2+50626}\right]\epsilon.
\]
This is bounded by $(10a+10000)\epsilon$ for every $a\ge0$:
the remaining square-root allowance is $6a+7745\ge0$, and its
squared surplus is
\[
 (6a+7745)^2-\bigl(4(2a+1015)^2+a^2+50626\bigr)
 =19a^2+76700a+55813499>0.
\]
Thus the dimension-independent recurrence
\[
 C_0=340000,\qquad C_{d+1}=10C_d+10000
\]
bounds every depth-$d$ target. It is global in $\mathcal R$; at
$\mathcal R=0$ all errors vanish, and rank-zero corners require no
division or different convention. The maximum depth is
$3$, so
\[
 C_1=3410000,\qquad L=C_{3}=341110000.
\]

The $45$ exceptional inputs comprise $19$ singleton roles and $26$ roles
of the $13$ paired classes.
The table records the parent-before-child order. A parent not
listed as a target is an ordinary input with its calibration retained.
\begingroup\small
\begin{longtable}{@{}llllrl@{}}
\caption{Parent-before-child input calibration forest.}\label{tab:214forestdata}\\
\toprule Class&Face&Target product&Parent&Depth&Type\\\midrule
\endfirsthead
\toprule Class&Face&Target product&Parent&Depth&Type\\\midrule
\endhead
$c$&12&$c\cdot wx$&$wx_{23}$&1&pair\\
$d$&12&$d\cdot h$&$h_{23}$&1&pair\\
$e$&12&$e\cdot q$&$q_{23}$&1&single\\
$l$&12&$l\cdot h$&$h_{23}$&1&single\\
$p$&12&$p\cdot Dp$&$Dp_{23}$&1&pair\\
$u$&12&$u\cdot wx$&$wx_{23}$&1&pair\\
$w$&12&$w\cdot z$&$z_{23}$&1&single\\
$G$&12&$G\cdot LG$&$LG_{23}$&1&pair\\
$J$&12&$J\cdot z$&$z_{23}$&1&pair\\
$f$&23&$q\cdot f$&$q_{12}$&1&single\\
$i$&23&$h\cdot i$&$h_{12}$&1&single\\
$k$&23&$h\cdot k$&$h_{12}$&1&pair\\
$n$&23&$h\cdot n$&$h_{12}$&1&single\\
$p$&23&$h\cdot p$&$h_{12}$&1&pair\\
$r$&23&$h\cdot r$&$h_{12}$&1&single\\
$u$&23&$wx\cdot u$&$wx_{12}$&1&pair\\
$v$&23&$N\cdot v$&$N_{12}$&1&pair\\
$x$&23&$z\cdot x$&$z_{12}$&1&single\\
$F$&23&$q\cdot F$&$q_{12}$&1&pair\\
$c$&23&$wx\cdot c$&$wx_{12}$&1&pair\\
$d$&23&$t\cdot d$&$t_{12}$&1&pair\\
$D$&12&$D\cdot p$&$p_{23}$&2&pair\\
$a$&12&$a\cdot u$&$u_{23}$&2&single\\
$g$&12&$g\cdot k$&$k_{23}$&2&single\\
$k$&12&$k\cdot k$&$k_{23}$&2&pair\\
$m$&12&$m\cdot k$&$k_{23}$&2&pair\\
$s$&12&$s\cdot n$&$n_{23}$&2&single\\
$v$&12&$v\cdot v$&$v_{23}$&2&pair\\
$y$&12&$y\cdot k$&$k_{23}$&2&pair\\
$F$&12&$F\cdot p$&$p_{23}$&2&pair\\
$H$&12&$H\cdot F$&$F_{23}$&2&single\\
$K$&12&$K\cdot u$&$u_{23}$&2&single\\
$L$&12&$L\cdot u$&$u_{23}$&2&single\\
$A$&23&$J\cdot A$&$J_{12}$&2&single\\
$y$&23&$l\cdot y$&$l_{12}$&2&pair\\
$D$&23&$l\cdot D$&$l_{12}$&2&pair\\
$G$&23&$G\cdot G$&$G_{12}$&2&pair\\
$J$&23&$c\cdot J$&$c_{12}$&2&pair\\
$b$&23&$J\cdot b$&$J_{12}$&2&single\\
$E$&12&$E\cdot G$&$G_{23}$&3&single\\
$I$&12&$I\cdot G$&$G_{23}$&3&pair\\
$m$&23&$m\cdot m$&$m_{12}$&3&pair\\
$o$&23&$F\cdot o$&$F_{12}$&3&single\\
$B$&23&$K\cdot B$&$K_{12}$&3&single\\
$I$&23&$F\cdot I$&$F_{12}$&3&pair\\
\bottomrule
\end{longtable}
\endgroup

The active depth distribution is 21 at depth 1, 18 at depth 2, 6 at depth 3.
Each parent input is constructed before its child. The selected inverse triangles are used only as spectral geometric controls at this stage; their inverse-output alias identities and source signs are established later.

All physical inputs retain a spectral projection-defect bound, with
the anchor also controlling reference $12$. The uniform error is
$L\epsilon=341110000\epsilon$. Choose
\[
 \begin{aligned}
 T&=6L=2046660000,&V&=23L=7845530000,\\
 G_0&=100L=34111000000,&G&=200L=68222000000.
 \end{aligned}
\]
The auxiliary scale $\gamma_0=G_0\epsilon$ covers output aliases.
The final scale $\gamma=2\gamma_0$ covers all paired leaf interfaces used after sign reconstruction.
These choices cover the spectral output and full virtual-output
estimates below, as well as
$8\sqrt2L+(4\sqrt2+16)/60<G/7$ for primitive and commutator meets.

\subsection{Spectral outputs and decoder interfaces}
Identify the two small corners using the reference12 unitary. For a
canonical triangle with unitary $U$ and coisometry $W:C\to F_1$,
$WW^*=I$ on $F_1$. Its top projection is the rectangular frame
\[
 P_{\rm top}=\frac16
 \begin{pmatrix}U\\2I\\W^*\end{pmatrix}
 \begin{pmatrix}U^*&2I&W\end{pmatrix}.
\]
The frame has squared norm $6I$. The sum of the two input projections
has norm at most $3/2$. Input approximation, \eqref{q214:topstable}, and
$\qnorm{p(H)-p(H_0)}\le(8/3)\qnorm{H-H_0}$ give
\[
 \qnorm{\Pi_{P_T}-P_{\rm top}}
 \le\left(\frac{2\sqrt2}{60}+\frac{16}{3}\,341110000\right)\epsilon
 \le2046660000\epsilon.
\]
Here $V_T$ uses thresholded triangle and input projections, with
common exact coordinates $F_1+C$; its aliases are the rounded rows.
Compress only the $13$ off-diagonal block:
\[
 \qnorm{2F_1V_TC-UW}\le2(3\cdot2046660000+4\cdot341110000)\epsilon
                                    \le15008840000\epsilon.
\]
For each nonempty ordinary output class choose the exact canonical
coisometry of one incident retained triangle. Define $A_0=W_{\rm ref23}$
from the canonical reference $23$ input. The reference triangle
has canonical product $I A_0=A_0$; it serves as the alias center for the
inverse controls and already supplied the initial geometric seed.
A star of aliases compares any two outputs in a class through at most
two rows. Introduce the two comparison scales
\begin{equation}\label{q214:gamma}
 \gamma_0=100L\epsilon=34111000000\epsilon,
 \qquad\gamma=2\gamma_0=G\epsilon,\quad G=68222000000.
\end{equation}
The spectral and affine comparisons give the coisometry error
$\gamma_0$, since $2\cdot15008840000+4\le34111000000$.

For each of the $19$ singleton classes without a
calibration output, put $A_c=U_cA_0$ for a left-only role
and $A_c=W_c$ for a right-only role. Its other incidences and eligibility
are checked in the manifest. For each paired class, put
$A_j=W_j$.
For every class set $B_c=A_cA_0^*$. These are contractions and
$B_{\varnothing}=I$ exactly.

Nonpaired classes retain their standard calibration or direct-definition
comparisons. In particular, whenever their respective physical roles
exist,
\[
 \qnorm{U_cA_0-A_c}\le\gamma_0,\qquad
 \qnorm{W_c-A_c}\le\gamma_0.
\]
Every paired class $j$ has an inverse triangle with canonical
product $U_jW_j$, aliased to the reference output $A_0$.
This establishes $\qnorm{U_jW_j-A_0}\le\gamma_0$ before any source
sign is used. For a literal triangle, the output class is unpaired:
its spectral-output and alias estimate gives
\[
 \qnorm{B_o-U_l(W_rA_0^*)}\le\gamma_0.
\]
Every internal source-word node has such an unpaired output. The
paired classes occur only as atomic leaves in the active literal DAG.
The source comparison below therefore proceeds directly from these
physical factors to reconstructed signs.

\subsubsection{Paired inverse controls and primitive signs}\label{q214:pairedcalibration}
The paired classes comprise the five primitive letters $D,F,G,J,u$
and the eight recovered letters $I,c,d,k,m,p,v,y$. Each has both
physical roles and exactly two output projections, its left and
right calibrations. Their multiplicities illustrate
that repeated use as a product input is compatible with this rule:
\begin{center}\small
\begin{tabular}{lrrrr}
\toprule Pair&Left uses&Right uses&Meet face&Inverse control\\\midrule
D&3&3&23&$D\cdot D$\\
F&7&6&23&$F\cdot F$\\
G&4&5&23&$G\cdot G$\\
J&4&3&12&$J\cdot J$\\
u&2&6&23&$u\cdot u$\\
\bottomrule\end{tabular}\end{center}

For each pair $j$, add a spectral triangle with inputs $U_j,W_j$ and
empty output class. Use one common reference spectral projection
with canonical product $I A_0=A_0$. Synchronize each inverse
triangle's virtual output to this reference output by one affine row.
This gives thirteen aliases. Every alias is a difference of complete virtual-output
forms, so its coordinate terms cancel.

The spectral-output and alias estimates at scale $\gamma_0$ imply
\begin{equation}\label{q214:pairedinverse}
 \qnorm{U_jW_j-A_0}\le\gamma_0.
\end{equation}
Define $A_j=W_j$ and $B_j=W_jA_0^*$. The right comparison is exact.
Multiplying \eqref{q214:pairedinverse} by $A_0^*$ and using
$A_0A_0^*=I$ yields
\begin{equation}\label{q214:pairedadjoint}
 \qnorm{B_j-U_j^*}\le\gamma_0.
\end{equation}

\paragraph{Right-controlled pairs.}
For $j=D,F,G,u$, the right primitive meet is between the graphs of
$W_j$ and $A_0$. Apply the rectangular-frame rounding argument of
Section~\ref{q214:rightprimitive} directly to the contraction $B_j$.
At scale
$\gamma_0=100L\epsilon$ it gives an exact symmetry $S_j$ with
\[
 \qnorm{B_j-S_j}\le\gamma_0/10.
\]
Equation~\eqref{q214:pairedadjoint}, taking adjoints and using $S_j=S_j^*$,
then gives
\[
 \qnorm{U_j-S_j}\le11\gamma_0/10,\qquad
 \qnorm{U_j-B_j}\le6\gamma_0/5\le\gamma.
\]
The direct definition eliminates any synchronization loss in this
primitive rounding step.

\paragraph{The left-controlled pair.}
For $j=J$, the left primitive meet first rounds the independently
approximating unitary $U_j$. The estimate in Section~\ref{q214:primitivemeet},
again with scale $\gamma_0$, gives
$\qnorm{U_j-S_j}\le\gamma_0/7$.
Consequently $\qnorm{U_j-U_j^*}\le2\gamma_0/7$, and
\eqref{q214:pairedadjoint} gives
\[
 \qnorm{U_j-B_j}\le9\gamma_0/7\le\gamma,
 \qquad\qnorm{B_j-S_j}\le8\gamma_0/7.
\]
There is no circularity: right pairs round their directly defined
corner contraction first, while the left pair rounds its physical
unitary first. In both cases the source allowances
$\qnorm{U_j-S_j}\le3\gamma$ and $\qnorm{B_j-S_j}\le4\gamma$ remain valid.
The later source reconstruction keeps these same symmetries.

These primitive estimates imply both physical corner interfaces with
the common allowance $\delta_j=4$. The eight derived pairs acquire
their signs by the direct oriented-word argument below. The signal letter $d$ is
paired; its full coisometry estimate is proved separately from its
inverse relation and reconstructed sign in Section~\ref{lem:214sourceoutput}.

At the specified witness $U_j=S_j$, $W_j=S_jA_0$ and $S_j^2=I$,
so $U_jW_j=A_0$. The reference projection and every inverse
triangle have their usual canonical top-projection witnesses; their
aliases vanish. The direct class decoders agree with the source labels.
The thirteen inverse projections and one reference projection use twelve
fewer variables than twenty-six separate calibration variables and realize
the same exact negative witness.

\subsection{Exact symmetries and direct word comparison}\label{q214:rounding}
Every source letter $j$ has the contraction $B_j=A_jA_0^*$, even if
one or both physical roles are absent. Paired primitive signs use the
canonical inputs already constructed and meets as above. All remaining source
signs are reconstructed in the order of the finite ledger.

\subsubsection{Direct symmetry signs}
All $11$ direct-sign letters use primitive meet energies.
Their selected physical inputs, on either face, are constructed before any
symmetry reconstruction. The two arguments below give the conservative
physical-left (when present) and output bounds
$3\gamma$ and $4\gamma$, so we keep $\delta_j=4$ for every direct
letter in the source recovery ledger.

\subsubsection{Primitive symmetry signs on the left face}\label{q214:primitivemeet}
For a paired primitive retain the same symmetry already constructed in
Section~\ref{q214:pairedcalibration}. The following general formulas also
construct the remaining primitive signs. For a primitive letter $j$, assume its physical $12$ input $Y_j$
still has an ordinary spectral calibration or a certified finite-forest
spectral approximation. This approximation must precede sign reconstruction and
must not use any symmetry property of $j$. Let $Y_0$ be the reference
$12$ effect. The preceding geometry gives canonical balanced graph
projections $P(U_j)$ and $P(I)$ with the safe raw-input bounds
\[
 \qnorm{Y_j-P(U_j)},\ \qnorm{Y_0-P(I)}
       \le\sqrt{2\mathcal R}+341110000\epsilon.
\]
The second estimate deliberately weakens the better reference bound.
Add one effect $T_j$ with energy
\[
 \mathcal E_j=\tr(Y_jY_0+5T_j-3T_j(Y_j+Y_0)).
\]
As in \eqref{q214:meetenergy}, the integer quadratic has trace
$8\mathcal E_j$ and contributes $2\mathcal E_j$ to $\mathcal R$.

For an arbitrary unitary $U$, direct block multiplication gives
\[
 \qnorm{[P(U),P(I)]}^2=\tfrac18\qnorm{U-U^*}^2.
\]
Applying the independent universal inequality of
Lemma~\ref{q214:commutationenergy}, and transferring both raw inputs to
their canonical graphs, gives
\[
 \begin{aligned}
 \qnorm{U_j-U_j^*}\le{}&4\sqrt2\sqrt{\mathcal E_j}
                  +16\sqrt{\mathcal R}
                  +8\cdot341110000\sqrt2\epsilon\\
 \le{}&\bigl(4\sqrt2+16+60\cdot8\cdot341110000\sqrt2\bigr)
                   \sqrt{\mathcal R}
 \le\gamma.
 \end{aligned}
\]
The last inequality follows from $\sqrt2<3/2$ and the exact
integers $L=341110000$, $G=68222000000$:
$8\sqrt2L+(4\sqrt2+16)/60<G$.
At $\mathcal R=0$ both sides vanish.

Define $S_j=\operatorname{sign}(\operatorname{Re}U_j)$, with
sign zero equal to one. For a unit-circle spectral value $e^{i\theta}$,
\[
 |e^{i\theta}-\operatorname{sign}(\cos\theta)|^2
 =2(1-|\cos\theta|)\le4\sin^2\theta
 =|e^{i\theta}-e^{-i\theta}|^2.
\]
Functional calculus proves that $S_j$ is a symmetry and
$\qnorm{U_j-S_j}\le\gamma$. The decoder interface
$\qnorm{B_j-U_j}\le\gamma$ then gives
$\qnorm{B_j-S_j}\le2\gamma$. Hence the common conservative bounds
$3\gamma$ for the physical label and $4\gamma$ for its decoder
remain valid, exactly as for a retained spectral involution triangle.

The left-face primitive meet letters are $J,w$.
All their physical left inputs have the approximations established in
Section~\ref{q214:referencefree}.

\paragraph{Virtual primitive squares.}
For a primitive $j$, the direct empty decoder is $B_{\varnothing}=I$.
When the manifest explicitly marks a literal square $jj$ as virtual,
\[
 \qnorm{B_j^2-I}\le\qnorm{B_j^2-U_j^2}+\qnorm{U_j^2-I}
 \le2\qnorm{B_j-U_j}+\qnorm{U_j-U_j^*}\le3\gamma.
\]
The first comparison telescopes between contractions; the equality
$\qnorm{U_j^2-I}=\qnorm{U_j-U_j^*}$ follows by unitary multiplication.
This supplies the stated $3\gamma$ class-square bound.

At the witness $U_j=S_j$ is already self-adjoint, so $P(S_j)$ commutes
with $P(I)$. Choose $T_j=P(S_j)P(I)$, a projection with zero energy.
Every primitive meet therefore preserves the specified witness and
remains a Hermitian integer quadratic on arbitrary contraction inputs.

\subsubsection{Primitive symmetry reconstruction on the right face}\label{q214:rightprimitive}
A primitive letter can instead use its physical $23$ input and the
reference $23$ input. Both are constructed before symmetry reconstruction.
Let $W,A_0:C\to F$ be their canonical coisometries, so
$WW^*=A_0A_0^*=I_F$. Their initial projections may differ. Define
\[
 V_W=2^{-1/2}\binom{I_F}{W^*},\qquad
 V_0=2^{-1/2}\binom{I_F}{A_0^*},\qquad
 B=WA_0^*,\qquad Z=V_W^*V_0=(I_F+B)/2.
\]
The columns are isometries and their range projections are the balanced
graphs $P_W,P_0$. Put $H=Z^*Z$. Rectangular trace cyclicity gives
\[
 \tr(P_WP_0)=\tr H,\qquad
 \tr(P_WP_0P_WP_0)=\tr H^2,\qquad
 \qnorm{[P_W,P_0]}^2=2\tr(H-H^2).
\]

The operator $B$ is a contraction.
Thus $0\le H\le I$ and
\[
 \operatorname{Re}Z-H=(I-B^*B)/4\ge0.
\]
Let $P=1_{[1/2,1]}(H)$ and define the exact symmetry $S=2P-I_F$.
Expansion of the squared norm and trace positivity give
\[
 \begin{aligned}
 \qnorm{Z-P}^2
 &=\tr H+\tr P-2\tr(P\operatorname{Re}Z)\\
 &\le\tr\bigl(H(I-P)+(I-H)P\bigr)\\
 &\le2\tr(H-H^2)=\qnorm{[P_W,P_0]}^2.
 \end{aligned}
\]
The first inequality uses $\tr(P(\operatorname{Re}Z-H))\ge0$. The second is scalar
functional calculus for $H$. Hence
\begin{equation}\label{q214:rightframeround}
 \qnorm{B-S}=2\qnorm{Z-P}\le2\qnorm{[P_W,P_0]}.
\end{equation}

Give the two raw effects their intersection energy $\mathcal E$.
The universal effect inequality and the $2$-Lipschitz commutator in
each contraction argument, together with threshold and geometry errors,
yield
\[
 \begin{aligned}
 \qnorm{B-S}
 &\le4\sqrt{\mathcal E}+8L\epsilon+8\sqrt{2\mathcal R}\\
 &\le\left(8L+\frac{4+8\sqrt2}{60}\right)\epsilon
 \le\gamma/10.
 \end{aligned}
\]
The last inequality follows from $\sqrt2<3/2$, $L\ge1$, and
$G\ge100L$. At zero energy the bound is interpreted non-strictly.
The right calibration and aliases, or a direct right decoder, give
$\qnorm{B_j-WA_0^*}\le\gamma$. Therefore
\[
 \qnorm{B_j-S}\le11\gamma/10\le4\gamma.
\]
For an unpaired primitive with a physical left role, its class
calibration gives $\qnorm{U_j-B_j}\le\gamma$, hence
$\qnorm{U_j-S}\le21\gamma/10\le3\gamma$.
For a paired primitive, retain the physical comparison already proved
in Section~\ref{q214:pairedcalibration}. The two signal letters use their stated
physical interfaces in Section~\ref{lem:214sourceoutput}.

For a certified virtual class square one also obtains
$\qnorm{B_j^2-I}\le2\qnorm{B_j-S}\le11\gamma/5\le3\gamma$.

The left and right primitive assignments are recorded below.
Each entry identifies an actual physical role, and all input approximations precede
the sign step.
\begin{center}\small
\begin{tabular}{llll}
\toprule Letter&Face&Physical input&Reference\\\midrule
B&23&\texttt{r1\_23\_Q}&\texttt{ref23\_Q}\\
D&23&\texttt{r2\_23\_Q}&\texttt{ref23\_Q}\\
F&23&\texttt{r3\_23\_Q}&\texttt{ref23\_Q}\\
G&23&\texttt{r4\_23\_Q}&\texttt{ref23\_Q}\\
J&12&\texttt{l6\_12\_Q}&\texttt{ref12\_Q}\\
n&23&\texttt{r15\_23\_Q}&\texttt{ref23\_Q}\\
o&23&\texttt{r16\_23\_Q}&\texttt{ref23\_Q}\\
t&23&\texttt{r20\_23\_Q}&\texttt{ref23\_Q}\\
u&23&\texttt{r21\_23\_Q}&\texttt{ref23\_Q}\\
w&12&\texttt{l25\_12\_Q}&\texttt{ref12\_Q}\\
z&23&\texttt{r25\_23\_Q}&\texttt{ref23\_Q}\\
\bottomrule\end{tabular}\end{center}

\subsubsection{Direct comparison of words with reconstructed signs}\label{q214:derivedword}
For every paired class put $A_j=W_j$ and $B_j=C_j=W_jA_0^*$.
The inverse control gives
\[
 \qnorm{C_j-U_j^*}\le\gamma_0\le\gamma.
\]
Its spectral and alias comparison gives
\[
 \qnorm{B_{[pq]}-U_{[p]}C_{[q]}}\le\gamma_0.
\]
For unpaired physical roles, the ordinary or singleton comparison
between that role and its class decoder still costs at most $\gamma_0$.
We will compare complete words directly with exact source symmetries.

\paragraph{Known signs and their leaf interfaces.}
For a reconstructed unpaired source letter retain
$\qnorm{B_j-S_j}\le\delta_j\gamma$.
Its available physical left or right corner input is within
$(\delta_j+1)\gamma$ of $S_j$.
For a reconstructed paired letter use the shared, stronger interface
\begin{equation}\label{q214:pairedknownsign}
 \qnorm{U_j-S_j}\le\delta_j\gamma,\qquad
 \qnorm{C_j-S_j}=\qnorm{B_j-S_j}\le\delta_j\gamma.
\end{equation}
The paired primitive construction supplies this with $\delta_j=4$.

For a nonempty active word in known signs, set
\[
 H(w)=3(|w|-1)+\sum_{j\text{ in }w}\delta_j,
\]
counting occurrences with multiplicity. Then
\begin{equation}\label{q214:directsignedword}
 \qnorm{B_{[w]}-S_w}\le H(w)\gamma.
\end{equation}
The one-letter case is the decoder interface. At an internal node
$w=pq$, both containing classes are unpaired whenever the corresponding
child is itself an internal word. Its physical input differs from
$B_{[p]}$ or $B_{[q]}$ by at most $\gamma_0$; use the inductive bound
for that child. An unpaired leaf has the same one-$\gamma$ allowance.
For a paired leaf, use \eqref{q214:pairedknownsign} directly, leaving that
allowance unused. The product and two input comparisons therefore cost
at most
\[
 (H(p)+H(q)+3)\gamma=H(w)\gamma.
\]

\paragraph{One unrounded paired target.}
Consider a retained recovery cut $u=v$ containing its target
$j$ once. Every other letter is already reconstructed. The certificate
checks that isolating $j$ in this specific cut gives an odd palindrome
$w$ in its predecessor symmetries. For each derived pair, the
target occurs inside a nontrivial product endpoint.

At the immediate parent of that occurrence, use the unknown unitary
$U_j$ when it is a left factor, and $U_j^*$ when it is a right factor.
The right substitution replaces
$C_j$ by $U_j^*$ at cost at most $\gamma_0$ from the inverse relation.
At every higher node the containing word is internal and therefore
unpaired. The preceding direct-word induction applies to this hybrid
word with zero source delta for its single unknown target. It gives
an error at most
\[
 e_j\gamma,\qquad e_j=3(|u|+|v|-2)
       +\sum_{k\ne j\text{ in the cut}}\delta_k.
\]
All factors surrounding the target are known unitaries. Multiplying by
their inverses isolates $U_j$ or $U_j^*$ near the exact palindrome
symmetry $S_j=S_w$. In the latter case take the adjoint of the whole
norm comparison. Thus
\[
 \qnorm{U_j-S_j}\le e_j\gamma,\qquad
 \qnorm{C_j-S_j}\le(e_j+1)\gamma.
\]
Choose the shared error $\delta_j=e_j+1$. Both interfaces in
\eqref{q214:pairedknownsign} now hold and can be used by later recoveries.
The actual occurrence orientations and predecessor order are included
in the exact source certificate and displayed here.
\begin{center}\small
\begin{tabular}{lllr r}
\toprule Target&Actual cut&Unknown leaf&$e_j$&$\delta_j$\\\midrule
v&$Ju=vJ$&$U_v$&18&19\\
y&$Fy=wF$&$U_y^*$&18&19\\
I&$FI=GFGF$&$U_I^*$&32&33\\
k&$Ey=yk$&$U_k^*$&62&63\\
m&$Fm=lF$&$U_m^*$&91&92\\
p&$Dp=lD$&$U_p^*$&91&92\\
c&$Jb=cJ$&$U_c$&898&899\\
d&$Fd=eF$&$U_d^*$&927&928\\
\bottomrule\end{tabular}\end{center}
 An unknown unpaired target uses its ordinary decoder comparison,
without this extra paired allowance.

\paragraph{Source and helper scope.}
After reconstruction, every source commutator meet sees physical
inputs within the existing $(\delta_j+1)\gamma$ allowance of their
exact signs. Hence the same left- and right-face coefficient
$4(\delta_a+\delta_b)+10$ applies. Ordinary source rows use
\eqref{q214:directsignedword}; all squared row bounds are computed with
the propagated deltas.

\subsubsection{Unpaired one-row recovered symmetries}
For a target $j$ with recorded odd palindrome $w=i_1\cdots i_m$,
all parent signs are available and the target occurs once in the
paid cut $a=b$. Define $S_j=w(S_{i_1},\ldots,S_{i_m})$.
It is a symmetry for every parent symmetry tuple.
For an unpaired target, keep its contraction $B_j$ in the hybrid word
and replace only known parent leaves. The same direct-word induction
bounds the cut difference by
\[
 \left(3(m-1)+\sum_{k=1}^m\delta_{i_k}\right)\gamma.
\]
The certificate checks the single-occurrence isolation: if the
relevant side is $L j R$ and the other side is $v$, free involution
reduction of $\operatorname{rev}(L)v\operatorname{rev}(R)$ is the
recorded palindrome. All surrounding factors are now exact unitaries.
Their inverses isolate $B_j$, proving
\[
 \qnorm{B_j-S_j}\le\delta_j\gamma,\qquad
 \delta_j=3(m-1)+\sum_{k=1}^m\delta_{i_k}.
\]
There are $23$ single-row palindrome recoveries,
two compound recoveries, and the special $h$ sign. Twenty-two
single-row recoveries and both compound steps precede $h$; only
$i=php$ follows it.

The missing source-letter left roles are \texttt{y0}, \texttt{z3\_\allowbreak 0}, \texttt{S3}, \texttt{T31}, \texttt{W12}, \texttt{H1\_\allowbreak 0}, \texttt{C1}, \texttt{Ggraph}, \texttt{Tau}.
The missing source-letter right roles are \texttt{x0}, \texttt{z2\_\allowbreak 0}, \texttt{p12\_\allowbreak 0}, \texttt{R31}, \texttt{M1}, \texttt{B1}, \texttt{T21}, \texttt{F2\_\allowbreak E2}, \texttt{F2\_\allowbreak w13}, \texttt{F2\_\allowbreak ladder\_\allowbreak L6}, \texttt{F2\_\allowbreak ladder\_\allowbreak h}.

\subsubsection{Reconstruction order for $h$}
Retain the eleven primitive controls $B,D,F,G,J,n,o,t,u,w,z$. Before assigning signs to $h$ or $i$, construct the following 24 symmetries in the displayed order. Each one-row recovery is paid by its selected cut; the two compound rows use the proof below. For the eight paired targets $v,y,I,k,m,p,c,d$, the direct signed-leaf comparison gives $\delta_j=e_j+1$, paying the extra decoder error after the physical-unitary recovery. Unpaired targets retain the ordinary contraction rule. Row indices refer to the 80-row presentation.
\begin{center}\small
\begin{tabular}{rllrr}
\toprule Step&Target&Exact definition&Row&$\delta_j$\\\midrule
1&E&GnG&67&18\\
2&H&JzJ&47&18\\
3&v&JuJ&45&19\\
4&x&zwz&60&18\\
5&y&FwF&61&19\\
6&A&FHF&49&32\\
7&I&FGFGF&44&33\\
8&N&vuv&76&48\\
9&L&IGI&73&76\\
10&k&yEy&63&63\\
11&l&GkG&56&77\\
12&m&FlF&58&92\\
13&p&DlD&55&92\\
14&q&FpF&54&106\\
15&s&mkm&64&253\\
16&a&LGNGL&77&608\\
17&b&IGNGI&78&884\\
18&c&JbJ&46&899\\
19&K&ucu&50&913\\
20&e&DcD&62&913\\
21&d&FeF&65&928\\
22&g&DdD&68&942\\
23&r&kgk&5&1074\\
24&f&qeq&66&1131\\
\bottomrule\end{tabular}\end{center}
At this stage all 35 coordinates other than $h,i$ are exact symmetries. In particular $r=kgk$, and $k,l,p$ are reconstructed without $h$. The selected source cuts and literal-product endpoints are retained, and every required input has been constructed before this table is used. Only afterward does Section~\ref{q214:hcontrol} define $S_h=\operatorname{sign}(\operatorname{Re}U_h)$ with $\delta_h=10915614$, and then $S_i=S_pS_hS_p$ from row 70 with $\delta_i=10915804$. These are the final two reconstructions, following the 24 displayed steps; the source theorem is invoked only after both.

\subsubsection{Two compound recoveries paid by helper residuals}\label{q214:compound}
We construct $a$ and $b$ by
compound identities and then recover $c$ and $K$.
Suppress the prefix $S$ on already reconstructed symmetries. Before the
compound steps the strict ledger has constructed
\[
 I=FGFGF,\qquad N=vuv,\qquad L=IGI,
 \qquad \delta_I=33,\quad\delta_N=48,\quad\delta_L=76.
\]
The letters $G,u,z$ are retained direct signs with errors $4\gamma$.
Put $P_\circ=GI$ and $T_\circ=LG$. The defining identity $L=IGI$ gives
$T_\circ=P_\circ^{*\,2}=P_\circ^{-2}$ exactly.

\paragraph{The two helper errors.}
For row 79 of the 80-row presentation, apply the direct signed-word comparison of
Section~\ref{q214:derivedword} to its selected cut and the already
constructed predecessor symmetries. Its graph and hybrid charge is
invariant under the selected split. On these exact
symmetries, cyclic/reverse invariance converts the residual to
$I-GzN$, giving
\[
 \qnorm{I-GzN}\le e_{79}\gamma,\qquad
 e_{79}=6+33+4+4+48=95.
\]
Multiplying by $G$ and $N$ gives $\qnorm{z-P_\circ N}\le95\gamma$.
Since $z,N$ are genuine symmetries, adjoints and conjugation by $N$ give
\[
 \qnorm{H}\le190\gamma,\qquad H=P_\circ^*N-NP_\circ.
\]
Consequently
\[
 T_\circ N-NT_\circ^*=P_\circ^*H+HP_\circ,
 \qquad\qnorm{T_\circ N-NT_\circ^*}\le380\gamma.
\]
The same argument for the selected cut of row 75 of the 80-row presentation
gives the comparison $GLG=uL$, with charge
\[
 \qnorm{GLG-uL}\le e_{75}\gamma,\qquad
 e_{75}=9+2\cdot4+2\cdot76+4=173.
\]
Right multiplication by $L$ gives $\qnorm{(GL)^2-u}\le173\gamma$.
Taking adjoints, using the already reconstructed symmetry $u=u^*$,
gives $\qnorm{T_\circ^2-u}\le173\gamma$.

\paragraph{Constructing $a$.}
Define the exact symmetry $S_a=T_\circ NT_\circ^*=LGNGL$.
In the selected cut of row 77, replace only the known parent
letters by their exact signs through the direct-word comparison, and
isolate the sole occurrence of the unpaired contraction $B_a$. The graph charge is $3\gamma$
and the parent errors are $4\gamma$ and $48\gamma$.
The certificate gives a word $R_a\in\{uN,Nu\}$ with
$\qnorm{B_a-R_a}\le55\gamma$. The helper estimates yield
\[
 \qnorm{uN-T_\circ NT_\circ^*}
 \le\qnorm{u-T_\circ^2}+\qnorm{T_\circ N-NT_\circ^*}
 \le553\gamma.
\]
Since $S_a$ is symmetric, $uN$ and its adjoint $Nu$ have the same
distance to it. Thus $\qnorm{B_a-S_a}\le608\gamma$, without using $B_a^2$ or
$B_a^*=B_a$.

\paragraph{Constructing $b$.}
Define the exact symmetry $S_b=P_\circ^*NP_\circ=IGNGI$.
The selected cut of row 78 gives
$R_b\in\{S_aLG,GLS_a\}$ by the same single-occurrence isolation.
Its hybrid comparison gives
\[
 \qnorm{B_b-R_b}\le(6+4+608+76)\gamma=694\gamma.
\]
Now $S_aT_\circ=T_\circ N$ and
$T_\circ N-P_\circ^*NP_\circ=P_\circ^*H$, so the correction is
at most $190\gamma$. The adjoint word has the same distance to the
symmetry $S_b$. Hence $\qnorm{B_b-S_b}\le884\gamma$ without
using $B_b^2$ or a primitive control for $b$.

\paragraph{Completing the stage without $h,i$.}
Resume the one-row rules in the certified predecessor order. For each
paired target, the oriented physical-leaf recovery adds one unit to
the hybrid allowance; unpaired targets retain their ordinary contraction
comparison. The resulting bounds are:
\[
 \begin{array}{c|rrrrrrrrr}
 j&a&b&c&K&e&d&g&r&f\\\hline
 \delta_j&608&884&899&913&913&928&942&1074&1131.
 \end{array}
\]
Here $c=JbJ$, $K=ucu$, $e=DcD$, $d=FeF$, $g=DdD$, $r=kgk$,
and $f=qeq$. In particular, the paired targets $c,d$ have hybrid
allowances $898,927$ and shared decoder bounds $899,928$, respectively.
The eleven primitive controls, 22 one-row recoveries (eight paired) and
these two compound recoveries give the 35 coordinates preceding $h,i$. At zero source defect every step
equals the specified source coordinate, preserving the prescribed witness.

The signal coordinates are already available, with physical errors
$3\gamma$ for $t$ and $(\delta_d+1)\gamma=929\gamma$ for $d$.

\subsubsection{Paid cuts and reconstruction order}
The checker isolates each target from its single occurrence in the selected cut, following the displayed $h$-free order and using only previously reconstructed symmetries.

For the compound rows the selected cuts remain $au=N$ and $ab=LG$.
They isolate $B_a$ near $Nu$ with coefficient 55 and $B_b$ near
$S_aLG$ with coefficient 694. The symmetric comparison $S_a$ permits
the adjoint orientation of $uN$ without asserting $B_a=B_a^*$.
The paid helper cuts remain $Lu=GLG$ and $IG=Nz$; their letters are
already reconstructed and their bounds are now 173 and 95.
The helper paths below, with the stated charges, are
checked separately on words omitting $h,i,C$.

\subsubsection{Recovering the symmetry of $h$ from paid source meets}\label{q214:hcontrol}
Both physical $h$ inputs and both ordinary $h$ calibrations enter the
214-variable construction. Construct all canonical inputs in the order
specified by the calibration forest before assigning a sign to $h$.
Let $H_0=U_h$ be the resulting canonical left unitary.
The retained calibration gives
$\qnorm{B_h-H_0}\le\gamma$.
The canonical $h$ graph is also available for its ordinary geometric
neighbors, independently of the source sign constructed below.

\paragraph{Four source meets before assigning the sign.}
For arbitrary unitaries $U,V$, block multiplication gives
\[
 \qnorm{U^*V-V^*U}^2=8\qnorm{[P(U),P(V)]}^2,
 \qquad P(U)=\tfrac12\begin{pmatrix}I&U\\U^*&I\end{pmatrix}.
\]
The universal effects inequality and the raw-input errors imply
\[
 \qnorm{U^*V-V^*U}
 \le4\sqrt2\sqrt{\mathcal E}+16\sqrt{\mathcal R}
       +8\sqrt2L\epsilon
 \le\gamma/7\le2\gamma/7.
\]
For the last comparison use $G=200 L$, $L=341110000$,
$\epsilon=60\sqrt{\mathcal R}$ and $\sqrt2<3/2$.
All terms vanish at zero error, and all norms use the ambient trace.

The four retained source meets involving $(h,a),(h,z),(h,y),(h,e)$
are all on the left face. Their row numbers in the 80-row presentation are respectively
8,17,39,9, and their active indices are 5,13,29,6.
Each other letter $s$ is already an exact symmetry, with
$\qnorm{U_s-S_s}\le(\delta_s+1)\gamma$.
Replacing the two occurrences of $U_s$ in the relative-unitary estimate
and multiplying by the symmetry $S_s$ gives
\begin{equation}\label{q214:presignconjugation}
 \qnorm{S_sH_0S_s-H_0^*}\le\eta_s\gamma,
 \qquad\eta_s=2\delta_s+3.
\end{equation}
Thus $(\eta_a,\eta_z,\eta_y,\eta_e)=(1219,11,41,1829)$.

\paragraph{Paid helper identities with no $h$ or $i$.}
The helper certificate contains 28 acyclic helper nodes and 303 literal
rewrite steps. Apply the charged-cap argument of Section~\ref{app:214coefficient} after
the exact substitutions $j=lpl$ and $M=LGL$; their defining
rows 57 and 74 vanish identically. Every active proof word and
every charged leaf omits $h,i,C$.

The exact positive weighted sums of those source-row coefficients are
\[
 \begin{array}{c|r}
 \text{residual}&\text{bound in units of }\gamma\\\hline
 {[c,l]}&872404\\
 {[c,p]}&54496\\
 Bu-cB&23611\\
 {[u,z]}&2816
 \end{array}
\]
Here and in the following identities letters denote these already
constructed exact symmetries. The paths and positive charge vectors are in
\path{quadratic_214/H_CONTROL_CERTIFICATE.json}; \path{quadratic_214/DERIVED_PAIR_SOURCE.json}
evaluates them with the source-row bounds used here. The independent
source audit checks the four displayed sums.

\paragraph{The first quantitative identity.}
Put $P=GI$ and $j=lpl$. The definitions give $LG=P^{-2}$ and
$a=P^{-2}NP^2$ exactly. The paid helper rows 75 and 79, with the
already proved bounds $e_{75}=173$ and $e_{79}=95$, give
\[
 \qnorm{u-P^{-4}}\le e_{75}\gamma,\qquad
 \qnorm{z-PN}\le e_{79}\gamma,\qquad
 \qnorm{NP-P^*N}\le2e_{79}\gamma.
\]
Since $u$ is an exact symmetry,
$\qnorm{P^8-I}\le2e_{75}\gamma$.
Move $N$ across the three factors of $P^3$ to obtain
\[
 \qnorm{az-P^{-5}}
 \le\qnorm{z-PN}+\qnorm{NP^3-P^{-3}N}
 \le7e_{79}\gamma.
\]
Taking fourth powers costs $28e_{79}\gamma$.
Comparing $P^{-20}$ with $P^{-4}$ costs $4e_{75}\gamma$.
Consequently
\begin{equation}\label{q214:hfirstidentity}
 \qnorm{u-(az)^4}\le(5e_{75}+28e_{79})\gamma=3525\gamma.
\end{equation}

\paragraph{The second quantitative identity.}
The row-59 residual gives both
$\qnorm{D-jy}\le278\gamma$ and $\qnorm{D-yj}\le278\gamma$;
the second follows by taking adjoints of exact symmetry words.
Row 62 has the conservative residual bound $1826\gamma$.
Telescoping $j=lpl$ yields
\[
 \qnorm{jcj-c}\le2\qnorm{[c,l]}+\qnorm{[c,p]}.
\]
Substitute $D\approx yj$ and $D\approx jy$ into the two positions
of $DcD$, then compare with $e$ using its paid residual. Conjugating
by $y$ gives
\begin{equation}\label{q214:hsecondidentity}
 \qnorm{c-yey}\le E_c\gamma,
 \qquad E_c=1826+2(278)+2(872404)+54496=1801686.
\end{equation}
The estimate remains valid if a defining palindrome makes one of
these residuals identically zero; the stated allowances are conservative.

\paragraph{The third quantitative identity.}
The paid rows 51 and 50 imply
$B\approx Kz$ and $K\approx ucu$, with errors
$924\gamma$ and $1826\gamma$. Thus
$\qnorm{B-ucuz}\le(924+1826)\gamma$.
Insert this comparison twice into the helper residual $Bu-cB$.
Commuting the adjacent $z,u$ in $ucuzu$ costs at most
$2816\gamma$ and reduces that word to $ucz$.
Right multiplication by $z$ then gives
\[
 \qnorm{uc-(uc)^{-2}}\le E_b\gamma,
 \qquad E_b=23611+2(924)+2(1826)+2816=31927.
\]
Multiplication by the exact unitary $(uc)^2$ proves
\begin{equation}\label{q214:hthirdidentity}
 \qnorm{(uc)^3-I}\le31927\gamma.
\end{equation}

\paragraph{An odd conjugation cycle forces a sign.}
By \eqref{q214:presignconjugation}, conjugation by each of $a,z$ sends
$H_0$ close to $H_0^*$. Their product $az$ sends $H_0$ close to
itself, with error $(\eta_a+\eta_z)\gamma$.
Four conjugations and \eqref{q214:hfirstidentity} therefore give
\[
 \qnorm{uH_0u-H_0}\le\eta_u\gamma,
 \qquad\eta_u=4(1219+11)+2(3525)=11970.
\]
Conjugation by the odd word $yey$ sends $H_0$ close to $H_0^*$.
Using \eqref{q214:hsecondidentity},
\[
 \qnorm{cH_0c-H_0^*}\le\eta_c\gamma,
 \qquad\eta_c=2(41)+1829+2(1801686)=3605283.
\]
Thus conjugation by $uc$ interchanges $H_0,H_0^*$ within
$(\eta_u+\eta_c)\gamma$; the assertion for the adjoint follows
by taking adjoints. Three iterations still interchange them.
Compare $(uc)^3$ with $I$ using \eqref{q214:hthirdidentity}, at cost
$2(31927)\gamma$ for its two appearances. It follows that
\[
 \qnorm{H_0-H_0^*}
 \le\bigl(3(11970+3605283)+2(31927)\bigr)\gamma
 =10915613\gamma.
\]
Define $S_h=\operatorname{sign}(\operatorname{Re}H_0)$, with
sign zero equal to one. The unitary functional-calculus inequality
proved for the primitive signs gives
\[
 \qnorm{H_0-S_h}\le\qnorm{H_0-H_0^*},\qquad
 \qnorm{B_h-S_h}\le10915614\gamma.
\]
Only now use row 70, whose cut is $pi=hp$, to define
$S_i=S_pS_hS_p$. Its ordinary comparison gives
\[
 \delta_h=10915614,\qquad
 \delta_i=6+2(92)+10915614=10915804.
\]
Neither $h^2=I$ nor $i^2=I$ was used before these steps.

\paragraph{Completion, source cost, and witness.}
All 37 source coordinates are now exact symmetries. Compute every
retained ordinary-row coefficient and all 27 source-meet coefficients
with these error bounds. With $h(w)$ as in \eqref{q214:word}, their squared sum is
\[
 c=\sum_{r\in\mathcal L}(h(a_r)+h(b_r))^2
       +\sum_{r\in\mathcal C}q_r^2
   =19184111188572071.
\]
The subsequent source application uses this value and the
coefficient $\lambda$ from Theorem~\ref{thm:quadratic214}.
The four pre-sign energies remain paid source rows after $h$ is
rounded. Using the same nonnegative energy in several comparison
estimates is legitimate: each estimate is bounded by the common
budget $\mathcal R$, and the full displayed $c$ pays the combined
source estimate.

Every construction is global in the defect and uses no positive-rank
division. At the specified witness the paid identities are exact, the
reconstruction palindromes recover the prescribed coordinates, and $H_0=h$
is a symmetry fixed by sign rounding. Thus the construction realizes
the specified negative witness.

All long words and signs here are comparison operators in the proof.
The polynomial consists of the specified source meets, paired
calibration terms, and scalar penalty; its degree is two.

\subsubsection{Application of the typed source theorem}
All $37$ retained coordinates $S_j$ are now symmetries
in the same finite corner. The direct-word estimate gives
\begin{equation}\label{q214:word}
 \qnorm{S_w-B_{[w]}}\le h(w)\gamma,\qquad
 h(w)=3(|w|-1)+\sum_{j\text{ in }w}\delta_j
 \quad(w\ne\varnothing),
\end{equation}
counting occurrences with multiplicity; set $h(\varnothing)=0$.
Every retained literal source cut is therefore bounded by
$(h(a_r)+h(b_r))\gamma$, using \eqref{qsrc:cuttrace} only after symmetry
reconstruction. Source meets use \eqref{q214:meetrowbound}; the squared
sum is \eqref{q214:mixedcapbound}. The source theorem is applied only to
this exact symmetry tuple.
For the signal letters we retain conservative physical bounds
\[
 \qnorm{U_t-S_t}\le3\gamma,\qquad
 \qnorm{U_d-S_d}\le929\gamma.
\]
The direct sign gives $b_t=3$; the paired recovery gives
$\delta_d=928$, and we choose $b_d=\delta_d+1=929$.
The needed full coisometry bounds are established separately below.

\subsection{Terminal commutation through an intersection effect}\label{q214:terminal}
For selected source rows $ab=ba$, one independent intersection effect
pays the commutator without requiring literal decoder endpoints for
$ab$ and $ba$. The graph is compiled directly for the remaining word constraints. The following
estimate is global on arbitrary effects.

\begin{lemma}\label{q214:commutationenergy}
For effects $0\le A,B,T\le I$ in a finite tracial algebra, put
\[
 \mathcal E(A,B,T)=\tr\bigl(AB+5T-3T(A+B)\bigr).
\]
Then $\mathcal E\ge0$ and
\[
 \qnorm{[A,B]}^2\le4\mathcal E(A,B,T).
\]
The trace need not be normalized. The factor four is optimal.
\end{lemma}
\begin{proof}
The function $\mathcal E-\qnorm{[A,B]}^2/4$ is separately concave in
$A$ and $B$: the first term is affine in either variable and the
squared Hilbert norm of the commutator is convex. Expressing each
effect as its spectral layer-cake integral of projections and applying
Jensen twice reduces the assertion to projections $P,Q$. One may first
pass to the finite tracial von Neumann representation, so the spectral
projections belong to the working algebra.

For fixed $P,Q$, set $H=5I-3(P+Q)$. Positivity of
$\tr(TH_+)$ and $\tr((I-T)H_-)$ gives
$\tr(TH)\ge-\tr(H_-)$. The negative spectral projection of $H$
attains equality. Thus
\[
 \min_{0\le T\le I}\mathcal E(P,Q,T)
   =\tr(PQ)-\tr\bigl((3(P+Q)-5I)_+\bigr).
\]
On a noncommuting two-dimensional fiber of the two-projection
decomposition, write
\[
 P=\begin{pmatrix}1&0\\0&0\end{pmatrix},\qquad
 Q=\begin{pmatrix}t^2&t\sqrt{1-t^2}\\
 t\sqrt{1-t^2}&1-t^2\end{pmatrix},\qquad 0\le t\le1.
\]
The eigenvalues of $P+Q$ are $1+t$ and $1-t$. With the
unnormalized two-dimensional trace, $\tr(PQ)=t^2$ and
$\qnorm{[P,Q]}^2=2t^2(1-t^2)$. The minimum energy is
\[
 e_{\min}(t)=\begin{cases}
 t^2,&0\le t\le2/3,\\
 (1-t)(2-t),&2/3\le t\le1.
 \end{cases}
\]
Both branches are nonnegative, and their commutator slacks factor as
\[
 4e_{\min}(t)-2t^2(1-t^2)=\begin{cases}
 2t^2(1+t^2),&0\le t\le2/3,\\
 2(1-t)^2(t^2+2t+4),&2/3\le t\le1.
 \end{cases}
\]
On every commuting scalar fiber the minimum and commutator are zero.
Off-diagonal entries of $T$ do not contribute to $\tr(TH)$.
Integrating the fibers and applying the earlier concavity proves the
claim. The ratio $2t^2(1-t^2)/((1-t)(2-t))$ tends to four as
$t\uparrow1$, proving optimality of the universal constant.
\end{proof}

\subsubsection{Exact integer energy and eligibility}
Under $A=(I+X_a)/2$, $B=(I+X_b)/2$, $T=(I+X_T)/2$, define
\begin{equation}\label{q214:meetenergy}
 \begin{aligned}
 J_{a,b,T}={}&10I-4X_a-4X_b+8X_T+\qac{X_a}{X_b}\\
             &-3\qac{X_T}{X_a}-3\qac{X_T}{X_b}.
 \end{aligned}
\end{equation}
Its trace is exactly $8\mathcal E(A,B,T)$. Thus adding this
integer polynomial to a budget of trace $4\mathcal R$ adds
$2\mathcal E$, not $\mathcal E$, to $\mathcal R$.

Each selected row is a certified atomic commutation relation, up to
the norm-preserving cyclic and inverse cap operations. Both letters
retain physical inputs on a common selected face, either $12$ or $23$.

\begin{center}\small
\begin{tabular}{rcrr@{\qquad}rcrr}
\toprule Row&Pair&Face&$q_{ab}$&Row&Pair&Face&$q_{ab}$\\\midrule
1&$d,e$&12&7374&22&$i,o$&23&43663242\\
4&$h,m$&12&43662834&23&$A,p$&23&506\\
6&$e,k$&12&3914&24&$A,o$&23&154\\
8&$a,h$&12&43664898&25&$x,y$&23&158\\
9&$e,h$&12&43666118&28&$d,z$&12&3738\\
11&$e,t$&12&3678&29&$k,z$&12&278\\
12&$g,t$&12&3794&30&$A,k$&23&390\\
13&$f,t$&23&4550&31&$a,k$&12&2694\\
14&$a,w$&12&2458&33&$A,B$&23&154\\
15&$v,x$&23&158&34&$B,k$&23&278\\
16&$w,y$&12&102&36&$B,t$&23&42\\
17&$h,z$&12&43662482&37&$B,p$&23&394\\
18&$A,h$&23&43662594&39&$h,y$&12&43662542\\
21&$h,q$&12&43662890&&&&\\
\bottomrule\end{tabular}\end{center}
The table lists every retained source row paid by an intersection energy. Row indices refer to the 80-row presentation. The manifest gives each physical input pair and its independent intersection effect.

For each selected row, introduce one independent effect $T_{ab}$ with
energy \eqref{q214:meetenergy}.

\subsubsection{Paying a left-face commutation residual after symmetry recovery}
Write $A_j$ for the actual input effect in the energy, $Q_j$ for its
spectral threshold projection, and $P(U_j)$ for the balanced $12$
graph projection of the canonical unitary $U_j$. The retained geometry
and recovery estimates give
\[
 \qnorm{A_j-Q_j}\le\sqrt{2\mathcal R},\qquad
 \qnorm{Q_j-P(U_j)}\le341110000\epsilon,\quad
 \qnorm{U_j-S_j}\le(\delta_j+1)\gamma.
\]
For a paired letter the direct signed-leaf interface is even stronger.
Every norm uses the ambient trace. Exact graph multiplication gives
\[
 \qnorm{P(U_j)-P(S_j)}=\qnorm{U_j-S_j}/\sqrt2,\qquad
 \qnorm{[P(S_a),P(S_b)]}^2=\tfrac18\qnorm{[S_a,S_b]}^2.
\]
The right-hand norm is corner-supported and is not separately
renormalized. For contractions, a commutator changes by at most twice
the sum of its input distances. Lemma~\ref{q214:commutationenergy} therefore
implies
\[
 \begin{aligned}
 \qnorm{[S_a,S_b]}\le{}&4\sqrt2\sqrt{\mathcal E(A_a,A_b,T_{ab})}
                    +16\sqrt{\mathcal R}\\
 &+\bigl(4(\delta_a+\delta_b)+8\bigr)\gamma
   +8\cdot341110000\sqrt2\epsilon.
 \end{aligned}
\]
Since $8\cdot341110000\sqrt2<G=68222000000$, $\gamma=G\epsilon$, and
$\epsilon=60\sqrt{\mathcal R}$, the last term is at most $\gamma$.
The $16\sqrt{\mathcal R}$ term retains the two raw-effect threshold
errors; the energy is never silently evaluated at rounded inputs.
The budget contains $2\mathcal E$, hence $\mathcal E\le\mathcal R$;
also $(4\sqrt2+16)\sqrt{\mathcal R}\le\gamma$. Consequently
\begin{equation}\label{q214:meetrowbound}
 \qnorm{[S_a,S_b]}\le q_{ab}\gamma,\qquad
 q_{ab}=4(\delta_a+\delta_b)+10.
\end{equation}

\subsubsection{Paying a source commutator on the right face}\label{q214:rightsource}
For a source row $ij=ji$, both inputs may instead be on face $23$.
They must have independently established ordinary or finite-forest
approximations. A mixed $12/23$ pair is not covered by this argument.
Let $W_i,W_j,A_0:C\to F$ be the constructed coisometries and put
\[
 P_0=A_0^*A_0,\qquad C_i=W_iA_0^*,\qquad
 D_i=W_i(I-P_0),
\]
with analogous definitions for $j$. Suppose the already reconstructed
source symmetries satisfy $\qnorm{C_i-S_i}\le e_i$ and
$\qnorm{C_j-S_j}\le e_j$. The coisometry identities imply
\[
 D_iD_i^*=I_F-C_iC_i^*,\qquad
 \qnorm{D_iD_i^*}\le2e_i.
\]
Trace cyclicity and Hilbert--Schmidt Cauchy--Schwarz give the linear
cross-product bound
\[
 \begin{aligned}
 \qnorm{D_iD_j^*}^2
 &=\tr\bigl((D_i^*D_i)(D_j^*D_j)\bigr)\\
 &\le\qnorm{D_i^*D_i}\qnorm{D_j^*D_j}
  =\qnorm{D_iD_i^*}\qnorm{D_jD_j^*}\le4e_ie_j.
 \end{aligned}
\]
The equality of the positive-product norms holds for rectangular maps
by the same trace cyclicity. Since
$W_iW_j^*=C_iC_j^*+D_iD_j^*$,
\[
 \qnorm{W_iW_j^*-S_iS_j}
 \le e_i+e_j+2\sqrt{e_ie_j}\le2(e_i+e_j).
\]

For the balanced graph projections $P_i,P_j$, direct block
multiplication gives
\[
 F[P_i,P_j]F=(W_iW_j^*-W_jW_i^*)/4.
\]
Compression and the preceding bound in both orders yield
\[
 \qnorm{[S_i,S_j]}\le4\qnorm{[P_i,P_j]}+4(e_i+e_j).
\]
The class calibration or direct right decoder gives
$e_i\le(\delta_i+1)\gamma$ and likewise for $j$.
Moving from the canonical graphs to the two raw input effects and
using their universal intersection energy gives
\[
 \begin{aligned}
 \qnorm{[S_i,S_j]}\le{}&8\sqrt{\mathcal E}
    +16L\epsilon+16\sqrt{2\mathcal R}\\
   &+[4(\delta_i+\delta_j)+8]\gamma\\
 \le{}&[4(\delta_i+\delta_j)+10]\gamma.
 \end{aligned}
\]
Indeed the extra scalar allowance is at most
$[16L+(8+16\sqrt2)/60]\epsilon<\gamma$ for $G\ge100L$, $L\ge1$;
the retained two-$\gamma$ allowance is conservative.
Thus the same charge $q_{ij}$ works on either face. The source
comparison is applied only after all source coordinates are exact
symmetries, so cyclic/reverse orientations preserve its norm.

At the witness $W_i=S_iA_0$ and $W_j=S_jA_0$ with the prescribed
initial space. If $S_i,S_j$ commute, their graph projections commute,
and their intersection projection gives zero energy. The same integer
quadratic is used, with its input names changed to the chosen face.
Every chosen face and both surviving physical inputs are listed in
the manifest and in the preceding source-meet table.

Let $\mathcal C$ index these commutation rows and let $\mathcal L$
index the remaining literal source rows. For $r\in\mathcal C$, let
$(\alpha_r,\beta_r)$ be its commuting letter pair, irrespective of
the cyclic orientation of its stored cap cut. With
$h(w)=\sum_{j\text{ in }w}\delta_j+3(|w|-1)$ as before,
the complete source comparison is now
\begin{equation}\label{q214:mixedcapbound}
 \tr Q_{66}(S)\le
 \left(\sum_{r\in\mathcal L}(h(a_r)+h(b_r))^2
       +\sum_{r\in\mathcal C}q_{\alpha_r\beta_r}^2\right)\gamma^2
 =19184111188572071\gamma^2.
\end{equation}
Every retained source row appears exactly once. The source theorem and its coefficient $\lambda$ apply to the exact
symmetry tuple reconstructed above.

\subsubsection{Witness and degree}
At the specified source witness $S_aS_b=S_bS_a$, so the two balanced
graph projections $Q_a,Q_b$ commute. Set $T_{ab}=Q_aQ_b$.
This is a projection, $T_{ab}(Q_a+Q_b)=2T_{ab}$, and
$\mathcal E(Q_a,Q_b,T_{ab})=0$. Thus the witness variable
$X_{T_{ab}}=2T_{ab}-I$ is a symmetry with zero penalty contribution.
The signal and the specified coordinate assignments realize the witness.
The energy \eqref{q214:meetenergy} is Hermitian, integral, and quadratic;
it contains no $E_3$ coordinate.

\subsection{Coefficient data}
In \path{quadratic214/eliminated_expanded.json}, a word is an ordered list of variable labels, with the empty list denoting the identity. A coefficient triple $[a,b,c]$ denotes $a+b\kappa+c\kappa^2$. Thus the constant coefficient is
\[
 3+521350\kappa+1999648800\kappa^2.
\]
The labels denote the contraction coordinates $X_v$, although they retain the names of the underlying projection gadgets. The incidence manifest and the collected coefficient array are two independently compared descriptions of the same polynomial.

\section{The characteristic-five source}\label{app:fiveproof}
A finite matrix corner has identity $q$ and positive tracial functional
$t$, with $p=t(q)$. Write $\|x\|_2^2=t(x^*x)$. If $p>0$, division by $p$
gives its normalized trace; if $p=0$, the trace inequalities are trivial.
A symmetry is a self-adjoint unitary.

The exact coefficient parameters are
\begin{equation}\label{fivemerge:eq:parameters}
 J_5=2^{\,2^{1122100}},\qquad \delta_5=2^{-J_5},\qquad
 \mu_5=2^{2J_5+30},\qquad C_{33}=2^{56},\qquad C=2^{57}.
\end{equation}
All exponentiations are integer exponentiations. In this appendix, $\mathcal G$ denotes the canonical phase-corrected source $\mathcal G^{\rm can}$ of Section~\ref{sec:degree24}:

\begin{equation}\label{fivemerge:eq:source-polynomial}
 \mathcal G=P_0+\mu_5\lambda(T)P,\qquad
 \lambda(T)=2C+64T^6,\qquad P=R_{\rm corr}+(30+84T^2)E.
\end{equation}
Its 101 real components have declared self-adjoint or skew-adjoint types.
Section~\ref{fivemerge:sec:components} defines $P_0$, the corrected caps $R_{\rm corr}$ and the recovery budget $E$; the source files specify all words and coefficients.

\begin{theorem}\label{fivemerge:thm:conditional}
Using the verified finite transfer and incidence
certificates specified in Sections~\ref{fivemerge:sec:bridge}--\ref{fivemerge:sec:components},
for every finite matrix corner every tuple of the 101 declared
real components of norm at most the specified $T$ satisfies
\[
 t(\mathcal G)\ge\tfrac12t(q).
\]
The coefficient recipe for $T$ has $T\ge4$. The conclusion holds with the specified coefficients of the polynomial.

The exact witness construction
in Section~\ref{fivemerge:sec:witness} gives a typed tuple in $M_2(A)$, for a finite
tracial algebra $A$, with all component norms at most one and normalized
trace $-1$. Its two cheap skew components are
$\mathsf J\otimes S_1$ and $\mathsf J\otimes S_2$, where
$\mathsf J=\left(\begin{smallmatrix}0&1\\-1&0\end{smallmatrix}\right)$ and
$S_1,S_2$ are commuting symmetries.
\end{theorem}

The proof has the following chain of scalar trace inequalities:
\begin{equation}\label{fivemerge:eq:chain}
\begin{gathered}
 \tr(F_{57}^0)+\mu_5Q_{57}\ge\tfrac34,\\
 Q_{57}\le2^{56}Q_{33},\quad
 t(F_{57}^0)=t(F_{33}^0),\quad Q_{33}\le2Q_{65},\\
 t(F_{33}^0)+2^{57}\mu_5Q_{65}\ge\tfrac34p,\\
 t(P)\ge\tfrac12Q_{65}(U)+64I,\qquad
 t(\mathcal G)\ge\left(\tfrac34-\frac1{4\mu_5}\right)p\ge\tfrac12p.
\end{gathered}
\end{equation}
Here the first line concerns normalized traces on correctly typed
unitary tuples; the subsequent lines use the corner trace and the
rounded tuple $U$. The third line follows from the first two.

\Needspace{6\baselineskip}
\subsection{The native 57-coordinate source}
\subsubsection{Presentation and defect conventions}
The defining file \path{degree24/f5_source57.json} fixes 57 unitary
coordinates and 170 paid residuals $r=L-R$. For a normalized finite
matrix trace set
\[
 Q_{57}=\sum_r\|r\|_2^2,\qquad
 \epsilon=\max_r\|r\|_2,\qquad F_{57}^0=-1+P^*Q+Q^*P.
\]
There are 18 root coordinates, seven action coordinates, 15
root-crossing and product-crossing helpers, nine integer-relation
helpers, six square helpers, and two signal coordinates. To explain
the root notation, put $t_0=1$ and use
\[
\begin{aligned}
 x&=e_{12}(1),& y&=e_{23}(1),&z_l&=e_{31}(t_l),\\
 c&=e_{13}(1),&d_l&=e_{32}(-t_l),&e_l&=e_{21}(t_l)
 \quad(0\le l\le3),
\end{aligned}
\]
with $p_{ab}=e_{31}(t_at_b)$ for $1\le a<b\le3$, over
$\fiveF_5[t_1,t_2,t_3]$. The sign in $d_l$ gives
$[x,z_l]=d_l$ and $[d_a,e_b]=p_{ab}^{-1}$ for $a\ne b$.
This exact realization explains the names; it is not the domain of the
finite-matrix theorem, which allows arbitrary correctly typed tuples.

The paid rows are grouped as follows:
\begin{center}
\begin{tabular}{@{}lr@{}}
\toprule
Row family&Count\\\midrule
Root commutations&60\\
Nine split root crossings&18\\
Six split product crossings&12\\
Integer commutations&6\\
Six split integer Steinberg relations&12\\
Three split opposite-root relations&6\\
Common-compression rows&42\\
Six split order-five relations&12\\
Signal rows&2\\\midrule
Total&170\\\bottomrule
\end{tabular}
\end{center}
Every side has length at most two. The split of $AB=CBA$ is
$H=CB$, $HA=AB$. The split of $ab^*a=b^*ab^*$ is
$H=ab^*$, $Ha=b^*H$. For each basic root $U$, order five is imposed
by $H=U^2$, $H^2=U^*$. The signal rows are $uP=hu$, $tQ=ht$,
where $u=u_{12}$ and $h=z_2$. Thus their unsplit errors are at most
$2\epsilon$, $3\epsilon$, and $3\epsilon$, respectively.

\subsubsection{A rational root spectral identity}
Let $U_1=e_{12}(\fiveF_5)$, $U_2=e_{23}(\fiveF_5)$, and
$U_3=e_{31}(V)$, where
$V=\fiveF_5+\fiveF_5t_1+\fiveF_5t_2+\fiveF_5t_3$. Write $P_i$ for their group
averages in a characteristic-zero group algebra. Each pair has a
presentation with one generator $X$, $m$ generators $Z_j$, and $m$
central generators $C_j$, all of order five, with
$XZ_j=C_jZ_jX$ and the other prescribed commutations. Here $m=1,4,4$;
normal forms $X^a\prod Z_j^{b_j}\prod C_j^{c_j}$ have exponents modulo
five and yield group order $5^{1+2m}$. Their distinctness follows from
the separately visible root and central parameters in the elementary
matrix realization.

For $A_{ij}=2-P_i-P_j$, the trivial central characters give commuting
projections. A nonzero central character reduces the nontrivial
commutator pairing to a five-dimensional Weyl block. Either the
kernel character makes one average vanish, or the two invariant lines
have overlap of absolute value $1/\sqrt5$. Consequently
\begin{equation}\label{fivemerge:eq:root-spectrum}
 \fivespec(A_{ij})\subset\{0,1,2,1-1/\sqrt5,1+1/\sqrt5\}.
\end{equation}
Define
\[
\begin{aligned}
 E_1(s)&=-s(s-2)(5(s-1)^2-1),\\
 E_2(s)&=s(s-1)(5(s-1)^2-1)/8,\\
 E_w(s)&=-\tfrac{25}{4}s(s-2)(s-1)^2,\\
 V_w(s)&=(s-1+29/65)E_w(s),\\
 m(s)&=s(s-1)(s-2)((s-1)^2-1/5).
\end{aligned}
\]
Direct rational polynomial division gives
\begin{equation}\label{fivemerge:eq:pair-sos}
 s^2-\tfrac{11}{20}s\equiv
 \tfrac9{20}E_1(s)^2+\tfrac{29}{10}E_2(s)^2+
 \tfrac{13}{8}V_w(s)^2+\tfrac1{650}E_w(s)^2\pmod{m(s)}.
\end{equation}
Put $Q_i=1-P_i$ and $D_H=(Q_1+Q_2+Q_3)/6$. Since $Q_i^2=Q_i$,
\[
 \sum_{i<j}A_{ij}^2=36D_H^2+6D_H,
 \qquad \sum_{i<j}A_{ij}=12D_H.
\]
Equations~\eqref{fivemerge:eq:root-spectrum}--\eqref{fivemerge:eq:pair-sos} therefore give
\begin{equation}\label{fivemerge:eq:root-sos}
\begin{split}
 D_H^2-D_H/60=\sum_{i<j}\bigl(&E_1(A_{ij})^2/80
 +29E_2(A_{ij})^2/360\\
 &+13V_w(A_{ij})^2/288+E_w(A_{ij})^2/23400\bigr).
\end{split}
\end{equation}
All weights are positive rational numbers. The root parameters have
characteristic five, but the coefficients, positivity, spectral
operators, and Hilbert-space estimates in this argument have
characteristic zero. No positive-order assertion is being made inside
$\fiveF_5$.

\subsubsection{Approximate transfer}\label{fivemerge:sec:analytic}
For a unitary list $w_1,\ldots,w_h$, use the self-adjoint averaging
Laplacian on normalized Hilbert--Schmidt space
\[
 L(y)=\frac1{2h}\sum_{j=1}^h
 \left(y-\frac{w_jyw_j^*+w_j^*yw_j}{2}\right).
\]
The root list uses exponent representatives $-2,-1,0,1,2$, lists each
$U_3$ element once, and repeats each $U_1,U_2$ element 125 times.
It has $h_H=1875$ and exact Laplacian $D_H$. Its selected words have
length at most eight. The integer list consists of its six elementary
shears, ordered $12,13,21,23,31,32$.

For approximate native tuples, the elementary estimates behind the
root transfer are
\[
 \|X^5-1\|_2\le3\epsilon,\qquad
 \|X^5Z-C^5ZX^5\|_2\le20\epsilon,\qquad
 \|C^5-1\|_2\le26\epsilon.
\]
Each finite pair-presentation test is bounded by $32\epsilon$. There
are $2m^2+3m+1$ tests per pair, hence 96 in total. The exactly checked
free pair-SOS residual has coefficient $\ell^1$ norm at most
$2390851927/5<2^{29}$ and word length at most 80. The collection
bound is $1024(L+1)^4<2^{36}$ for $L\le80$.
Here is a concrete collection bound. Replace each initial inverse by
four positive letters, using at most $L$ power-test replacements, so
there are $n\le4L$ letters. At most $n^2$ adjacent root swaps sort the
noncentral letters; every crossing inserts at most one signed central
letter. There are at most $n+n^2$ letters thereafter. Moving and sorting
the central letters costs at most $(n+n^2)^2$ commutation replacements;
reducing powers costs at most $n+n^2$ more. Signed central inverses may
be moved using the same relation and cancelled freely. Thus the total
is bounded by $1024(L+1)^4$. Each presentation replacement costs at
most $32\epsilon$ in unitary 2-norm. Passing to conjugation operators
costs a factor two, since
$\|uau^*-vav^*\|_2\le2\|u-v\|_2$ for $\|a\|\le1$.

Choose symmetric free lifts of the group averages, of coefficient
$\ell^1$ norm at most one. The pair residual vanishes in the presented
group algebra; collecting every word to its normal form cancels its
coefficients. The three additional idempotence remainders
$Q_i^2-Q_i=P_i^2-P_i$ each have coefficient norm at most two. The
$1/36$ factor in the summed identity gives
\[
 \frac{3(2^{29}+2)}{36}\,2^{36}\,64\epsilon
 <2^{70}\epsilon.
\]
The positive square terms therefore give
\begin{equation}\label{fivemerge:eq:root-approx}
 \langle(L_H^2-L_H/60)a,a\rangle\ge-2^{70}\epsilon
 \qquad(\|a\|\le1).
\end{equation}
The finite checker verifies the 96 presentation-test matches, the
central fifth-power paths, the rational polynomial identity and norm
bound, and the integer comparisons. The collection argument above is
what transfers the exact finite-group identity to arbitrary tuples.

\subsubsection*{The imported analytic results}
The spectral correction and reconstruction are the results of
Wang--Zhi~\cite[Theorem 3.4, Corollary 4.3, Theorem 4.1, and
Lemma 5.1]{WangZhi2026} (version 1). Their inputs here are
\[
 (h_H,\lambda_H,C_H)=(1875,1/60,2^{70}),\qquad
 (h_K,\lambda_K,C_K)=(6,1/144,128).
\]
The second estimate is proved from the published integer theorem in
\eqref{fivemerge:eq:integer-approx}. The correction and reconstruction results
traced in Wang--Zhi include Alekseev--Liu--Thom~\cite[Theorems 4.3 and
5.2]{AlekseevLiuThom2026}; the normalization inputs include Thom~\cite[Proposition 3.1,
Lemma 4.1, and Theorem 4.2]{Thom2026}.

The tolerance comparison uses $\theta=2^{-140096}$ and the
Lemma 5.1 majorants
\[
 T_0=2^{1121200},\qquad M\le M_0=2^{32768T_0}.
\]
Both spectral inputs have $h\le2^{14}$ and
$2^{-9}\le\lambda/2\le1/4$. The error coefficient $C$ enters that construction only in its final
tolerance formula. Changing its bound to $C+1\le2^{71}$ therefore gives
\[
 \epsilon_{\rm in}\ge2^{-4M_0-71}\ge2^{-8M_0},\qquad
 8M_0=2^{3+2^{1121215}}<J_5.
\]
Hence $\delta_5<\epsilon_{\rm in}$ for both spectral inputs.

For normalization put
\[
 k=2^{-9},\quad r=2^{-17335},\quad
 \eta=8192r/k,\quad z=2^{-17500}.
\]
Apply Wang--Zhi Corollary 4.3 to the two spectral estimates. The
inequalities $3z\le r$ and $6\sqrt{2^{14}}z\le r$ give its Poincar\'e
and basic-root commutator bounds at the constants used in the published
normalization argument. A contraction in the first approximate
commutant has basic-root commutators at most $r$. Each $d_l,e_l$ is,
up to error $2\epsilon$, a four-letter basic-root commutator, so its
commutator costs at most $4r+4\epsilon$. Each product root is similarly
a four-letter $d/e$ commutator, giving $16r+20\epsilon$.
The compression rows add $2\epsilon$, hence the transported bound is
$16r+22\epsilon\le38r$. Averaging words have length at most eight,
so their energy square root is at most $(8/2)38r=152r$, below the
$168r$ allowance in the cited proof. This proves its forward near
inclusion with parameter $\eta$; Theorem 4.1 supplies the reverse
near inclusion with the same prescribed constants.

More explicitly, let $A$ be the root approximate-commutant algebra.
The Poincar\'e bound just used gives, for each first action $V=u_{ab}$
and each contraction $a\in A$,
\[
 \operatorname{dist}_2(V^*aV,A)
 \le152r/\sqrt{k}+r\le\eta.
\]
Theorem 4.1 reverses these inclusions to
$VAV^*\subseteq_{2^{-20}}A$. Set $c_0=u^*t$. Matching compression
targets give basic-root commutators of $c_0$ at most $2\epsilon$.
Its root energy square root is at most $8\epsilon$, and hence
\[
 \operatorname{dist}_2(c_0,A)
 \le8\epsilon/\sqrt{k}+r\le512r/k.
\]
Conditional expectations provide contraction approximants. The reverse
inclusion therefore places $tu^*=uc_0u^*$ within
$512r/k+2^{-20}$ of a contraction in $A$. For $h=z_2$ this yields
\[
 \|[h,tu^*]\|_2\le r+1024r/k+2^{-19}\le5\cdot2^{-20}.
\]
Multiplication by unitaries identifies that commutator norm with the
difference below:
\begin{equation}\label{fivemerge:eq:signal-close}
 \|u^*hu-t^*ht\|_2\le5\cdot2^{-20}
 \qquad(\epsilon<\delta_5/128).
\end{equation}
Only the fact that $h=z_2$ is a basic unitary root is used here.

\subsubsection{The integer adaptation is exact}\label{fivemerge:sec:integer-adaptation}
The integer part is independent of the characteristic-five root
calculation. In a split Steinberg relation let
$D=H-cb$ and $R=Ha-ab$. In the free $*$-algebra,
\[
 ab-cba=-R+Da.
\]
For an opposite-root relation, with $b$ denoting the appropriate
adjoint, let $D=H-ab$ and $R=Ha-bH$. Then
\[
 aba-bab=-Da+R+bD.
\]
Unitary multiplication preserves the tracial 2-norm. Together with the
six direct commutations, the 15 unsplit relations therefore have
maximum error at most $3\epsilon$. The independent checker matches
all 24 paid native rows and both formal identities.
Wang--Zhi~\cite[Lemma 3.3 and equation (3.12)]{WangZhi2026} (version 1) now gives
\begin{equation}\label{fivemerge:eq:integer-approx}
 \langle(L_K^2-L_K/144)a,a\rangle
 \ge-26(3\epsilon)\ge-128\epsilon\qquad(\|a\|\le1).
\end{equation}
Non-strict inequalities follow by increasing the error parameter and
letting it decrease to its endpoint. The resulting coefficient 78
already suffices; 128 is the retained conservative coefficient. This step
uses the cited spectral inequality and the displayed residual identities;
no separate integer spectral SOS file is required.

\subsubsection{Deduction of the native matrix gap}
\begin{proposition}\label{fivemerge:prop:native-gap}
Using the results and estimates above, every native 57 unitary matrix
tuple satisfies
\[
 \tr(F_{57}^0)+\mu_5Q_{57}\ge\tfrac34.
\]
\end{proposition}
\begin{proof}
If $Q_{57}<\delta_5^2/2^{14}$, every paid defect is strictly smaller
than $\delta_5/128$. The signal rows and~\eqref{fivemerge:eq:signal-close} imply
\[
 \|P-Q\|_2\le5\cdot2^{-20}+2\epsilon<\tfrac12.
\]
For unitaries $P,Q$,
$\tr(F_{57}^0)=1-\|P-Q\|_2^2\ge3/4$.
Otherwise $\mu_5Q_{57}\ge2^{16}$ by~\eqref{fivemerge:eq:parameters}, whereas
$\tr(F_{57}^0)\ge-3$. This also proves the conclusion.
\end{proof}
The proof uses no involutivity of the signal $h$; its unitarity and
basic-root role are what matter.

\Needspace{6\baselineskip}
\subsection{The compact and split sources}\label{fivemerge:sec:bridge}
\subsubsection{The compact 33-coordinate bridge}
The defining file \path{degree24/source33.json} has 33 coordinates:
21 symmetries and 12 genuine unitaries, with 78 paid rows. The
unitaries are
\[
\begin{gathered}
 x_0,y_0,z_{0,0},z_{1,0},z_{2,0},z_{3,0},p_{12,0},\quad
 H_{1,0},K_{xy,00},c_0,d_{0,0},\mathrm{square}_{x_0}.
\end{gathered}
\]
All other declared coordinates are symmetries. Its objective is
\[
 F_{33}^0=-1+VA^*VA+A^*VAV,\qquad A=z_{1,0}.
\]
The added square rows impose order five on $x_0$, and the paid
centralizers $[x_0,c_0]=[y_0,c_0]=1$ are essential for the signed
root calculus.

The finite arbitrary-tuple bridge asserts
\begin{equation}\label{fivemerge:eq:compact-bridge}
 Q_{57}\le C_{33}Q_{33},\qquad t(F_{57}^0)=t(F_{33}^0),
 \qquad C_{33}=2^{56}.
\end{equation}
The defining witness for this estimate is the typed rewrite certificate
\path{degree24/full_transfer_audit.json}.
A proof node with integer cost $k$ bounds its equality defect by
$k\sqrt{Q_{33}}$. Its rules are typed free cancellation, adjoint,
multiplication, transitivity, and cyclic relator reframing. The last
rule uses $\|L-R\|_2=\|LR^*-q\|_2$ for unitary words and invariance
under conjugation and inversion. All non-free leaves are paid rows and
all references point backwards. Summing squared row charges yields
\[
 1\,055\,531\,165\,250\,989<2^{56}.
\]
For transparency, the principal ledgers are as follows.

The signed root DAG has 2,005 nodes. It covers all 102 native root,
product and order-five rows with squared charge $2\,548\,771$.
The nine crossing charges are
\[
 2,25,49,45,57,47,95,103,159,
\]
and the six ordered product charges are
$201,279,176,296,254,321$.

For $x,y,z_0,z_1,z_2,z_3$, the fifth-power charges are
$3,8,13,23,18,23$.
These paths retain literal adjoints, rather than converting inverse
roots to roots.

The integer branch passes through auxiliary correctly typed tuples,
using finite-corner spectral-sign choices and finite word ledgers:
\[
 Q_{63}\le11Q_{33},\quad Q_{160}\le Q_{63},\quad
 Q_{157}\le22869994Q_{160},\quad Q_{48}\le985Q_{157}.
\]
Consequently $Q_{48}\le247796384990Q_{33}$. The 15 integer
relations have reflection paths of cost at most 11. Their
compression-compatible replacements have comparison costs
$0,0,2,8,7,3$, and hence each balanced integer relation has bound
\[
 (11\sqrt{247796384990}+48)\sqrt{Q_{33}}
 <2^{23}\sqrt{Q_{33}}.
\]
For a unitary $W$, the spectral-sign
step uses only
\[
 \|W-\fivesgn(\operatorname{Re}W)\|_2^2
 \le\tfrac12\|W-W^*\|_2^2.
\]

For common compression, use $u=R_{21}S_1$ and its six conjugates by
\[
 P_{21}=1,\ P_{12}=W_{12},\ P_{31}=W_{23},\
 P_{13}=W_{23}W_{12},\ P_{23}=W_{13},\ P_{32}=W_{12}W_{23}.
\]
The 42 goals have maximum cost 90. Set $t=Vu$ and $h=A$, exchanging
the parameter labels 1 and 2 to match the native convention, and set
$P=u^*hu$, $Q=t^*ht$. The objective comparison is exact because
\[
 P^*Q=u^*h^*VhVu,\qquad Q^*P=u^*Vh^*Vhu,
\]
and cyclicity gives $t(F_{57}^0)=t(F_{33}^0)$. There are 81 nonzero
root rows, 15 integer rows, 42 compression rows, and 32 free rows,
accounting for all 170 native rows.

\subsubsection{The expanded 65-coordinate source and trace scaling}
Rotate each $3+1$ row to $2+2$ by multiplying by a unitary inverse.
Of the 78 compact rows, 33 are split using 32 active shared helpers;
the unused allocated helper $H_{67}$ is discarded. The resulting
presentation has 65 coordinates, 36 genuine unitaries, and 111 paid
rows. For a split equality the squared triangle inequality gives a
factor at most two. Summing and allowing the unsplit rows the same
factor yields
\begin{equation}\label{fivemerge:eq:split}
 Q_{33}\le2Q_{65}.
\end{equation}
A shared-helper component with an odd adjoint cycle is a symmetry;
otherwise its helper is a genuine unitary.

Normalize the corner trace by $p$, apply Proposition~\ref{fivemerge:prop:native-gap},
and multiply back by $p$. Equations~\eqref{fivemerge:eq:compact-bridge}
and~\eqref{fivemerge:eq:split} give, for every correctly typed expanded tuple,
\begin{equation}\label{fivemerge:eq:expanded-gap}
 t(F_{33}^0)+C\mu_5Q_{65}\ge\tfrac34p,\qquad C=2^{57}.
\end{equation}
This establishes the unitary input used in the component-domain proof.

\subsection{The fixed component polynomial}\label{fivemerge:sec:components}
The eight recovery seeds are
\[
 F_{\rm graph},G_{\rm graph},J_{\rm hybrid},R_{21},\mathrm{Tau},
 V,W_{23},x_0.
\]
Only $x_0$ is a genuine-unitary seed. The 57 unique-target gates
$c=W_c$ are acyclic, have 46 words of length two and 11 of length
three, and have layer sizes
$3,2,5,5,6,7,11,3,7,4,3,1$. Each uses only earlier coordinates and
follows from a paid row. Their complete data are in
\path{degree24/expanded_source_and_recovery.json}.

Assign grade one to $V,G_{\rm graph},H_{21},H_{22},H_{23},H_{36}$,
and grade zero to the other coordinates. Every row is homogeneous
modulo two. Substitute $-iZ$ for an untwisted grade-one coordinate and
$iZ^*$ for its adjoint, then discard each common unit-modulus row phase.
The relative signs are real. The phase ledger has two negative row
multipliers and no negative gate multipliers. A graded symmetry is
skew-adjoint and has square $-q$ at the exact witness. The graded
objective is
\begin{equation}\label{fivemerge:eq:graded-objective}
 P_0=-q-VZ^*VZ-Z^*VZV,\qquad Z=z_{1,0},\quad V^*=-V.
\end{equation}
A genuine-unitary coordinate is written $Z=S_Z+K_Z$, with
$S_Z^*=S_Z$, $K_Z^*=-K_Z$; its adjoint is $S_Z-K_Z$.
Thus the 65 complete coordinates use $29+2\cdot36=101$ real components,
63 even and 38 odd.

For the definitions below, undo the phase substitution first; it is a
unit-modulus change of coordinates and preserves the relevant trace
identities. For a literal starred word $W$, define
\[
 L_W=\sum_{a\text{ occurring in }W}(Z_a^*Z_a-q),
\]
counting every occurrence, including repetitions and stars. The gate
polynomial and corrected cap are
\begin{align*}
 g_c&=Z_c^*Z_c+q-Z_c^*W_c-W_c^*Z_c+L_{W_c},\\
 \operatorname{cap}(u,v)&=2q-u^*v-v^*u+L_u+L_v,\\
 R_{\rm corr}&=\sum_{111\text{ paid rows }u-v}\operatorname{cap}(u,v).
\end{align*}
Set
\[
 k=2(T+1)^2,\qquad b_2=1+2T^2,\qquad b_3=2+T^2+3T^4.
\]
Initialize $a_v=1$ at every complete coordinate. Process the 57 gates
in reverse acyclic order: for target $c$ put $\beta_c=ka_c$, and add
$\beta_cb_{|W_c|}$ to $a_v$ for each input occurrence of $v$ in $W_c$.
Then the fixed budget is
\begin{equation}\label{fivemerge:eq:budget}
 E=\sum_{s\text{ seed}}a_s(Z_s^*Z_s-q)^2+\sum_c\beta_cg_c.
\end{equation}
Apply the phase substitution to this entire polynomial, not just the
objective. The collected component dictionary
\path{degree24/graded_bounded_support_polynomial.json} has nine
objective terms and 1,025 corrected-penalty terms. These data, the phase
ledger, the recovery order, and the occurrence counts are the finite
coefficient inputs to the next section.

\subsection{Component-domain coercivity}\label{fivemerge:sec:rounding}
Component bounds of $T$ give $\|S_Z+K_Z\|\le2T$, not $T$. The short-word estimates below retain the prescribed $T$, $k$, $b_2$, $b_3$, $K_{\rm corr}$ and $\lambda(T)$ by using the coefficient surplus.

In the proof, $P$ is the corrected penalty in \eqref{fivemerge:eq:source-polynomial}, $P_0$ is the ungraded objective, and $Q_{65}$ is a scalar squared-defect sum. Untwisting and retwisting preserve the trace and give the graded statement.

Let $q$ be the identity of a finite matrix corner with positive trace $t$,
$p=t(q)$, and write $\|x\|_r=t(|x|^r)^{1/r}$.
Untwist the phase convention by multiplying each grade-one coordinate by
$-i$. This leaves the norms, the unitarity defects, and the traces of the
corrected caps and gate terms unchanged, with the prescribed row phases.
The resulting coordinates are of the ungraded symmetry or unitary type.
For every complete coordinate $Z_a$, the stated component bounds imply
$\|Z_a\|\le 2T$; for a symmetry-type coordinate the stronger bound $T$
holds. Let $U_a$ be a unitary polar completion, or the sign of $Z_a$ at a
symmetry coordinate, and define
\[
 I_a=\|Z_a^*Z_a-q\|_2^2,\quad I=\sum_a I_a,\quad
 D_a=\|Z_a-U_a\|_2^2.
\]
Singular-value calculus gives $D_a\le I_a$ and
$\|Z_a-U_a\|_4^2\le\sqrt{I_a}$; it also gives
$I_a\le(2T+1)^2D_a\le4(T+1)^2D_a$.
These statements remain valid for starred letters, since the left and
right unitarity-defect squared traces coincide.

For two letters $A,B$ (possibly starred), write $D_A=A-U_A$ and similarly
for $B$. By trace cyclicity, Cauchy--Schwarz, and the above $4$-norm bound,
\[
 \|D_A B\|_2^2
 \le \|D_A\|_2^2+
 \|D_A^*D_A\|_2\,\|BB^*-q\|_2
 \le \tfrac32 I_A+\tfrac12 I_B.
\]
Thus the telescoping identity $AB-U_AU_B=D_AB+U_AD_B$ gives
\begin{equation}\label{fivemerge:eq:component-new-two}
 \|AB-U_AU_B\|_2^2\le 3(I_A+I_B).
\end{equation}
For three letters, set $D_{AB}=AB-U_AU_B$. Then
$\|D_{AB}\|_4^2\le8T^2\sqrt{I_A}+2\sqrt{I_B}$, and, applying the same
trace estimate to $D_{AB}C$,
\begin{align}
 \|ABC-U_AU_BU_C\|_2^2
 &\le (6+8T^2)I_A+8I_B+(4+8T^2)I_C \notag\\
 &\le (6+8T^2)(I_A+I_B+I_C). \label{fivemerge:eq:component-new-three}
\end{align}
Here $T\ge1$, and the factor order is retained.

For a word $W=T_1\cdots T_m$, put
$L_W=\sum_i(Z_{a_i}^*Z_{a_i}-q)$, counting starred occurrences by their
underlying coordinate. As an identity under $t$,
\[
 t(W^*W-q-L_W)=\sum_{i<j}t\bigl((T_i^*T_i-q)
 T_{i+1}\cdots T_{j-1}(T_jT_j^*-q)
 T_{j-1}^*\cdots T_{i+1}^*\bigr).
\]
For $m=2$ this has absolute value at most $\frac12\sum_i I_{a_i}$;
for $m=3$ it is at most $(\frac12+2T^2)\sum_i I_{a_i}$.
Combining these estimates with \eqref{fivemerge:eq:component-new-two} and
\eqref{fivemerge:eq:component-new-three}, the sum of the absolute trace remainder
and the squared rounding error is at most $b'_m\sum I_{a_i}$, with
\[
 b'_2=\tfrac72,\qquad b'_3=\tfrac{13}{2}+10T^2.
\]
The recovery data have 46 length-two and 11 length-three gates.
For every $T\ge4$, the prescribed recursion coefficients satisfy
\[
 4b'_2\le b_2:=1+2T^2,\qquad
 4b'_3\le b_3:=2+T^2+3T^4.
\]
For the gate $c=W_c$ of length $m$, its corrected polynomial satisfies
$t(g_c)=\|Z_c-W_c(Z)\|_2^2-t(W_c^*W_c-q-L_{W_c})$. Since $U_c$ is a
closest unitary and $W_c(U)$ is unitary,
\[
 \frac{I_c}{4}\le (T+1)^2D_c
 \le k\left(t(g_c)+b'_m\sum_{i\in W_c}I_i\right)
 \le k\,t(g_c)+\frac{k b_m}{4}\sum_{i\in W_c}I_i,
 \qquad k=2(T+1)^2.
\]
The gate trace may be negative; the displayed inequality nevertheless holds
with its sign and is the one to be summed.

Initialize $a_v=1$ at all 65 coordinates, and process the 57 gates in
reverse acyclic order as in the fixed polynomial: $\beta_c=k a_c$, and
add $\beta_c b_m$ to $a_v$ for every input occurrence. For the resulting
$E=\sum_s a_s(Z_s^*Z_s-q)^2+\sum_c\beta_cg_c$, summing the preceding
signed estimates and cancelling the acyclic charges gives
\begin{equation}\label{fivemerge:eq:component-E}
 t(E)\ge \frac14 I.
\end{equation}
Indeed the nonseed terms cancel to $\frac14\sum_c I_c$, while the seed
coefficient is $(3a_s+1)/4\ge1/4$.

Let $Q_{65}(Z)$ and $Q_{65}(U)$ be the scalar sums of the 111 squared
row defects. Every side has length at most two. The specified occurrence
counts for the 111 paid rows have at most 4 length-one and 21 length-two
occurrences per coordinate. The two-letter remainder and rounding bounds
above therefore imply, conservatively,
\[
 Q_{65}(Z)\le t(R_{\rm corr})+21I,\qquad
 Q_{65}(U)\le 2Q_{65}(Z)+4(4+3\cdot21)I
 =2Q_{65}(Z)+268I.
\]
With $K_{\rm corr}=30+84T^2$, put
$a=t(P)=t(R_{\rm corr}+K_{\rm corr}E)$. Equations above give
\begin{align*}
 a&\ge Q_{65}(Z)+(K_{\rm corr}/4-21)I\\
 &\ge \tfrac12Q_{65}(U)+(21T^2-590/4)I
 \ge \tfrac12Q_{65}(U)+64I
 \qquad(T\ge4).
\end{align*}
In particular $a\ge0$. The selected $T$ is at least four: the cheap
coordinate decoder $B_V=Q_0^2(XY-YX)$ has coefficient $\ell^1$ norm
$2Q_0^2$ and $Q_0>1$ in the specified integer recipe.

Assume now the native 57 finite-matrix gap and the certified compact bridge
and helper splitting; together they give
\[
 t(P_0(U))+C\mu_5Q_{65}(U)\ge 3p/4,\qquad C=2^{57}.
\]
The two objective words have length four. Using only the complete norm
bound $2T$, ordered telescoping and $\|x\|_1\le\sqrt p\|x\|_2$ give
\[
 |t(P_0(Z))-t(P_0(U))|\le64T^3\sqrt{pI}
 \le8T^3\sqrt{pa}.
\]
Consequently, with $\lambda(T)=2C+64T^6$,
\[
 \begin{aligned}
 t(\mathcal G)&\ge\tfrac34p-C\mu_5Q_{65}(U)
 -8T^3\sqrt{pa}+\mu_5\lambda(T)a\\
 &\ge \tfrac34p-8T^3\sqrt{pa}+64\mu_5T^6a
 \ge \left(\tfrac34-\frac1{4\mu_5}\right)p\ge\tfrac12p.
 \end{aligned}
\]
This proves the component-bound reduction with the fixed coefficients.

Together with~\eqref{fivemerge:eq:expanded-gap}, this proves the matrix conclusion
of Theorem~\ref{fivemerge:thm:conditional}.

\Needspace{6\baselineskip}
\subsection{The exact negative witness}\label{fivemerge:sec:witness}
\subsubsection{A faithful group construction}
Let
\[
 B=\fiveF_5[t_1,t_2,t_3],\quad
 \mathcal L=\fiveF_5[t_1^{\pm1},t_2^{\pm1},t_3^{\pm1}],\quad
 H=\fiveEL_3(B),\quad G=\fiveGL_3(\mathcal L)\rtimes\fiveGL_3(\fiveZ).
\]
An integer matrix $A$ acts on exponent columns by $t^n\mapsto t^{An}$;
thus $(M,A)(N,B)=(M\alpha_A(N),AB)$. Inverse transpose
$\vartheta(M,A)=(M^{-T},A)$ is an involutive automorphism, commutes with
these substitutions, and preserves $H$. On $\Omega=G*_H G$, use the
same $\vartheta$ on both factors and let $\sigma$ exchange the factors.
They commute, and
\[
 \Gamma=\Omega\rtimes(C_{2,\sigma}\times C_{2,\vartheta})
\]
contains $\Omega$ and both copies of $G$ injectively. The enlargement
to $\fiveGL$ is needed: the compact presentation uses exponent flips of
determinant $-1$ and Laurent diagonal matrices whose determinants are
nonconstant Laurent units. It imposes no quotient relation on the
specified amalgam.

Use the displayed elementary roots in the first factor. For integer
sign flips $s_j$ and $u_{ij}(k)=I+kE_{ij}$, set
\[
\begin{aligned}
 S_j&=(I,s_j),& R_{ij}&=(I,u_{ij}(1)s_j),&
 T_{ij}&=(I,u_{ij}(-1)s_j),\\
 W_{ij}&=R_{ij}R_{ji}R_{ij},&
 B_j&=(\operatorname{diag}(1,1,t_j),s_j),&
 C_j&=(\operatorname{diag}(t_j^{-1},1,1),s_j).
\end{aligned}
\]
The remaining named values are
\[
\begin{gathered}
 M_1=R_{32}R_{21}R_{32},\quad D_{\rm sign13}=S_1S_3,
 \quad J_{\rm hybrid}=C_2R_{13},\\
 \mathrm{Tau}=\vartheta,\quad
 F_{\rm graph}=(-P_{13})\vartheta,\quad
 G_{\rm graph}=(-P_{12})\vartheta,
\end{gathered}
\]
where $P_{ij}$ is the constant permutation matrix exchanging $i,j$.
For $w_{ij}=e_{ij}(1)e_{ji}(-1)e_{ij}(1)$,
\[
 -P_{13}=w_{12}^2w_{13},\qquad -P_{12}=w_{13}^2w_{12}.
\]
These explicit elementary decompositions place both signed permutations
in $H$. They square to one and are fixed by inverse transpose; thus the
graph elements are symmetries and commute with $\sigma$.

Put $D=W_{12}R_{12}S_2$ and $V=D\sigma D^{-1}$. Direct exponent
calculation gives
\[
\begin{aligned}
 D^{-1}x_0D&=x_0,&D^{-1}z_{2,0}D&=e_{31}(t_1),\\
 D^{-1}p_{12,0}D&=e_{31}(t_2),&D^{-1}z_{3,0}D&=e_{31}(t_3).
\end{aligned}
\]
The right sides lie explicitly in $H$. Hence $V$ commutes with the four
corresponding values. The pure exponent action $D$ commutes with
$\mathrm{Tau}$ and both graph elements, and those commute with $\sigma$.
This proves all seven rows involving $V$ and, in particular,
$[V,G_{\rm graph}]=1$. The other 71 compact rows are exact formal
Laurent identities in the embedded first-factor subgroup together with
$\vartheta$. The independent witness verification checks those identities,
the 12 genuine-unitary values, and the 20 remaining structural
involutions. No trace claim is inferred from the folded substitution
$V=1$.

\subsubsection{Normal form gives the negative trace}
Let $A=z_{1,0}$ and $h=D^{-1}AD=e_{31}(t_1^{-1}t_2)$. Both $h$ and
$h^{-1}$ lie outside $H$: their $(3,1)$ entries have a nonzero negative
$t_1$ exponent, whereas every matrix in $H$ has polynomial entries.
Conjugating the two objective words by $D^{-1}$ gives
\[
 \sigma h^{-1}\sigma h=h_2^{-1}h_1,\qquad
 h^{-1}\sigma h\sigma=h_1^{-1}h_2.
\]
Each is an alternating two-syllable word with both syllables outside
$H$. The amalgam normal-form theorem makes them nonidentity
\cite[Section 2.3, Theorems 2.5--2.6]{Miller2002}. The faithful left regular
representation of $\Gamma$ has canonical trace zero on every
nonidentity group element. Therefore all compact residuals vanish and
$\tau(F_{33}^0)=-1$ in the finite von Neumann algebra $A_0=L(\Gamma)$.
Here ``finite'' is tracial von Neumann finiteness, not finite dimension.

Assign shared helpers their representative group words. The
signed sharing graph proves all 111 expanded rows and the 57 recovery
gates, including the eight helper-symmetry types. For grade $g$, replace
a group value $u$ by $\mathsf J^g\otimes\lambda(u)$. The relative
powers of $\mathsf J$ give exactly the recorded phase signs. Splitting
a unitary into $(U+U^*)/2$ and $(U-U^*)/2$ produces components of their
declared types and norms at most one.

At this tuple all complete coordinates are unitary, all row and gate
defects vanish, and every occurrence correction vanishes. Thus
$E=R_{\rm corr}=P=0$. The two factors of $\mathsf J$ in each graded
objective word contribute $\mathsf J^2=-I$, cancelling the explicit
minus signs in~\eqref{fivemerge:eq:graded-objective}. Its normalized trace is
$-1$. The cheap coordinates are
$\mathsf J\otimes\lambda(V)$ and
$\mathsf J\otimes\lambda(G_{\rm graph})$ with commuting symmetry
factors. This proves the witness clause of Proposition~\ref{hyp:source24}, independently of its matrix clause.

\section{Quantitative proof of the operator-valued encoder}\label{app:encoderproof}\label{app:24data}
The estimates below prove Proposition~\ref{thm:encoder24} from the specified finite rational certificates. All spaces carry maximum operator norms; right inverses act between the resulting real Banach spaces. The isotropic identities keep the cheap channels fixed throughout the iteration. The defining arrays and word lists are identified below.
\subsection{Finite spectral weights and exact physical moments}
Fix $N=800$ and distinct rational numbers
\[
 \lambda_i=\frac{5N+1+i}{7N}\quad(0\le i<N),\qquad a_i=1/N.
\]
Let $E_i>0$ in $\mathcal B$, with $\sum_iE_i=1$. Form the isometry
$V:\mathcal B\to\mathcal B^N$ with entries $V_i=E_i^{1/2}$, and put
$D=\operatorname{diag}(\lambda_i)$, $Q_E=VV^*$ and
\[
 Y'_{ij}=\frac{E_i^{1/2}E_j^{1/2}}{(\lambda_i+\lambda_j)^2}.
\]
The exact Sylvester identity gives $e(D,Y')=Q_E$.
The integral formula
\[
 Y'=\int_0^\infty u e^{-uD}Q_E e^{-uD}\,du
\]
shows $Y'\ge0$ and $\enorm{Y'}\le1/(4\lambda_{\min}^2)<49/100$.
Also $5/7<\lambda_i\le6/7$.
Let $V_0=N^{-1/2}(1,\ldots,1)^\mathsf T$ and $q=V_0V_0^*$.
For weights sufficiently close to $a_i1$, the invertible operator
\[
 K_E=VV_0^*+(1-VV^*)(1-V_0V_0^*)
\]
has unitary polar factor $U$ satisfying $UV_0=V$.
Indeed $K_EV_0=V$ and $K_E^*K_EV_0=V_0$.
Set $X=U^*DU$, $Y=U^*Y'U$. Then $e(X,Y)=q$ exactly.

For a word $w=X^{b_0}YX^{b_1}\cdots YX^{b_r}$ its corner value,
identified with $\mathcal B$, is the finite polynomial
\begin{equation}\label{enc:eq:moment}
 \mathfrak m_w(E)=
 \sum_{i_0,\ldots,i_r}
 \frac{\prod_{j=0}^r\lambda_{i_j}^{b_j}}
      {\prod_{j=0}^{r-1}(\lambda_{i_j}+\lambda_{i_{j+1}})^2}
 E_{i_0}\cdots E_{i_r}.
\end{equation}
This follows directly by multiplying the entries of $Y'$ and placing
$V^*,V$ on the outside; no differentiation of a square root is required.
For $k=1,2$ in particular,
\begin{equation}\label{enc:eq:cheapmoment}
 \mathfrak m_{X^kY-YX^k}(E)=\sum_{i,j}C^k_{ij}E_iE_j,
 \qquad C^k_{ij}=\frac{\lambda_i^k-\lambda_j^k}
                         {(\lambda_i+\lambda_j)^2}.
\end{equation}
These two odd maps are exactly quadratic, not quadratic approximations.

An even raw polynomial $P_w$ is $w$ for a palindrome and $w+w^*$
otherwise. An odd raw polynomial is $P_w=w-w^*$.
The list of 110 even words is \path{degree24/selected_even_words.json};
the list of 36 ordinary odd words is \path{degree24/ordinary_odd_words.json}.
All even words have length at most eight, and all ordinary odd words
have length between three and seven.

\subsection{An exact common isotropic slice}
Represent $W=\{z\in\mathbb R^N:\sum z_i=0\}$ by its first $n=N-1$
coordinates, with last coordinate minus their sum. Let $A$ be the
$110\times n$ derivative of the even moments at $E_i=a_i1$ and let
$\Omega_k$ be the pullback to $W$ of the alternating form $C^k$.
The certificate prescribes rational columns
\[
 H\in\mathbb Q^{n\times240},\quad d\in\mathbb Q^n,
 \quad F=(f_1,f_2)\in\mathbb Q^{n\times2}
\]
satisfying, for $k=1,2$,
\begin{align}
 AH&=[I_{110}\ 0],&Ad&=AF=0,\nonumber\\
 H^\mathsf T\Omega_kH&=0,&H^\mathsf T\Omega_kd&=H^\mathsf T\Omega_kF=0,
 \label{enc:eq:isotropic}\\
 F^\mathsf T\Omega_kF&=0,&d^\mathsf T\Omega_kF&=e_k^\mathsf T.\nonumber
\end{align}
Here and below the final coordinate is restored by the zero-sum rule.

The exact definition is a finite sequential rational lift of the
residues modulo the prime $p_*=1000003$ specified in
\path{degree24/extended_h240_povm_N800_s2.npz}, with pivot records in
\path{degree24/extended_h240_povm_N800_s2.json}. At each $H_j$, impose its
specified $A$ value and all earlier $H$ orthogonality equations.
Then solve for $d$ with zero $A$ value and every $H$ orthogonality;
then solve for $f_1,f_2$ with their two pairing values and all earlier
orthogonality equations. Every free coordinate is fixed to its saved
residue in $\{0,\ldots,p_*-1\}$; solve the recorded pivot square system
over $\mathbb Q$. All 243 stages and all 594 accumulated equations
are retained. Each recorded pivot minor is nonzero modulo $p_*$; the final pivot minor is $3521$.
Induction shows that each rational solve exists, has $p_*$-integral entries, and reduces to the saved array, proving \eqref{enc:eq:isotropic} over $\mathbb Q$.

Put $D_2=\operatorname{diag}(1,-1)$ and
$E_2=\left(\begin{smallmatrix}0&1\\1&0\end{smallmatrix}\right)$.
For self-adjoint controls $C_1,\ldots,C_{240}\in\mathcal B$, use
\begin{equation}\label{enc:eq:weights}
 E_i=a_i1+\varepsilon d_iD_2/2
       +\varepsilon f_{i1}(E_2\otimes S_1)
       +\varepsilon f_{i2}(E_2\otimes S_2)+\sum_lH_{il}C_l.
\end{equation}
The weights sum to one. Equations~\eqref{enc:eq:isotropic} cancel every
quadratic term in \eqref{enc:eq:cheapmoment} except anchor-partner terms.
Since $[D_2/2,E_2]=\mathsf J$,
\begin{equation}\label{enc:eq:exactcheap}
 \mathfrak m_{XY-YX}(E)=\varepsilon^2\mathsf J\otimes S_1,
 \quad
 \mathfrak m_{X^2Y-YX^2}(E)=\varepsilon^2\mathsf J\otimes S_2.
\end{equation}
These equalities hold for every control tuple, including noncommuting
controls. Thus all subsequent ordinary-output corrections preserve the
two cheap outputs exactly.

\subsection{An exactly even-centered reference}
Let $c_w\in\mathbb Q$ be the even scalar moment at $E_i=a_i1$.
Set $C_l=0$ for $l>110$. The derivative of the first 110 even outputs
with respect to the other controls at the origin is the identity.
Their derivative in $\varepsilon$ is zero. A uniform contraction
argument, quantified below, gives a unique analytic reference
$C^\flat(\varepsilon,S_1,S_2)=O(\varepsilon^2)$ at which every even
moment equals $c_w1$ exactly.

Let $p_s$ be the joint spectral projections of $S_1,S_2$,
$s=(s_1,s_2)\in\{\pm1\}^2$. The iteration defining $C^\flat$ starts
at zero and stays in the real algebra generated by $D_2,E_2S_1,E_2S_2$.
In each sector its weights are real symmetric $2\times2$ matrices.
Their linear term is
\[
 T_i(s)=d_iD_2/2+(f_{i1}s_1+f_{i2}s_2)E_2.
\]
A reversal-even quadratic expression in these matrices is scalar,
because $D_2^2=E_2^2=I$ and $D_2E_2+E_2D_2=0$.
Thus the reference coefficient $C_2(s)$ at order $\varepsilon^2$ is
scalar and depends only on $s_1s_2$. In particular $C_2(s)=C_2(-s)$.
Every odd quadratic moment is a sum of commutators, so its derivative
at this scalar correction vanishes within a diagonal sector and between
opposite sectors. Hence the correction does not change the sensitive
order-two ordinary-odd jets. Between single-flip sectors the rank proof
uses only order one, which an order-two correction cannot affect.

\subsection{Leading minors and all Peirce components}
Differentiate the 110 even and 36 ordinary odd outputs with respect to
all 240 controls at the reference. In a diagonal joint sector, decompose
self-adjoint $2\times2$ inputs into three real-symmetric matrix components
with self-adjoint algebra coefficients and one $\mathsf J$ component
with a skew-adjoint algebra coefficient. The two corresponding real
scalar derivative blocks use these leading orders:
\begin{center}
\begin{tabular}{lll}
\toprule
Input block & Even rows & Ordinary odd rows\\
\midrule
Real symmetric & order 0 & order 1\\
$\mathsf J$ component & order 0 & traceless order 1, scalar order 2\\
\bottomrule
\end{tabular}
\end{center}
For a single-flip pair of sectors use all four matrix-unit components;
even rows have order zero and all ordinary odd rows have order one.
Independent reconstruction gives the following full-row minors:
\begin{center}
\begin{tabular}{llrr}
\toprule
Sector & Input type & Rows & Rank\\
\midrule
$(+,+)$ & real symmetric &366&366\\
$(+,+)$ & $\mathsf J$ component &218&218\\
$(+,-)$ & real symmetric &366&366\\
$(+,-)$ & $\mathsf J$ component &218&218\\
$(+,+)\to(+,-)$ & all matrix units &584&584\\
$(+,+)\to(-,+)$ & all matrix units &584&584\\
\bottomrule
\end{tabular}
\end{center}
Every row is retained. The exact modular matrices and their pivot
columns are specified in \path{degree24/anchor_leading_jet_matrices.npz}
and \path{degree24/anchor_leading_jet_ranks.json}.
The rational lift described above makes these nonzero modular minors nonzero
rational leading minors.

These six cases cover all joint Peirce blocks. Simultaneous sign reversal
is conjugation by $D_2$. Opposite sectors $s,-s$ reduce to a diagonal
sector by multiplying one matrix side by $D_2$, since
$E_i^\flat(-s)=D_2E_i^\flat(s)D_2$.
Reversing a sector pair is handled by adjoints. The two remaining
single-flip classes are exactly the two 584-row cases. Complexifying
the diagonal real-symmetric and $\mathsf J$ blocks gives all matrix
components for opposite-sector inputs; pairing adjoint blocks then
restores self-adjoint controls. Zero joint projections remove
blocks. Applying the finite scalar right inverses entrywise in each
Peirce block works in every $C^*$-algebra; all projection multiplications
are contractions. The finite sum over at most 16 blocks therefore has
a dimension-independent bound.

After dividing the specified output rows by their powers
$1,\varepsilon,\varepsilon^2$, each pivot matrix has an invertible
constant term. The quantitative estimates below give a right inverse
$R_\varepsilon$ of the full control derivative satisfying
\begin{equation}\label{enc:eq:rightinverse}
 \enorm{R_\varepsilon}\le K_{\rm inv}\varepsilon^{-2}.
\end{equation}

The six centered leading-jet blocks have the following nonzero minor residues, obtained by reconstructing their entries and performing exact finite-field elimination:
\begin{center}\small
\begin{tabular}{@{}lrrr@{}}
\toprule Block&Rows (rank)&Columns&Minor modulo \(1000003\)\\\midrule
Skew, sector \((1,1)\)&218&240&467093\\
Symmetric, sector \((1,1)\)&366&720&424753\\
Skew, sector \((1,-1)\)&218&240&740276\\
Symmetric, sector \((1,-1)\)&366&720&700354\\
Mixed, sector \((1,-1)\)&584&960&884863\\
Mixed, sector \((-1,1)\)&584&960&884863\\\bottomrule
\end{tabular}
\end{center}
The finite Taylor-center verification also reconstructs all 36 ordinary odd rows and the four sector arrays through orders two, three and four.

\subsection{Low-degree rational odd centering}
Every ordinary odd reference value in a joint sector is
$\mathsf J\beta_j(\varepsilon,s)$ with a real scalar $\beta_j$.
Conjugation by $D_2$ reverses $s$ and reverses $\mathsf J$, so
$\beta_j(\varepsilon,-s)=-\beta_j(\varepsilon,s)$.
Every odd function on the four sign pairs has the form
\[
 \beta_j(\varepsilon,s)=a_j(\varepsilon)s_1+b_j(\varepsilon)s_2.
\]
Conjugation by $\mathsf J$ sends $\varepsilon$ to $-\varepsilon$
and fixes $\mathsf J$. Uniqueness of even centering implies that
$a_j,b_j$ are even functions of $\varepsilon$. The odd map has no
constant or linear term at the scalar reference. Hence
\[
 a_j=a_{j2}\varepsilon^2+a_{j4}\varepsilon^4+O(\varepsilon^6),
 \qquad
 b_j=b_{j2}\varepsilon^2+b_{j4}\varepsilon^4+O(\varepsilon^6).
\]
All four coefficients are rational. They are defined by finite Taylor
recursion: set the degree-two and degree-three even moment coefficients
to zero using $AH_{1:110}=I$, then evaluate the odd moments through
degree four. No $C_4$ is needed because the odd derivative at the scalar
reference vanishes. The saved modular recursion is independently verified;
the characteristic-zero coefficients are defined by the same rational
recursion on the lifted arrays.

Define ordinary centered odd polynomials
\begin{equation}\label{enc:eq:oddcenter}
 \begin{split}
 \widehat P_j=P_j
 &-(a_{j2}+a_{j4}\varepsilon^2)(XY-YX)\\
 &-(b_{j2}+b_{j4}\varepsilon^2)(X^2Y-YX^2).
 \end{split}
\end{equation}
Their reference residual is $O(\varepsilon^6)$, and their derivative
on the control slice is unchanged by \eqref{enc:eq:exactcheap}.
All ordinary odd raw degrees are at least three, so no raw degree
increases. Even centering is $\widehat P_w=P_w-c_w$.

\subsection{Explicit integer bounds and the fixed parameter}
This deliberately large recipe removes every unspecified smallness choice.
The following integers depend only on the saved arrays and $N,p_*$.
Set $n=N-1$, $v=243$, and
\[
 D_0=((12N)!)^{25},\quad R_0=p_*,\quad
 B_j=36^nD_0R_j^{n+1},\quad
 R_{j+1}=n!(1+np_*)^{B_j^n}\quad(0\le j<v).
\]
For a word with $r$ occurrences of $Y$ and $x$ of $X$, $r+x\le8$,
canceling powers of $7N$ bounds the initial denominator exponent by
$3r+\max(x-2r,0)\le24$. Thus $D_0$ clears every initial scalar
constraint coefficient.
Write $R=R_v$ and $D_*=D_0R!$. The rational-lift Cramer induction
bounds every numerator and denominator of a first-$n$-coordinate entry
by $R$, and every extended entry in magnitude by $nR$.
Indeed each new constraint coefficient has a common denominator bounded
by the products of at most $n$ previous denominators and $D_0$;
clearing them gives the bound $B_j$. Free coordinates are at most $p_*$.
Cramer's rule on a square subsystem of size at most $n$ bounds both
its determinant and replaced-column determinants. Clearing at most $n$ resulting denominators is bounded by
\[
 n!\bigl[(1+np_*)B_j^n\bigr]^n.
\]
Our fixed recurrence is larger than this Cramer bound. Indeed put $a=1+np_*\ge2$. Since $B_j\ge36^n\ge2n^2$, $n=799\ge3$, and $B_j\le a^{B_j}$, we have
\[
 a^nB_j^{n^2}\le a^{n+n^2B_j}\le a^{B_j^n}.
\]
Also $36^n\ge36n$, so the stated $B_j$ covers the coefficient-clearing bound. This proves the recurrence inductively, including $p_*$-integrality from the nonzero modular minors, while retaining the specified decoder scales.

Set
\[
 M_0=10^6N^{12}(R+1)^4,\qquad
 B_{\rm piv}=584!(M_0+1)^{584}D_*^{40\cdot584}.
\]
Every leading derivative entry has denominator dividing $D_*^{40}$
and magnitude at most $M_0$. To see this, a raw word contributes at
most $2N^9$ summands, has at most nine weight positions, and its
order-two derivative uses at most 28 anchor choices. The $H$ and
anchor entries have magnitudes at most $NR$ and $3NR$.
There are at most 16 spectral denominator factors and at most 11
weight, anchor, and $H$ factors. The row sums and differences cost
at most two. These counts are covered by the stated bounds.
A nonzero rational $r\times r$ minor, $r\le584$, thus has determinant
magnitude at least $D_*^{-40r}$. Bounding its cofactors by the Leibniz
formula gives inverse infinity norm at most $B_{\rm piv}$.

Define
\begin{align*}
 K_{\rm ref}&=110NR,&K_{\rm all}&=240NR,&G_0&=3NR,\\
 L_2&=144N^9,& b&=2+L_2G_0^2,\\
 Z_0&=8\max(N,L_2K_{\rm ref})(G_0+K_{\rm ref}b+1),&\rho&=1/Z_0,\\
 J_{\max}&=18N^9K_{\rm all},&K_{\rm inv}&=2^{24}\cdot240B_{\rm piv},\\
 H_2&=L_2K_{\rm all}^2,&C_{\rm rem}&=8N^9Z_0^6.
\end{align*}
On the weight unit ball the Hessian norm of every raw moment is at
most $144N^9$: at most $2N^9$ terms have degree at most nine, and
$9\cdot8$ second derivative choices suffice. The even-reference
contraction works holomorphically for complex $|\varepsilon|\le\rho$,
with $\enorm{C^\flat}\le b|\varepsilon|^2$ and weight norms below one.
In detail the first derivative is the identity, the first anchor
term vanishes, and the Hessian estimate bounds the derivative error
by $L_2K_{\rm ref}(G_0|\varepsilon|+K_{\rm ref}b|\varepsilon|^2)$.
The definition of $Z_0$ makes this less than $1/2$, keeps the fixed-point
ball invariant, and makes all real weights positive.

The full control derivative is bounded on this disk by $J_{\max}$.
Cauchy bounds show that each rescaled leading matrix differs from its
constant matrix in infinity norm, for $|\varepsilon|\le\rho/2$, by at most
$1600J_{\max}Z_0^3|\varepsilon|$. The row-summed Cauchy factors are $4\cdot218=872$,
$2\cdot366=732$ and $2\cdot584=1168$; the uniform factor 1600
therefore covers the largest sensitive rescaled block and its geometric
series remainder. Thus its pivot inverse has norm at most $2B_{\rm piv}$ if
\[
 |\varepsilon|\le(3200B_{\rm piv}J_{\max}Z_0^3)^{-1}.
\]
The generous $K_{\rm inv}$ covers all Peirce decompositions, matrix
component norms, and assembly factors, proving \eqref{enc:eq:rightinverse}.
The odd scalar centers have magnitude at most $2N^9$ on the reference
disk. Their even-in-$\varepsilon$ Taylor remainder after degree four
is at most $C_{\rm rem}\varepsilon^6$ for
$|\varepsilon|\le\rho/2$, including both $a$ and $b$ terms.

Now fix the integer and rational parameter
\begin{equation}\label{enc:eq:Qchoice}
 \begin{split}
 Q_0={}&1+2Z_0+3200B_{\rm piv}J_{\max}Z_0^3
           +4NK_{\rm all}(K_{\rm inv}+1)\\
       &+32C_{\rm rem}H_2K_{\rm inv}^2,
 \qquad \varepsilon=Q_0^{-1}.
 \end{split}
\end{equation}
The reference residual after odd centering is at most
$\varepsilon^4/(32H_2K_{\rm inv}^2)$.
Work in output coordinates $u$ with control correction $R_\varepsilon u$.
The derivative of this map at zero is the identity; its Hessian is at
most $H_2K_{\rm inv}^2\varepsilon^{-4}$.
On the ball of radius
\[
 r=\frac{\varepsilon^4}{8H_2K_{\rm inv}^2},
\]
its nonlinear remainder has Lipschitz constant at most $1/8$.
The reference residual is at most $r/4$; a target of norm at most $r/4$
therefore gives a contraction self-map of this ball.
The reference weights on the correction ball differ from $N^{-1}I$ by at most
\[
 G_0\varepsilon+K_{\rm ref}b\varepsilon^2
 \le(G_0+K_{\rm ref}b)\varepsilon<\frac1{16N},
\]
because $\varepsilon<(2Z_0)^{-1}$ and $Z_0\ge8N(G_0+K_{\rm ref}b+1)$. Moreover $r\le\varepsilon^3$, since $\varepsilon<1$ and $H_2,K_{\rm inv}\ge1$. The correction changes a weight by at most
\[
 K_{\rm all}\|R_\varepsilon\|r
 \le K_{\rm all}K_{\rm inv}\varepsilon
 <\frac1{4N},
\]
using $Q_0>4NK_{\rm all}(K_{\rm inv}+1)$. Thus
$\max_i\|E_i-N^{-1}I\|<5/(16N)$ and $E_i>11/(16N)I$.
This also justifies the polar-unitary step explicitly: the scalar square-root comparison gives $\|V-V_0\|<5/16$, and hence
\[
 \|K_E-I\|\le\|V-V_0\|+\|VV^*-V_0V_0^*\|<15/16<1.
\]
So $K_E$ is invertible and its polar factor is a unitary, with no finite-dimensional restriction on the auxiliary algebra.

A common denominator for all rational Taylor centers is $D_*^{4096}$.
One direct exponent induction bounds anchor coefficients by denominator
$D_*^2$, then $C_2$ by $D_*^{34}$ and $HC_2$ by $D_*^{35}$,
then $C_3$ by $D_*^{331}$ and $HC_3$ by $D_*^{332}$.
Odd evaluation through order four has denominator dividing $D_*^{3004}$;
extracting $a,b$ by half-sums adds at most one power of $D_*$.
This overcounts every spectral factor and all nine weight factors,
and also covers the even scalar centers.
Set
\begin{equation}\label{enc:eq:scales}
 R_{\rm dec}=32H_2K_{\rm inv}^2D_*^{4096},\quad
 A_{\rm ord}=R_{\rm dec}Q_0^4,\quad A_{\rm cheap}=Q_0^2.
\end{equation}
Then $1/A_{\rm ord}\le r/4$ and the decoders
\[
 B_w=A_{\rm ord}(P_w-c_w),\quad
 B_j=A_{\rm ord}\widehat P_j,\quad
 B_V=Q_0^2(XY-YX),\quad
 B_{G_{\rm graph}}=Q_0^2(X^2Y-YX^2)
\]
have integer coefficients. The $\varepsilon^2$ coefficients in
\eqref{enc:eq:oddcenter} are cleared by $Q_0^4$ and $D_*^{4096}$.
Equations~\eqref{enc:eq:exactcheap} give their exact cheap outputs, and
the preceding contraction gives every prescribed ordinary output.
This proves Proposition~\ref{thm:encoder24}.

\Needspace{8\baselineskip}
\subsection{Complete active decoder-channel assignment}
The 101 rows below give each source coordinate, component type ($S$ self-adjoint, $K$ skew-adjoint), raw degree, weighted degree including selector cost three, and raw word.
\small
\begin{longtable}{@{}lllll@{}}
\toprule Coordinate&Component&Raw degree&Weight&Raw word\\\midrule
\endfirsthead
\toprule Coordinate&Component&Raw degree&Weight&Raw word\\\midrule
\endhead
x0&S&3&6&\texttt{XXX}\\
x0&K&3&6&\texttt{XYY}\\
y0&S&4&7&\texttt{XXXX}\\
y0&K&4&7&\texttt{XXYY}\\
z0 0&S&5&8&\texttt{XXXXX}\\
z0 0&K&5&8&\texttt{XXXYY}\\
z1 0&S&4&7&\texttt{XXXY}\\
z1 0&K&4&7&\texttt{XYYY}\\
z2 0&S&5&8&\texttt{XXXXY}\\
z2 0&K&5&8&\texttt{XXYYY}\\
z3 0&S&6&9&\texttt{XXXXXX}\\
z3 0&K&6&9&\texttt{XXXXYY}\\
p12 0&S&6&9&\texttt{XXXXXY}\\
p12 0&K&6&9&\texttt{XXXYYY}\\
S1&S&1&4&\texttt{X}\\
S3&S&5&8&\texttt{XXXYY}\\
R13&S&4&7&\texttt{XXYY}\\
R21&S&2&5&\texttt{XX}\\
R31&S&4&7&\texttt{XYXY}\\
R32&S&4&7&\texttt{XYYX}\\
T31&S&5&8&\texttt{XXYXY}\\
W12&S&3&6&\texttt{XXY}\\
W13&S&3&6&\texttt{XYY}\\
W23&S&2&5&\texttt{XY}\\
Kxy 00&S&6&9&\texttt{XXXXYY}\\
Kxy 00&K&6&9&\texttt{XXYXYY}\\
H1 0&S&6&9&\texttt{XXXYXY}\\
H1 0&K&6&9&\texttt{XXYYYY}\\
M1&S&4&7&\texttt{XYYY}\\
V&K&2&5&\texttt{XY}\\
c0&S&5&8&\texttt{XXYYX}\\
c0&K&5&8&\texttt{XYXYY}\\
d0 0&S&5&8&\texttt{XXYYY}\\
d0 0&K&5&8&\texttt{XYYYY}\\
B2&S&5&8&\texttt{XYXYY}\\
C1&S&4&7&\texttt{YXYY}\\
C2&S&1&4&\texttt{Y}\\
Fgraph&S&3&6&\texttt{YXY}\\
Ggraph&K&3&6&\texttt{XXY}\\
Tau&S&3&6&\texttt{YYY}\\
Dsign13&S&5&8&\texttt{XYYXY}\\
Jhybrid&S&2&5&\texttt{YY}\\
T21&S&4&7&\texttt{YYYY}\\
square x0&S&6&9&\texttt{XXXYYX}\\
square x0&K&6&9&\texttt{XYXYYY}\\
H1&S&6&9&\texttt{XXXYYY}\\
H1&K&6&9&\texttt{XYYXYY}\\
H3&S&6&9&\texttt{XXYXYY}\\
H3&K&6&9&\texttt{XYYYXY}\\
H4&S&6&9&\texttt{XXYYXY}\\
H4&K&6&9&\texttt{XYYYYY}\\
H6&S&7&10&\texttt{XXXXXXX}\\
H6&K&7&10&\texttt{XXXXXYY}\\
H7&S&7&10&\texttt{XXXXXXY}\\
H7&K&7&10&\texttt{XXXXYYY}\\
H9&S&7&10&\texttt{XXXXXYY}\\
H12&S&7&10&\texttt{XXXXYXY}\\
H12&K&7&10&\texttt{XXXYXYY}\\
H14&S&5&8&\texttt{XYYYX}\\
H14&K&7&10&\texttt{XXXYYYY}\\
H15&S&6&9&\texttt{XXYYYX}\\
H15&K&7&10&\texttt{XXYXYYY}\\
H21&S&6&9&\texttt{XXYYYY}\\
H21&K&7&10&\texttt{XXYYXYY}\\
H22&S&7&10&\texttt{XXXXYYX}\\
H22&K&7&10&\texttt{XXYYYXY}\\
H23&S&6&9&\texttt{XYXYYY}\\
H23&K&7&10&\texttt{XXYYYYX}\\
H24&S&6&9&\texttt{XYYXYY}\\
H24&K&7&10&\texttt{XXYYYYY}\\
H26&S&6&9&\texttt{XYYYXY}\\
H26&K&6&9&\texttt{YXYYYY}\\
H27&S&6&9&\texttt{XYYYYX}\\
H29&S&7&10&\texttt{XXXXYYY}\\
H29&K&7&10&\texttt{XYXYYYX}\\
H30&S&5&8&\texttt{XYYYY}\\
H36&S&7&10&\texttt{XXXYXYY}\\
H36&K&7&10&\texttt{XYXYYYY}\\
H39&S&5&8&\texttt{YXYYY}\\
H39&K&7&10&\texttt{XYYXYYY}\\
H40&S&7&10&\texttt{XXXYYXY}\\
H40&K&7&10&\texttt{XYYYXYY}\\
H45&S&6&9&\texttt{XYYYYY}\\
H45&K&7&10&\texttt{XYYYYXY}\\
H47&S&6&9&\texttt{YXYXYY}\\
H55&S&6&9&\texttt{YXYYXY}\\
H58&S&5&8&\texttt{YYXYY}\\
H58&K&5&8&\texttt{YXYYY}\\
H60&S&6&9&\texttt{YXYYYY}\\
H61&S&7&10&\texttt{XXXYYYX}\\
H61&K&7&10&\texttt{XYYYYYY}\\
H63&S&6&9&\texttt{YYXYYY}\\
H63&K&7&10&\texttt{YXYXYYY}\\
H68&S&8&11&\texttt{XXXXXXYY}\\
H72&S&6&9&\texttt{YYYYYY}\\
H77&S&8&11&\texttt{XXXXXYXY}\\
H77&K&7&10&\texttt{YXYYYYY}\\
H86&S&7&10&\texttt{XXXYYYY}\\
H86&K&7&10&\texttt{YYXYYYY}\\
H89&S&5&8&\texttt{YYYYY}\\
H89&K&6&9&\texttt{YYXYYY}\\
\bottomrule
\end{longtable}
\normalsize

\enlargethispage{2\baselineskip}


\begin{thebibliography}{99}
\bibitem{Connes1976}
A.~Connes, \emph{Classification of injective factors. Cases \(II_1\), \(II_\infty\), \(III_\lambda\), \(\lambda\ne1\)}, Ann. of Math. (2) \textbf{104} (1976), no.~1, 73--115.

\bibitem{JiEtAl2020}
Z.~Ji, A.~Natarajan, T.~Vidick, J.~Wright, and H.~Yuen, \emph{MIP* = RE}, arXiv:2001.04383 (2020).

\bibitem{Helton2002}
J.~W. Helton, \emph{``Positive'' noncommutative polynomials are sums of squares}, Ann. of Math. (2) \textbf{156} (2002), no.~2, 675--694.

\bibitem{BurgdorfKlep2012}
S.~Burgdorf and I.~Klep, \emph{The truncated tracial moment problem}, J. Operator Theory \textbf{68} (2012), no.~1, 141--163.

\bibitem{KlepSchweighofer2008}
I.~Klep and M.~Schweighofer, \emph{Connes' embedding conjecture and sums of Hermitian squares}, Adv. Math. \textbf{217} (2008), no.~4, 1816--1837.

\bibitem{BurgdorfDykemaKlepSchweighofer2014}
S.~Burgdorf, K.~Dykema, I.~Klep, and M.~Schweighofer, \emph{Addendum to ``Connes' embedding conjecture and sums of Hermitian squares''}, Adv. Math. \textbf{252} (2014), 805--811.

\bibitem[Ershov and Jaikin-Zapirain(2010)]{ErshovJaikinZapirain2010}
M.~Ershov and A.~Jaikin-Zapirain, \emph{Property (T) for noncommutative universal lattices}, Invent. Math. \textbf{179} (2010), no.~2, 303--347.

\bibitem[Ozawa(2016)]{Ozawa2016}
N.~Ozawa, \emph{Noncommutative real algebraic geometry of Kazhdan's property (T)}, J. Inst. Math. Jussieu \textbf{15} (2016), no.~1, 85--90.

\bibitem[Miller(2002)]{Miller2002}
C.~F. Miller III, \emph{Combinatorial Group Theory}, lecture notes, March 5, 2002, Section 2.3, Theorems 2.5--2.6. \url{https://www.macs.hw.ac.uk/~lc45/Teaching/kggt/miller.pdf}.

\bibitem[Wang and Zhi(2026b)]{WangZhi2026}
J.~Wang and L.~Zhi, \emph{An explicit polynomial counterexample to Connes' embedding conjecture}, arXiv:2610.01536 (2026).

\bibitem[Kun and Thom(2026)]{KunThom2026}
G.~Kun and A.~Thom, \emph{Nonsofic wreath products of residually finite groups}, arXiv:2608.06222 (2026).

\bibitem[Thom(2026)]{Thom2026}
A.~Thom, \emph{A conditional construction of a nonhyperlinear group and the centralizer problem}, arXiv:2609.27795 (2026).

\bibitem[Alekseev et~al.(2026)]{AlekseevLiuThom2026}
V.~Alekseev, J.~Liu, and A.~Thom, \emph{Spectral gap, internality of relative commutants, and non-hyperlinear groups}, arXiv:2609.37761 (2026).

\bibitem[Bercovici et~al.(2010)]{BercoviciEtAl2009}
H.~Bercovici, B.~Collins, K.~Dykema, W.~S. Li, and D.~Timotin, \emph{Intersections of Schubert varieties and eigenvalue inequalities in an arbitrary finite factor}, J. Funct. Anal. \textbf{258} (2010), no.~5, 1579--1627.

\end{thebibliography}
\end{document}